\documentstyle[12pt]{amsart}

\numberwithin{equation}{section}
\newtheorem{theorem}[equation]{Theorem}
\newtheorem{lemma}[equation]{Lemma}

\newtheorem{prop}[equation]{Proposition}
\newtheorem{cor}[equation]{Corollary}
\theoremstyle{definition}
\newtheorem{defn}[equation]{Definition}
\newtheorem{example}[equation]{Example}
\theoremstyle{remark}
\newtheorem{remark}[equation]{Remark}

\newcommand{\QQbar}{\overline{{\mathbb Q}}}
\newcommand{\kvbar}{{\overline{k}}_v}

\renewcommand{\O}{{\mathcal O}}

\newcommand{\hhat}{\hat{h}}
\newcommand{\Hhat}{\hat{H}}

\newcommand{\Haar}{{\rm Haar}}
\newcommand{\into}{\hookrightarrow}     

\def\Berk{{\mathop{\rm Berk}}}
\def\BDV{{\mathop{\rm BDV}}}

\def\Gal{{\mathop{\rm Gal}}}

\def\supp{{\mathop{\rm supp}}}
\def\div{{\mathop{\rm div}}}

\def\det{{\mathop{\rm det}}}
\def\Det{{\mathop{\rm Det}}}
\def\dim{{\mathop{\rm dim}}}

\def\vol{{\mathop{\rm vol}}}
\def\deg{{\mathop{\rm deg}}}

\def\Spec{{\mathop{\rm Spec}}}

\def\Res{{\mathop{\rm Res}}}

\def\CH{{\mathop{\it CH}}}
\def\SL{{\mathop{\rm SL}}}
\def\SU{{\mathop{\rm SU}}}

\def\<{\langle }
\def\>{\rangle }
\def\({(\!(}
\def\){)\!)}

\def\AA{{\Bbb A}}
\def\BB{{\Bbb B}}
\def\CC{{\Bbb C}}

\def\QQ{{\Bbb Q}}
\def\EE{{\Bbb E}}
\def\NN{{\Bbb N}}
\def\RR{{\Bbb R}}

\def\PP{{\Bbb P}}

\def\v1{{\vec{1}}}

\def\BF1{{\mathbf 1}}

\def\cB{{\mathcal B}}
\def\cC{{\mathcal C}}

\def\cF{{\mathcal F}}
\def\cG{{\mathcal G}}

\def\cO{{\mathcal O}}

\def\hcO{{\widehat{\mathcal O}}}

\def\kbar{{\overline{k}}}

\def\Qbar{{\overline{\QQ}}}

\def\Sg{S_{\gamma}}

\def\tF{{\tilde{F}}}

\def\tk{{\tilde{k}}}
\def\tz{{\tilde{z}}}
\def\tw{{\tilde{w}}}

\def\blacksquare{\Box}
\newenvironment{proof}{\noindent{\textbf{Proof: }}}{\hfill $\blacksquare$}

\DeclareMathSymbol{\varnothing} {\mathord}{AMSb}{"3F} 

\newenvironment{notation}[0]{%
  \begin{list}%
    {}%
    {\setlength{\itemindent}{0pt}
     \setlength{\labelwidth}{4\parindent}
     \setlength{\labelsep}{\parindent}
     \setlength{\leftmargin}{5\parindent}
     \setlength{\itemsep}{0pt}
     }%
   }%
  {\end{list}}

\begin{document}

\title[Equidistribution and potential theory]
{Equidistribution of Small Points, Rational Dynamics, and Potential Theory}


\author[Matthew Baker]{Matthew H. Baker}
\author[Robert Rumely]{Robert Rumely}
\email{mbaker@@math.uga.edu \\
       rr@@math.uga.edu}
\address{Department of Mathematics,
         University of Georgia, Athens, GA 30602-7403, USA}
\address{Department of Mathematics,
         University of Georgia, Athens, GA 30602-7403, USA}


\thanks{The authors' research was supported by NSF Research Grant DMS-0300784.
The first author was also supported in part by an NSF Postdoctoral
Research Fellowship.  The authors would like to thank Laura DeMarco,
Xander Faber, and the anonymous referee for helpful comments on an earlier version of this manuscript.}

\begin{abstract}
  Given a dynamical system associated to a rational function $\varphi(T)$ 
  on $\PP^1$ of degree at least 2 with coefficients in a number field $k$, 
  we show that for each place $v$ of $k$, there is a unique probability 
  measure $\mu_{\varphi,v}$ on the Berkovich space $\PP^1_{\Berk,v} / \CC_v$ 
  such that if $\{ z_n \}$ is a sequence of points in $\PP^1(\kbar)$ whose
  $\varphi$-canonical heights tend to zero, then the $z_n$'s and their
  Galois conjugates are equidistributed with respect to $\mu_{\varphi,v}$.  
  In the archimedean case, $\mu_{\varphi,v}$ coincides with the well-known 
  canonical measure associated to $\varphi$.  This theorem generalizes a 
  result of Baker-Hsia \cite{BakerHsia} when $\varphi(z)$ is a polynomial.
  
  The proof uses a polynomial lift $F(x,y) = (F_1(x,y),F_2(x,y))$ of $\varphi$ to construct 
  a two-variable Arakelov-Green's function $g_{\varphi,v}(x,y)$ for each $v$.  
  The measure $\mu_{\varphi,v}$ is obtained by taking the Berkovich 
  space Laplacian of $g_{\varphi,v}(x,y)$, using a theory developed
  in \cite{RumelyNotes}.  
  The other ingredients in the proof 
  are (i) a potential-theoretic energy minimization principle 
  which says that $\iint g_{\varphi,v}(x,y) \, d\nu(x) d\nu(y)$ 
  is uniquely minimized over all probability measures $\nu$ on $\PP^1_{\Berk,v}$
  when $\nu = \mu_{\varphi,v}$, and 
  (ii) a formula for {\it homogeneous transfinite diameter} 
  of the $v$-adic filled Julia set $K_{F,v} \subset \CC_v^2$ 
  in terms of the resultant $\Res(F)$ of $F_1$ and $F_2$.
  The resultant formula, which generalizes a formula of DeMarco \cite{DeMarco},
  is proved using results from \cite{RLV} 
  about Chinburg's {\it sectional capacity}.  
  A consequence of the resultant formula is that the product  
  of the homogeneous transfinite diameters over all places is 1.  
\end{abstract}


\maketitle

\medskip


Let $k$ be a number field,
and let $\varphi(T) \in k(T)$ be a rational function of degree $d \ge 2$.
In this paper we investigate the equidistribution properties of small points
relative to the canonical dynamic height $\hhat_{\varphi}(z)$.  
We show that for each place $v$ of $k$, there
is a probability measure $\mu_{\varphi,v}$ such that
if $\{z_n\}$ is a sequence of distinct points
in $\PP^1(\kbar)$ satisfying 
$\hhat_{\varphi}(z_n) \rightarrow 0$,
then the Galois conjugates of the $z_n$ 
(regarded as embedded in $\PP^1(\CC_v)$)
are equidistributed relative to $\mu_{\varphi,v}$.
More precisely, if $\delta_n$ is the discrete probability measure
supported equally on the conjugates of $z_n$, then the sequence $\delta_n$
converges weakly to $\mu_{\varphi,v}$ for each $v$.
When $v$ is archimedean, $\mu_{\varphi,v}$ is the well-known canonical measure 
on $\PP^1(\CC)$ supported on the Julia set of $\varphi$ which was
constructed by Lyubich \cite{Lyubich} and Freire-Lopes-Ma{\~n}{\'e} \cite{FLM}.
(See \cite[\S 4]{Milnor} for the definition of the Julia set of a rational map.)

When $v$ is nonarchimedean, $\mu_{\varphi,v}$ is a measure
on the Berkovich space $\PP^1_{\Berk,v}$ over $\CC_v$
constructed by the authors in \cite{RumelyNotes}.  
It has the same invariance properties relative to $\varphi$ 
as the canonical measure in the archimedean case.

Conceptually, the proof is very simple.  Its main ingredients are
an energy-minimization principle, established at each place $v$
for the Arakelov Green's function $g_{\varphi,v}(x,y)$
assocated to $\mu_{\varphi,v}$, and two global inequalities,
an upper bound coming from the assumption that 
$\hhat_{\varphi}(z_n) \rightarrow 0$, 
and a lower bound coming from the product formula.
Combining these ingredients yields equidistribution
simultaneously at all places $v$.

At nonarchimedean places, 
the Arakelov Green's function $g_{\varphi,v}(x,y)$ is new.
We construct it by lifting $\varphi$ to a
polynomial map $F: \CC_v^2 \rightarrow \CC_v^2$,
and using the homogeneous local height associated
to the filled Julia set $K_{F,v}$ of this lift.
This approach was inspired by DeMarco \cite{DeMarco},
who introduced the {\it homogeneous capacity} $c^0(K)$ for
sets $K \subset \CC^2$, and proved for archimedean $v$ that
\begin{equation*}
c^0(K_{F,v}) \ = \ |\Res(F)|^{-1/d(d-1)} \ .
\end{equation*}
Although the homogeneous capacity does not easily generalize
to nonarchimedean places,
for arbitrary $v$ 
we introduce the closely related {\it homogeneous transfinite diameter}
$d^0_{\infty}(K_v)$ for sets $K_v \subset \CC_v^2$.
We generalize DeMarco's formula by showing that
\begin{equation}
d^0_{\infty}(K_{F,v}) \ = \ |\Res(F)|_v^{-1/d(d-1)} \label{FJJ1}
\end{equation}
for each $v$.  This is proved by relating the homogeneous transfinite diameter
to the {\it sectional capacity} studied in \cite{Ch} and \cite{RLV}.
The fact that $\prod_v d^0_{\infty}(K_{F,v}) = 1$,
which follows from the product formula applied to (\ref{FJJ1}),
is the key to the global lower bound mentioned above.

\smallskip

A philosophical idea which we hope to promote is the use of
Berkovich spaces as a natural setting for nonarchimedean Arakelov theory
and equidistribution theorems.  This point of view has been most 
strongly espoused by A.~Chambert-Loir \cite{CL}.  
The foundational results concerning potential theory
on the Berkovich projective line which are used in this paper can be found in 
\cite{RumelyNotes}.  Many of these results are proved for Berkovich curves of arbitrary 
genus in the doctoral thesis of A.~Thuillier \cite{Th}, 
a recent student of Chambert-Loir.  


\smallskip


\smallskip

P.~Autissier \cite{Au} has proved the archimedean part of the  
dynamical equidistribution theorem
using Arakelov-theoretic methods.
A proof of the nonarchimedean (Berkovich space) part of
the theorem, also based on ideas from Arakelov theory, 
has been announced by Chambert-Loir \cite{CL}. 

 C.~Favre and J.~Rivera-Letelier have also announced a 
proof of the dynamical equidistribution theorem.  Their 
preprints \cite{FRL}, \cite{FRL2} give
a proof of Theorem~\ref{DynamicalEquidistributionTheorem}
together with another construction
of the canonical measure on $\PP^1_{\Berk,v}$
attached to a rational map $\varphi$.  
They also prove a Berkovich space analogue of Theorem~\ref{LyubichTheorem} below.
The technical foundation for their work can be found in the 
monograph by Favre and Jonsson \cite{FJ}, 
in Rivera-Letelier's thesis \cite{R-L1}, 
and in a manuscript of Rivera-Letelier \cite{RLTFJ}.  
As with our approach, Favre and Rivera-Letelier's proof of 
of Theorem~\ref{DynamicalEquidistributionTheorem} 
is ultimately based on the product formula and 
an adelic energy-minimization theorem.  
However, there are also a number of differences between the
two proofs.

Finally we note that in the nonarchimedean case,
the construction of J.~Piniero, L.~Szpiro and T.~Tucker \cite{PST},  
which works scheme-theoretically with blowups of models of 
$\PP^1/\Spec(\cO_v)$ attached to iterates of $\varphi(T)$,
yields a sequence of discrete measures which 
can be shown to converge to the canonical measure 
on $\PP^1_{\Berk,v}$.    

\section{Notation}
\label{section:Notation}

\medskip

We set the following notation and normalizations, which will be used
throughout the paper unless otherwise noted.

\begin{notation}
\item[$k$]
a number field.
\item[$\O_k$]
the ring of integers of $k$.
\item[$M_k$]
the set of places of~$k$. 
\item[$k_v$]
the completion of $k$ at $v$.
\item[$\O_v$]
the ring of integers in $k_v$.
\item[$q_v$]  the order of the residue field of $k_v$.  If $v$ is archimedean,
we put $q_v = e$ if $k_v \cong \RR$, and $q_v = e^2$ if $k_v \cong \CC$.
\item[$\CC_v$] 
the completion of a fixed algebraic closure $\kvbar$ of $k_v$.  
Throughout the paper, we fix a choice of an embedding of $\kbar$ into $\CC_v$
for each $v \in M_k$ (though all of our conclusions will be
independent of the choices made).  If $v$ is nonarchimedean, we write
$\hcO_v$ for the ring of integers of $\CC_v$.  
\item[$|x|_v$] 
the canonical absolute value on $k_v$ given by
the modulus of additive Haar measure.
If $| x|'_v$ is the unique absolute value on~$k$ in the equivalence class of~$v\in
M_k$ that extends the standard absolute value on the
completion~$\QQ_v$, then $|x|_v = (|x|'_v)^{[k_v : \QQ_v]}$.
With this normalization, the product formula holds in the form
$\prod_v |\alpha|_v = 1$ for each $\alpha \neq 0$ in $k$.  
Each $|x|_v$ extends uniquely to an absolute value on $\CC_v$, the completion
of the algebraic closure of $k_v$.

\item[$h$]
the absolute logarithmic Weil height
$h:\PP^n(\QQbar)\to\RR$, defined for $[x_0:\dots:x_n]\in\PP^n(k)$ by
\[
  h\bigl([x_0:\dots:x_n]\bigr) = \frac{1}{[k:\QQ]} \sum_{v \in M_k} 
  \log \max \{ |x_0|_v, \ldots, |x_n|_v \}.
\]
\item[$\varphi$] A rational function on $\PP^1$ defined over $k$.
\end{notation}


\section{Overview}

%
%

%
%


\subsection{An equidistribution result for rational functions on $\PP^1$.}
\label{RationalFunctionSection}

\medskip

Let $\varphi : \PP^1 \to \PP^1$ be a finite morphism (i.e., a nonconstant rational
function) of degree $d \geq 2$ defined over the number field $k$.
Iterating $\varphi$ gives rise to a dynamical system on $\PP^1(\CC_v)$
for all places $v$ of $k$.  When $v$ is archimedean, 
this type of dynamical system has been extensively studied
since the pioneering work of Fatou and Julia in the early 20th century.
Just as one defines the N{\'e}ron-Tate canonical height on an elliptic
curve by iteration, one can define the {\em dynamical height} 
\[
\hhat_\varphi : \PP^1(\kbar) \to \RR
\]
attached to the rational function $\varphi$ by the rule
\[
\hhat_\varphi(z) \ = \ \lim_{n\to\infty} \frac{1}{d^n} h(\varphi^{(n)}(z)).
\]
Here $\varphi^{(n)}$ denotes the $n$-fold iterate $\varphi \circ \cdots \circ \varphi$.
By a general result of Call and Silverman \cite{CS}, the hypothesis $d
\geq 2$ guarantees that the above limit exists.

The dynamical height $\hhat_\varphi$ is
uniquely characterized by the following two properties:
\begin{itemize}
\item[(1)] The difference $|\hhat_\varphi - h|$ is bounded.
\item[(2)] $\hhat_\varphi \circ \varphi = d \cdot \hhat_\varphi$.
\end{itemize}
It follows from \cite {CS}
that $\hhat_\varphi(z) \geq 0$ for all $z \in \PP^1(\kbar)$, and
$\hhat_\varphi(z) = 0$ if and only if $z$ is {\em preperiodic} for
$\varphi$, meaning that the orbit $\{ \varphi^{(n)}(z) : n \in \NN \}$ of
$z$ under iteration of $\varphi$ is a finite set.
Additionally, we have $\hhat_\varphi(\sigma z) = \hhat_\varphi(z)$ for all
$z \in \PP^1(\kbar)$ and all $\sigma \in \Gal(\kbar / k)$.

\medskip

%


If $\varphi(z) = z^2$, then $\hhat_\varphi$ is the usual
logarithmic Weil height $h$ on $\PP^1(\QQbar)$.  Another
well-known height which can be defined by dynamical methods is the 
N\'eron-Tate canonical height on an elliptic curve.  If $k$ is a number field and
$E / k$ is an elliptic curve with Weierstrass equation
$y^2 = f(x)$, let $\varphi$ be the degree 4 rational function on $\PP^1$
given by $x \circ [2]$.  Then for $P \in E(\kbar)$ we have
$\hhat(P)= \hhat_\varphi(x(P))$.

\medskip

For any rational function $\varphi$ on $\PP^1$ of degree $d\geq 2$
defined over $\CC$, Lyubich \cite{Lyubich}, and independently Freire,
Lopes, and Ma{\~n}{\'e} \cite{FLM}, constructed a natural probability
measure $\mu_\varphi$ attached to the dynamical system $\{ \varphi^{(n)} : n
\in \NN \}$.  We will refer to the measure $\mu_\varphi$ as the {\em
canonical measure} attached to $\varphi$.  In order to characterize $\mu_\varphi$, we recall
the following definition.  A point $z_0 \in \PP^1(\CC)$ is said to be
{\em exceptional} if the set $\{ \varphi^{(-n)}(z_0) : n \in \NN \}$ of
backward iterates of $z_0$ is finite.  It is known (see \cite{Milnor})
that there are at most 2 exceptional points for $\varphi$ in
$\PP^1(\CC)$.  Proofs of the following theorem can be found in
\cite{Lyubich}, \cite{FLM}, and \cite{HP}.

\begin{theorem}
\label{LyubichTheorem} 

There exists a probability measure $\mu_\varphi$ $($independent of
$z_0$$)$ such that:

$A)$
For any non-exceptional point $z_0 \in \PP^1(\CC)$, let $\delta_n$ be
the probability measure
\[
\frac{1}{d^n} \sum_{\varphi^{(n)}(z) = z_0} \delta_z,
\]
where the points in the sum are counted with multiplicities and
$\delta_z$ denotes the Dirac measure on $\PP^1(\CC)$ giving mass 1 to
the point $z$.  Then the sequence of measures $\delta_n$ converges
weakly to $\mu_\varphi$.

$B)$ $\mu_\varphi$ is the unique measure on $\PP^1(\CC)$ with no point
masses such that
$\varphi^*(\mu_{\varphi}) = d \cdot \mu_{\varphi}$ as $(1,1)$-currents.

\end{theorem}

When $\varphi$ is a polynomial, Theorem~\ref{LyubichTheorem} was
originally proved by Brolin, and the measure $\mu_\varphi$ is known in
that case as {\em Brolin's measure}.  Brolin's measure coincides with
the equilibrium measure (in the sense of potential theory)
on the Julia set of $\varphi$.

\medskip

We will now briefly recall the construction of the 
Berkovich space $\PP^1_{\Berk,v}$
associated to the projective line over $\CC_v$, where $v$ is a nonarchimedean
place of $k$.

The Berkovich unit disc $\cB(0,1)$ is 
the set of all continuous multiplicative seminorms on the Tate algebra
$\CC_v \< T \>$ (see \cite[\S 1.4]{Be}, \cite[\S 1]{RumelyNotes}).
Examples of elements of $\cB(0,1)$ include
the evaluation seminorms $[f]_a = |f(a)|_v$ for $a \in \CC_v$ with
$|a|_v \le 1$;  $\sup$ norms $[f]_{B(a,r)} = \sup_{z \in B(a,r)} |f(z)|_v$
for discs $B(a,r) = \{z \in \CC_v : |z-a|_v \le r \}$;
and limit norms associated to nested sequences of discs
$B(a_1,r_1) \supset B(a_2,r_2) \supset \cdots$, defined by 
\begin{equation*}
[f]_x \ = \ \lim_{i \rightarrow \infty} [f]_{B(a_i,r_i)} \ .
\end{equation*}
A theorem of Berkovich says that all continuous multiplicative seminorms
on $\CC_v \< T \>$ arise in this way.
Following Chambert-Loir \cite{CL}, we call the
point $\zeta_0 \in \cB(0,1)$ corresponding to the Gauss norm
$\|f\| = [f]_{B(0,1)}$ the {\em Gauss point}.
Given a point $x \in \cB(0,1)$ corresponding either to a disc $B(a,r)$ or to a point 
$a = B(a,0)$ (which can be thought of as a degenerate disc),
there is a path $\{ [ \ ]_{B(a,t)} : r \le t \le 1 \}$
connecting $x$ to the Gauss point.  Given a collection of discs,
the union of the corresponding paths forms a subtree of $\cB(0,1)$
rooted at $\zeta_0$.  From this, one sees that $\cB(0,1)$ is an infinitely
branched real tree, with countably many branches emanating from each point
corresponding to a disc with radius $r \in |\CC_v^{\times}|_v$.

As a set, the Berkovich projective line $\PP^1_{\Berk,v}$ over $\CC_v$
is obtained by gluing together two copies of  $\cB(0,1)$.
It is made into a topological space by equipping 
it with the {\em Gelfand topology}, the weakest topology such that
each set of the form
\[
U_{a,b}(f) \ = \ \{x \in \PP^1_{\Berk,v} : a < [f]_x < b \}
\]
for $a, b \in \RR$ and $f \in \CC_v(T)$ is open.
The space $\PP^1_{\Berk,v}$ is also equipped with a sheaf of rings $\cO_{X}$,
constructed using localizations of Tate algebras;
see \cite{Be} for details.  There is a natural inclusion
$\PP^1(\CC_v) \subset \PP^1_{\Berk,v}$ (which associates to a point of
$\PP^1(\CC_v)$ the corresponding ``evaluation seminorm'') that induces the
usual (ultrametric) topology on $\PP^1(\CC_v)$, 
and $\PP^1(\CC_v)$ is dense in $\PP^1_{\Berk,v}$ under this inclusion. 

If $\varphi(T) \in \CC_v(T)$ is a nonconstant rational function,
then $\varphi$ acts on $\PP^1_{\Berk}$ by
\[
[f]_{\varphi(x)} \ = \ [f \circ \varphi]_x
\]
for all $f \in \CC_v(T)$.  This coincides with the usual action of
$\varphi$ on $\PP^1(\CC_v) \subset \PP^1_{\Berk,v}$.

As a topological space,
$\PP^1_{\Berk,v}$ is compact, Hausdorff, and path-connected, in
contrast with $\PP^1(\CC_v)$, which is completely disconnected and not even
locally compact.   Thus $\PP^1_{\Berk,v}$ is a much more suitable space
for doing measure theory and potential theory than $\PP^1(\CC_v)$.
The space $\PP^1_{\Berk,v}$ is also metrizable, although there is not a canonical
metric on it.


\begin{remark}
If $v$ is archimedean, one can define $\PP^1_{\Berk,v}$ over $\CC$ in
a similar way using continuous multiplicative seminorms 
on $\CC \< T \>$.
By the Gelfand-Mazur theorem, every such seminorm
arises from evaluation at a point.
Thus $\PP^1_{\Berk,v}/\CC$ is isomorphic to $\PP^1(\CC)$.
\end{remark}

In Theorem~\ref{LyubichTheorem},
note that if $z_n \in \varphi^{(-n)}(z_0)$, then
\[
\hhat_\varphi(z_n) = \frac{1}{d^n} \hhat_\varphi(z_0) \to 0
\]
as $n\to\infty$.  Also, note that if $\varphi$ is defined over the number
field $k$ and $z_0 \in k$, then the set $\varphi^{(-n)}(z_0)$ is stable
under $\Gal(\kbar / k)$.  We will prove the following adelic equidistribution
theorem, motivated by Theorem~\ref{LyubichTheorem} and by the archimedean 
equidistribution theorems of Bilu \cite{Bilu} and Szpiro-Ullmo-Zhang \cite{SUZ}.

\begin{theorem}[Main Theorem]
\label{DynamicalEquidistributionTheorem} 
For each place $v \in M_k$, there exists a canonical probability
measure $\mu_{\varphi,v}$ on the Berkovich space $\PP^1_{\Berk,v} / \CC_v$
such that the following holds: Suppose $z_n$ is a sequence of distinct
points of $\PP^1(\kbar)$ with $\hhat_\varphi(z_n) \to 0$.  For $v \in
M_k$, let $\delta_n$ be the discrete probability measure on the
Berkovich space $\PP^1_{\Berk,v} / \CC_v$ supported equally on the Galois
conjugates of $z_n$.  Then the sequence of measures $\delta_n$
converges weakly to $\mu_{\varphi,v}$ for all $v \in M_k$.
\end{theorem}


\medskip




When $\varphi$ is a polynomial, the archimedean part of
Theorem~\ref{DynamicalEquidistributionTheorem} was proved by
Baker-Hsia in \cite{BakerHsia}.  The present paper provides a
conceptual simplification of their method, and 
applies to arbitrary rational functions.  We note that the case of a
rational function is more difficult than the polynomial
case, due to the absence of a fixed pole at infinity.  A weaker
version of the nonarchimedean part of
Theorem~\ref{DynamicalEquidistributionTheorem}, formulated in terms of
``pseudo-equidistribution'', was also proved for the polynomial case
in Baker-Hsia in \cite{BakerHsia}.  Here we clarify the
meaning of pseudo-equidistribution by using Arakelov Green's functions
and Berkovich spaces.



When $\varphi(z) = z^2$,
the archimedean part of Theorem~\ref{DynamicalEquidistributionTheorem}
specializes to (and was motivated by) the following well-known result of Bilu:

\begin{theorem}[Bilu \cite{Bilu}]
\label{BiluTheorem}
Let $z_n$ be a sequence of distinct points in $\PP^1(\QQbar)$, and
suppose that $h(z_n) \to 0$.  Let $\delta_n$ be the discrete
probability measure on $\PP^1(\CC) = \CC \cup \{ \infty \}$ which is
supported with equal mass at each Galois conjugate of $z_n$.  Then the
sequence of measures $\delta_n$ converges weakly to the uniform
probability measure $\mu_{S^1}$ on the unit circle $\{ |z| = 1 \}$.
\end{theorem}

For previous explorations of the
relationship between Bilu's theorem and potential theory, see
\cite{Bombieri} and \cite{RumelyBilu}.


\medskip

\section{Adelic dynamics on $\PP^1$}
\label{AdelicDynamicsSection}

\subsection{Dynamical heights associated to rational functions.}
\label{DynamicalHeightSection}

\medskip

Recall that $\varphi : \PP^1 \to \PP^1$ is a rational function of degree
$d \geq 2$ defined over a number field $k$.

The map $\varphi$ can be represented in
homogeneous coordinates as
\[
\varphi([z_0:z_1]) \ = \ [F_1(z_0,z_1):F_2(z_0,z_1)]
\]
for some homogeneous polynomials $F_1,F_2 \in k[x,y]$ of degree $d$
with no common linear factor over $\kbar$.
(Note that since $F_1$ and $F_2$ factor into linear terms over $\kbar$, 
$F_1$ and $F_2$ have a common factor over $\kbar$ if and only if they
have a common linear factor over $\kbar$.)
The polynomials $F_1,F_2$ are uniquely determined by $\varphi$ up to
multiplication by a common scalar $c \in k^*$.

Dehomogenizing by setting $z = z_1 / z_0$, we obtain 
\[
\varphi(z) \ = \ \frac{f_2(z)}{f_1(z)}
\]
with $f_i \in k[z]$ and $\max \{ {\rm deg}(f_1), {\rm deg}(f_2) \} = d$.

\medskip

We will often want to work with the degree $d$ homogeneous polynomials
$F_1$ and $F_2$, so we now fix a choice of $F_1,F_2 \in k[X,Y]$.  This
allows us to consider the mapping
\[
F = (F_1,F_2) : \AA^2(\kbar) \to \AA^2(\kbar)
\]
as a global lifting of $\varphi$.

Let $\Res(F) := \Res(F_1,F_2)$ denote the homogeneous resultant of the
polynomials $F_1$ and $F_2$ (see e.g. \cite[\S 6]{DeMarco}).
Since $F_1$ and $F_2$
have no common linear factor over $\kbar$, we have $\Res(F) \neq 0$,
and $F(z_0,z_1)=(0,0)$ if and only if $(z_0,z_1)=(0,0)$.

In the archimedean case, define
$\|(z_0,z_1)\|_v = \sqrt{|z_0|^2+|z_1|^2}^{[k_v:\RR]}$;
in the nonarchimedean case, put
$\|(z_0,z_1)\|_v := \max \{ |z_0|_v,|z_1|_v \} $.
We begin with the following simple lemma.

\begin{lemma} \label{Lem1}
For each place $v$ of $k$, 
there are constants $0 < C_v \le D_v$ such that for all
$z \in \CC_v^2$,
\begin{equation} \label{FA2}
C_v (\|z\|_v)^d \ \le \ \|F(z)\|_v \ \le \ D_v (\|z\|_v)^d \ .
\end{equation}
For all but finitely many $v$, we may take $C_v = D_v = 1$.
\end{lemma}

\begin{proof}
First suppose $v$ is archimedean.  Identify $\CC_v$ with $\CC$. Since
$\partial B_v(1) = \{(x,y) \in \CC^2 : \max(|x|,|y|) = 1\}$ is compact,
and since the only common zero of $F_1(z)$ and $F_2(z)$ is the origin,
the constants 
\begin{equation*}
C_v  =  \min_{z \in \partial B_v(1)} \|F(z)\|_v \ , \quad
D_v  =  \max_{z \in \partial B_v(1)} \|F(z)\|_v 
\end{equation*}
satisfy $0 < C_v \le D_v$.
By homogeneity, (\ref{FA2}) holds for all $z \in \CC^2$.

Now let $v$ be nonarchimedean.
Write 
\begin{equation*}
\partial B_v(1) = \{(x,y) \in \CC_v^2 : \max(|x|_v,|y|_v) = 1\} \ .
\end{equation*}
Let $C_v^{\prime} := |\Res(F)|_v$, and let $D_v$ be an upper bound for
the absolute values of the coefficients of $F_1$ and $F_2$.  

By a well-known property of the resultant of two homogeneous
polynomials of degree $d$ \cite[\S 5.8]{VdW}, there exist polynomials 
$g_1(x,y)$, $g_2(x,y)$, and $h_1(x,y)$, $h_2(x,y)$, which are
homogeneous of degree $d-1$ in $x$ and $y$ and whose coefficients lie
in $k$, such that
\begin{eqnarray*}
g_1(x,y) F_1(x,y) + g_2(x,y) F_2(x,y) & = & \Res(F) x^{2d-1} \ , \\
h_1(x,y) F_1(x,y) + h_2(x,y)F_2(x,y)  & = & \Res(F) y^{2d-1} \ .
\end{eqnarray*}
For each $(x,y) \in \partial B_v(1)$, all of $|g_1(x,y)|_v$,
$|g_2(x,y)|_v$, $|h_1(x,y)|_v$ and $|h_2(x,y)|_v$ are $\le C_v^{\prime \prime}$
for some constant $C_v^{\prime \prime}>0$ independent of $(x,y)$.  
Furthermore, we may take $C_v^{\prime \prime} = 1$ for almost all $v$.
By the ultrametric inequality,
\begin{eqnarray*}
C_v^{\prime} |x|_v^{2d-1} & \le &
C_v^{\prime \prime} \max(|F_1(x,y)|_v,|F_2(x,y)|_v) \ , \\
C_v^{\prime} |y|_v^{2d-1} & \le &
C_v^{\prime \prime} \max(|F_1(x,y)|_v,|F_2(x,y)|_v) \ .
\end{eqnarray*}
Put $C_v := C_v^{\prime} / C_v^{\prime \prime}$. 
Then for each $z = (x,y) \in \partial B_v(1)$,
\begin{equation*}
C_v \cdot \max(|x|_v,|y|_v)^{2d-1} \ \le \ \|F(z)\|_v \ .
\end{equation*}
However, if $z \in \partial B_v(1)$ then $\max(|x|_v,|y|_v) = \|z\|_v = 1$,
so $C_v \le \|F(z)\|_v$ for all $z \in \partial B_v(1)$.
The first inequality in (\ref{FA2}) follows by homogeneity.  The second
follows trivially by the ultrametric inequality.

Finally, since the resultant and the coefficients of the $F_i$
are elements of $k$, independent of $v$,
we can take $C_v = D_v = 1$ for all but finitely many $v$.
\end{proof}

\begin{remark}
\label{Lem1Remark}
If $v$ is nonarchimedean and $F_1, F_2$ have $v$-integral
coefficients, then we may choose $g_1(x,y)$, $g_2(x,y)$, and
$h_1(x,y)$, $h_2(x,y)$ to have $v$-integral coefficients as well.
In this case we can take $C_v = |\Res(F)|_v$ and $D_v = 1$.
\end{remark}

\begin{cor} \label{Cor1}
For each $v$, there are radii $0 < r_v \le R_v$ such that
for each $z \in \CC_v^2$ with $\|z\|_v \le r_v$,
we have $\|F(z)\|_v \le \|z\|_v \cdot (\|z\|_v/r_v)^{d-1}$,
and for each $z$ with $\|z\|_v \ge R_v$,
we have $\|F(z)\|_v \ge \|z\|_v \cdot (\|z\|_v/R_v)^{d-1}$.  In particular,
\begin{eqnarray*}
F(B_v(r_v)) & \subseteq & B_v(r_v) \ , \\
F(\CC_v^2 \backslash B_v(R_v)) & \subseteq & \CC_v^2 \backslash B_v(R_v) \ .
\end{eqnarray*}
For all but finitely many $v$, we can take $r_v = R_v = 1$.
\end{cor}

\begin{proof}  Let $C_v$ and $D_v$ be as in Lemma \ref{Lem1}, and put
$r_v = D_v^{-1/(d-1)}$,  $R_v = C_v^{-1/(d-1)}$.  If $\|z\|_v \le r_v$, then
\begin{eqnarray*}
\|F(z)\|_v & \le & D_v \|z\|_v^d  \  = \ r_v^{-(d-1)} \|z\|_v^d \\
& = & \|z\|_v \cdot (\|z\|_v/r_v)^{d-1} \ .
\end{eqnarray*}
Similarly if $\|z\|_v \ge R_v$, then
\begin{eqnarray*}
\|F(z)\|_v & \ge & C_v \|z\|_v^d \ = \ R_v^{-(d-1)} \|z\|_v^d \\
& = & \|z\|_v \cdot (\|z\|_v/R_v)^{d-1} \ .
\end{eqnarray*}
For each $v$ with $C_v = D_v = 1$, we have $r_v = R_v = 1$.
\end{proof}


\medskip

Recall that the global dynamical height
$\hhat_\varphi : \PP^1(\kbar) \to \RR$ is defined by
\[
\hhat_\varphi(z) = \lim_{n \to \infty} \frac{1}{d^n} h(\varphi^{(n)}(z)).
\]
The choice of a global lifting $F$ of $\varphi$ allows us to decompose
the global dynamical height into a sum of local heights as follows.

For $v \in M_k$ and $z = (z_0,z_1) \in \CC_v^2 \backslash \{ 0 \}$,
define the {\em homogeneous local dynamical height} $\Hhat_{F,v} : \CC_v^2
\backslash \{ 0 \} \to \RR$ by
\[
\Hhat_{F,v}(z) := \lim_{n\to\infty} \frac{1}{d^n} \log \| F^{(n)}(z) \|_v.
\]

By convention, we define $\Hhat_{F,v}(0,0) := -\infty$.

\begin{lemma}
The limit $\lim_{n\to\infty} \frac{1}{d^n} \log \| F^{(n)}(z) \|_v$
exists for all $z \in \CC_v^2 \backslash \{ 0 \}$, and 
$\frac{1}{d^n} \log \| F^{(n)}(z) \|_v$ converges uniformly on
$\CC_v^2 \backslash \{ 0 \}$ to $\Hhat_{F,v}(z)$.
\end{lemma}

\begin{proof}
The proof is by a standard telescoping series argument
(see e.g. \cite{HP}).
Define
\begin{eqnarray*}
T_j(z) & := & 
\frac{1}{d^{j+1}} \log \| F^{(j+1)}(z) \|_v
- \frac{1}{d^j} \log \| F^{(j)}(z) \|_v \\
& = & \frac{1}{d^j} \left( \frac{1}{d} \log \| F^{(j+1)}(z) \|_v
- \log \| F^{(j)}(z) \|_v \right).
\end{eqnarray*}
By (\ref{FA2}) we see that
\[
C_v (\| F^{(j)}(z) \|_v)^d
\ \le \ \|F^{(j+1)}(z)\|_v \ \le \
D_v (\| F^{(j)}(z)\|_v)^d
\]
for all $z \in \CC_v^2$. Applying these inequalities to the sequence
$T_j$ yields the estimate
\[
|T_j(z)| \ \leq \ \frac{C}{d^{j+1}},
\]
where $C := \max \{ \log D_v, -\log C_v \}$.

It follows easily that $H_n := \sum_{j = 0}^n T_j$ is a Cauchy
sequence, and that the series defining $\Hhat_{F,v}(z)$
converges uniformly on $\CC_v^2 \backslash \{ 0 \}$.
\end{proof}


\medskip

\begin{remark}
\label{BoundedDifferenceRemark}
The proof yields the explicit bound
\[
| \Hhat_{F,v}(z) - \log \| z \|_v | \ \leq \ \frac{C}{d - 1},
\]
valid for all $z \in \CC_v^2 \backslash \{ 0 \}$.
\end{remark}

Note that by the definitions of the local and global canonical
heights, if $x \in \PP^1(k)$, then for any representation
$x = [x_0:x_1]$ with $x_0, x_1 \in k$, we have
\[
\hhat_\varphi(x) = \frac{1}{[k:\QQ]} \sum_{v \in M_k} \Hhat_{F,v}(x_0,x_1).
\]
By the product formula, the right side is independent of the
choice of lifting.

Also, note that the definition of $\Hhat_{F,v}$ is independent of the
norm used to define it.  This follows easily from the equivalence of
norms on $\CC_v^2$.  

\medskip

The homogeneous local dynamical height $\Hhat_{F,v}$
has the following properties, and in fact is uniquely
characterized by them:

\begin{itemize}
\item[(LH1)] The difference $|\Hhat_{F,v}(z) - \log \|z \|_v|$ is bounded.
\item[(LH2)] $\Hhat_{F,v}(F(z)) =  d \cdot \Hhat_{F,v}(z)$.
\item[(LH3)] $\Hhat_{F,v}$ {\em scales logarithmically}, i.e.,
for all $c \in \CC_v^*$,
\[
\Hhat_{F,v}(cz) \ = \ \Hhat_{F,v}(z) + \log |c|_v \ .
\]
\end{itemize}

\subsection{The filled Julia set.}

\medskip

By definition, the {\emph{filled Julia set}} $K_{F,v}$ of $F$ in $\CC_v^2$
is the set of all $z \in \CC_v^2$ for which the iterates $F^{(n)}(z)$
remain bounded.  
Clearly $F^{-1}(K_{F,v}) = K_{F,v}$, and the same is true for each $F^{(-n)}$.
Since all norms on $\CC_v^2$ are equivalent, the set $K_{F,v}$ is
independent of which norm is used to define it.


By Corollary~\ref{Cor1}, we have $B_v(r_v) \subseteq K_{F,v}$, so $K_{F,v}$ cannot be
too small.  Moreover:

\begin{lemma}
\label{FilledJuliaIntersectionLemma}
With $R_v$ as in Corollary~\ref{Cor1}, we have
\[
F^{(-1)}(B_v(R_v)) \supseteq F^{(-2)}(B_v(R_v)) \supseteq \cdots
\]
and
\begin{equation*}
K_{F,v} \ = \ \cap_{n=1}^{\infty} F^{(-n)}(B_v(R_v)).
\end{equation*}
\end{lemma}

\begin{proof}
The assertion that  $F^{(-n)}(B_v(R_v)) \supseteq
F^{(-n-1)}(B_v(R_v))$ for all $n\geq 1$
is equivalent to the statement that if $F^{(n)}(z) \not\in B_v(R_v)$ then
$F^{(n+1)}(z) \not\in B_v(R_v)$.
This follows from Corollary~\ref{Cor1}, which shows that 
$F(\CC_v^2 \backslash B_v(R_v)) \subseteq \CC_v^2 \backslash
B_v(R_v)$.

It also follows from Corollary~\ref{Cor1} that if $\|z\|_v > R_v$, then
\begin{equation*}
\lim_{n \rightarrow \infty} \|F^{(n)}(z)\|_v \ = \ \infty \ .
\end{equation*}
Thus, $K_{F,v} \subset B_v(R_v)$.
By iteration, $K_{F,v} \subset F^{(-n)}(B_v(R_v))$ for each $n$.

The fact that
$K_{F,v} \ = \ \cap_{n=1}^{\infty} F^{(-n)}(B_v(R_v))$
now follows, since if $z \notin K_{F,v}$, there is some $n$ for which 
$\|F^{(n)}(z)\|_v > R_v$, and so $z \notin F^{(-n)}(B_v(R_v))$.
\end{proof}

\medskip

The filled Julia set $K_{F,v}$ can be thought of as the `unit ball'
with respect to the dynamical local height $\Hhat_{F,v}$:

\begin{lemma}
For each place $v$ of $k$,
\[
K_{F,v} = \{ z \in \CC_v^2 \; : \; \Hhat_{F,v}(z) \leq 0 \}.
\]
\end{lemma}

\begin{proof}
If $z \in K_{F,v}$ then there exists $M>0$ such that $\| F^{(n)}(z)
\|_v \leq M$ for all $n$, and therefore
$\Hhat_{F,v}(z) \leq \lim_{n\to\infty} \frac{1}{d^n} \log M = 0$.

Conversely, suppose $z \not\in K_{F,v}$.  Then for $n_0$ sufficiently
large, $\beta := \| F^{(n_0)}(z) \|_v > R_v$.  Let
$\alpha := \beta / R_v > 1$.  Then by Corollary~\ref{Cor1} and
induction on $n$, it follows that
\[
\| F^{(n + n_0)}(z) \|_v > \beta \cdot \alpha^{d^n - 1}
\]
for all $n\geq 0$.  

Therefore 
\[
\Hhat_{F,v}(z) \ \geq \ \lim_{n\to\infty}
\frac{1}{d^{n+n_0}} \left( (d^n - 1)\log \alpha + \log \beta \right) 
\ = \ \frac{1}{d^{n_0}} \log \alpha \ > \ 0.
\]
\end{proof}

\medskip

In general, it is difficult to describe $K_{F,v}$ explicitly.
However, the following lemma shows that it is `trivial'
for all but finitely many $v$.

\begin{lemma}
\label{UnitPolydiskLemma}
Suppose $F_1,F_2 \in \O_k[x,y]$.
If $v$ is a nonarchimedean place of $k$ such that
$|\Res(F)|_v = 1$, then $K_{F,v} = B(0,1)^2$ is the unit polydisk in
$\CC_v^2$ and $\Hhat_{F,v}(z) = \log \| z \|_v$ for all $z \in \CC_v^2$.
\end{lemma}

\begin{proof}
By Remark~\ref{Lem1Remark}, it follows that
$\|F(z)\|_v = \|z\|_v^d$, and therefore 
$\|F^{(n)}(z)\|_v = \|z\|_v^{d^n}$ for all $n\geq 1$.  The result
follows immediately.
\end{proof}

\medskip

\subsection{The homogeneous transfinite diameter.}

\medskip

Let $v$ be a place of $k$, and let $K \subset \CC_v^2$ be a nonempty bounded set.
For $z = (z_0,z_1), w = (w_0, w_1) \in \CC_v^2$, put
\[
z \wedge w \ := \ z_0 w_1 - z_1 w_0 \ .
\]
By analogy with the classical transfinite diameter, for $n\geq 2$ we let
\[
d^0_n(K) \ := \ \sup_{z_1,\ldots,z_n \in K}
\left( \prod_{i \neq j} |z_i \wedge z_j|_v  \right)^{\frac{1}{n(n-1)}}.
\]


\begin{lemma}
\label{HomogeneousMonotonicLemma}
  The sequence of nonnegative real numbers $d^0_n(K)$ is
  non-increasing.  In particular, the quantity $d^0_\infty(K) :=
  \lim_{n\to\infty} d^0_n(K)$ is well-defined.
\end{lemma}

\begin{proof}
We claim that $d_n^0(K) \ge d_{n+1}^0(K)$ for all $n$.
The proof is the same as for the classical transfinite diameter:
write $P_n(z_1, \ldots, z_n) = \prod_{i \neq j} (z_i \wedge z_j)$, take
$\varepsilon > 0$, and choose $w_1, \ldots,w_{n+1} \in K$ with
$|P_{n+1}(w_1, \ldots, w_{n+1})|_v
\ge (d_{n+1}^0(K)-\varepsilon)^{n(n+1)}$.
For each $\ell = 1, \ldots, n+1$ write $\hat{w}_{\ell} =
(w_1, \ldots, w_{\ell-1}, w_{\ell+1}, \ldots, w_{n+1})$.  
By definition, 
$d_n^0(K)^{n(n-1)} \ge |P_n(\hat{w}_{\ell})|_v$ for each $\ell$.
It follows that
\begin{eqnarray*}
d_n^0(K)^{(n+1)n(n-1)}
& \ge & \prod_{\ell=1}^{n+1} |P_n(\hat{w}_{\ell})|_v 
\ = \ \big(\prod^{n+1}_{i \neq j} 
              |w_i \wedge w_j|_v \big)^{n-1} \\
& \ge & (d_{n+1}^0(K)-\varepsilon)^{(n+1)n(n-1)} \ .
\end{eqnarray*}
This holds for each $\varepsilon > 0$, so $d_n^0(K) \ge d_{n+1}^0(K)$.
\end{proof}

\medskip

We call $d^0_\infty(K)$ the {\em homogeneous transfinite diameter} of
$K$.

\medskip

We will now show that when $v$ is archimedean and $K \subset \CC^2$ is compact,
$d^0_\infty(K)$ coincides with the {\em homogeneous capacity} introduced
by DeMarco \cite{DeMarco}.  We recall the definition.

\begin{defn}
If $K \subseteq \CC^2$ is compact and nonempty, define $c^0(K)$ by 
\[
-\log c^0(K) \ := \ \inf_{\nu \in \PP(K)} I(\nu),
\]
where $\PP(K)$ is the space of probability measures supported on $K$,
and 
\[
I(\nu) \ := \ -\iint_{K \times K} \log |z\wedge w| \, d\nu(z) d\nu(w).
\]
The quantity $c^0(K)$ is called the {\em homogeneous capacity} of $K$.
\end{defn}

\begin{lemma} \label{EqualityLemma}
If $v$ is archimedean and
$K \subset \CC^2$ is compact, then $d^0_\infty(K) = c^0(K)^{[k_v:\RR]}$.
\end{lemma}

\begin{proof}
Note that $[k_v:\RR] = 1$ or $2$, according as
$k_v \cong \RR$ or $k_v \cong \CC$.
The power $[k_v:\RR]$ arises because of our normalization of absolute values:
$|x|_v = |x|^{[k_v:\RR]}$ for all $x \in  \CC$.  If we replace $|x|_v$ with
$|x|$ in the definition of $d^0_{\infty}(K)$, it suffices to show that
$d^0_{\infty}(K) = c^0(K)$.

By a general fact about measures proved in Lemma~\ref{LimInfLemma} below, 
if $\nu_n$ is a probability measure supported equally on
$z_1,\ldots,z_{N_n} \in K$, with $z_i \wedge z_j \neq 0$ for all $i\neq j$,
and if $\nu_n \to \nu$ weakly on $K$, then

\begin{equation}
\label{LimInfWedgeEquation}
\liminf_{n\to\infty}
\iint_{\CC^2 \times \CC^2 \backslash \text{(Diag)}}
 -\log |z \wedge w| \, d\nu_n(z) d\nu_n(w)
\ \geq \ \iint_{\CC^2 \times \CC^2}  -\log |z \wedge w| d\nu(z) d\nu(w) \ .
\end{equation}

Let $n\geq 2$, and define $D_n = -\log d^0_n(K)$.
Given any $z_1,\ldots,z_n \in K$, we have
\[
\frac{1}{n(n-1)} \sum_{i \neq j} -\log |z_i \wedge z_j| \ \geq \ D_n
\]
by definition.
Integrating this inequality against an arbitrary measure $\mu \in \PP(K)$,
we obtain
\[
\frac{1}{n(n-1)} \sum_{i \neq j} \iint -\log |z_i \wedge z_j| 
     \, d\mu(z_i) d\mu(z_j) \ \geq \ D_n
\]
for all $n$, and therefore $c^0(K) \leq d^0_\infty(K)$.

\medskip

For the other direction,
choose $w_1,\ldots,w_n \in K$ such that
\[
n(n-1) D_n \ = \ \sum_{i \neq j} -\log |w_i \wedge w_j| \ .
\]

Without loss of generality, we may assume that $w_i \wedge w_j \neq 0$
for all $i \neq j$.
Define the measure $\nu_n$ to be the discrete measure on $K$ supported
equally on each of the points $w_i$, i.e.,
\[
\nu_n \ := \ \frac{1}{n} \sum_i \delta_{w_i}.
\]
By passing to a subsequence if necessary, we may assume that $\nu_n$
converges weakly to some probability measure $\nu$ on $K$.
Noting that
\[
\frac{n-1}{n} D_n \ = \ \frac{1}{n^2} \sum_{i \neq j} -\log |w_i \wedge w_j|
\ = \ \iint_{\CC^2 \times \CC^2 \backslash \text{(Diag)}}
          -\log |z \wedge w| \, d\nu_n(z) d\nu_n(w),
\]
it follows from (\ref{LimInfWedgeEquation}) and the definition of $c^0(K)$ that
\[
\lim_{n\to\infty} D_n \ \geq \ I(\nu) \ \geq \ -\log c^0(K),
\]
so that $d^0_\infty(K) \leq c^0(K)$ as desired.
\end{proof}

\medskip

We have introduced the homogeneous transfinite diameter by analogy
with the relation  between the classical transfinite diameter
and logarithmic capacity over $\CC$.
The integral defining the homogeneous capacity is difficult to extend
to nonarchimedean places.  However, the transfinite diameter generalizes
directly.

We will now give a formula for $d^0_{\infty}(K_{F,v})$ in terms of
resultants, motivated by the following result 
(\cite[Theorem~1.5]{DeMarco}):

\begin{theorem}[DeMarco]
\label{DeMarcoCapacityTheorem}
Suppose $F = (F_1,F_2) : \CC^2 \to \CC^2$ 
for some homogeneous polynomials $F_1,F_2 \in \CC[x,y]$ of degree $d$
with no common linear factor, 
and let $K_F$ be the filled Julia set of $F$ in $\CC^2$.
Then 
\begin{equation}
\label{DeMarcoFormula}
c^0(K_F) \ = \ |\Res(F)|^{-1/d(d-1)}.
\end{equation}
\end{theorem}

The proof given in \cite{DeMarco}, which involves both algebraic and
analytic ingredients, and does not carry over easily to the nonarchimedean setting.
One of our main results is the following adelic generalization of
DeMarco's theorem:

\begin{theorem}
\label{AdelicDeMarcoCapacityTheorem}
Suppose $F = (F_1,F_2) : \AA^2(\kbar) \to \AA^2(\kbar)$ 
for some homogeneous polynomials $F_1,F_2 \in k[x,y]$ of degree $d$
with no common linear factor over $\kbar$.  For $v \in M_k$, 
let $K_{F,v}$ be the filled Julia set of $F$ in $\CC_v^2$.
Then 
\begin{equation}
\label{AdelicDeMarcoFormula}
d^0_\infty(K_{F,v}) \ = \ |\Res(F)|_v^{-1/d(d-1)}.
\end{equation}
\end{theorem}

The proof, which is given in Section
\ref{HomogeneousSectionalCapacitySection},
requires the development of a considerable amount of capacity-theoretic
machinery.  It is completely independent of DeMarco's proof.
The product formula yields the following corollary, 
a key ingredient in our proof of Theorem~\ref{DynamicalEquidistributionTheorem}:

\begin{cor}
\label{CapacityProductFormula}
\[
\sum_{v \in M_k} \log d^0_\infty(K_{F,v}) \ = \ 0 \ .
\]
\end{cor}


\begin{remark}
When $\varphi(z)$ is a polynomial of degree $d$ with leading
coefficient $a_d$, and $K_{\varphi,v} \subseteq \CC_v$ is the $v$-adic
filled Julia set of $\varphi$ (as defined in \cite{BakerHsia}),
formula (\ref{AdelicDeMarcoFormula})
specializes to the formula
\[
c(K_{\varphi,v}) \ = \ |a_d|_v^{-1/(d-1)},
\]
which was first proved in \cite{BakerHsia}.
\end{remark}



\subsection{The Arakelov-Green's function and Local heights.}

In this subsection we will construct a
two-variable Green's function $g_{\varphi,v}(z,w)$ for the dynamical system associated
to $\varphi$.  It arises as a function on $\CC_v^2$ which
is invariant under scaling, and therefore descends to a function on $\PP^1(\CC_v)$.
We will see that the descended function gives a continuously varying
one-parameter family (indexed by $w \in \PP^1(\CC_v)$) of Call-Silverman local height
functions.

\smallskip
For notational convenience, write $c_v(F) := |\Res(F)|_v^{-1/d(d-1)}$.

If $v \in M_k$ and $z, w \in \CC_v^2$ are linearly independent over
$\CC_v$, define
\begin{equation}
\label{GreenDefinition}
G_{F,v}(z,w) \ := \ -\log |z\wedge w|_v + \Hhat_{F,v}(z) + \Hhat_{F,v}(w)
+ \log c_v(F).
\end{equation}
Recall that in the archimedean case,
$\SU(2) = \{\theta \in SL(2,\CC) : {^t\overline{\theta}} \theta = 1 \}$ is the
group preserving the both the norm $\|z\|_v$ on $\CC^2$ and the
alternating product $z \wedge w$, while in the nonarchimedean
case if $\hat{\cO}_v$ denotes the ring of integers of $\CC_v$, then
$SL(2,\hcO_v)$ is the group preserving $\|z\|_v$ and
$z \wedge w$ on $\CC_v^2$.
Write $\theta(F) = \theta \circ F \circ \theta^{-1}$.
  
We note the following properties of $G_{F,v}$.

\begin{lemma} { \ \ \ }

$A)$ $G_{F,v}$ is {\em doubly scale-invariant}, in the sense
that if $\alpha, \beta \in \CC_v^*$, then
\[
G_{F,v}(\alpha z, \beta w) \ = \ G_{F,v}(z,w).
\]

$B)$ For $\gamma \in \CC_v^*$, we have
\[
G_{\gamma F,v}(z, w) \ = \ G_{F,v}(z,w).
\]

$C)$ If $v$ is archimedean, then for each $\theta \in \SU(2)$,
\[
G_{\theta(F),v}(\theta(z),\theta(w)) \ = \ G_{F,v}(z,w) \ .
\]

\quad If $v$ is nonarchimedean, then for each $\theta \in SL(2,\hcO_v)$
\[
G_{\theta(F),v}(\theta(z),\theta(w)) \ = \ G_{F,v}(z,w) \ .
\]
\end{lemma}

\begin{proof}
Part A) follows immediately from the fact that 
\[
\log |\alpha z \wedge
\beta w|_v = \log |z \wedge w|_v + \log |\alpha|_v + \log |\beta|_v
\]
and from the fact that $\Hhat_{F,v}$ scales logarithmically.

Part B) follows from the following two easily verified facts:
\begin{itemize}
\item[(a)] $\Hhat_{\gamma F,v}(z) = \Hhat_{F,v}(z) + \frac{1}{d-1}
  \log |\alpha|_v.$
\item[(b)] $|\Res(\gamma F)|_v = |\Res(F)|_v |\gamma|_v^{2d}$.
\end{itemize}

Part C) follows from the fact that the group $\SU(2)$
(resp. $SL(2,\hcO_v)$) preserves $z \wedge w$, $\|z\|_v$, and $|\Res(F)|_v$.
To see that $\Res(\theta(F)) = \Res(F)$, 
note first that  
manipulating the determinant defining $\Res(\theta \circ F)$ 
shows that $\Res(\theta \circ F) = \Res(F)$.  
On the other hand, 
if  $F_1(z) = \prod_{i=1}^d (z \wedge a_i)$
and $F_2(z) = \prod_{j=1}^d (z \wedge b_j)$, then
$\Res(F) = \pm \prod_{i,j} (a_i \wedge b_j)$.  A simple computation 
shows that $\theta^{-1}(z) \wedge a_i = z \wedge \theta(a_i)$
and $\theta^{-1}(z) \wedge b_j = z \wedge \theta(b_j)$.
Since $\theta(a_i) \wedge \theta(b_j) = a_i \wedge b_j$, it follows that
$\Res(F \circ \theta^{-1}) = \Res(F)$.
\end{proof}

\medskip

In particular, $G_{F,v}$ descends to a well-defined
function $g_{\varphi,v}(z,w)$ on $\PP^1(\CC_v)$:
for any $z, w \in \PP^1(\CC_v)$ and any lifts $\tz, \tw \in \CC_v^2$ 
\begin{equation} \label{gFvDef} 
g_{\varphi,v}(z,w) \ = \ -\log |\tz\wedge \tw|_v + \Hhat_{F,v}(\tz) 
                      + \Hhat_{F,v}(\tw) + \log c_v(F).
\end{equation}
If $z\neq w$ then the right-hand side of (\ref{gFvDef}) is finite; if
$z=w$ then we define $g_{\varphi,v}(z,z) := +\infty$.

We will now establish another fact needed for the proof
of Theorem~\ref{DynamicalEquidistributionTheorem}.  Define
\begin{equation} \label{GammaDef} 
\Gamma_{\varphi,v} \ = \
\liminf_{n \to \infty} \inf_{z_1,\ldots,z_n \in \PP^1(\CC_v)}
\frac{1}{n(n-1)} \sum_{i \neq j} g_{\varphi,v}(z_i,z_j).
\end{equation}



\begin{lemma}[Positivity] 
\label{DynamicalPositivityLemma}  
For each $v\in M_k$, we have $\Gamma_{\varphi,v} \geq 0$.
\end{lemma}

\begin{proof}
Let $\epsilon > 0$.
Choose the lifts of the points $z_i$ in the statement of the Lemma so that
\[
- \epsilon \ < \ \Hhat_{F,v}(\tilde{z_i}) \ \leq \ 0
\]
for all $i$.  This is possible because $\Hhat_{F,v}$ scales
logarithmically and the set $\{ \log |\alpha|_v \; : \; \alpha \in \CC_v^*\}$
is dense in $\RR$.

In particular,  $\tilde{z_i} \in
K_{F,v}$ for all $i$.  By the definition of the homogeneous
transfinite diameter,
\[
\liminf_{n \to \infty} \inf_{z_1,\ldots,z_n \in \PP^1(\CC_v)}
\frac{1}{n(n-1)} \sum_{i \neq j} - \log |\tilde{z}_i \wedge
\tilde{z}_j|_v \ \geq \ - \log d^0_\infty (K_{F,v}) \ .
\]
But $d^0_\infty(K_{F,v}) = c_v(F)$ by
Theorem~\ref{AdelicDeMarcoCapacityTheorem}.
Therefore we obtain the inequality
\[
\liminf_{n \to \infty} \inf_{z_1,\ldots,z_n \in \PP^1(\CC_v)} 
\frac{1}{n(n-1)} \sum_{i \neq j} g_{\varphi,v}(z_i,z_j) \ > \ -2\epsilon\ .
\]
Since $\epsilon > 0$ was arbitrary, this gives the desired result.
\end{proof}

\begin{remark}  Later, in Corollary \ref{NormConstantCor}, we will see
via a global argument that $\Gamma_{\varphi,v} = 0$ for each $v$. 
\end{remark} 

Next we will show that $g_{\varphi,v}(z,w)$
forms a one-parameter family of Call-Silverman local heights
(see \cite{CS},\cite{CG}).  Recall that a function
$\hhat_{\varphi,v,D} : \PP^1(\CC_v) \setminus {\rm supp}(D) \to \RR$
is called a Call-Silverman {\em canonical local height function}
for $\varphi$, relative to the divisor $D$, if
it is a Weil local height associated to $D$, 
and if there exists
a rational function $f$ on $\PP^1$ over $\CC_v$ with 
${\rm div}(f) =  \varphi^*D - d \cdot D$ such that
\begin{equation*}
\hhat_{\varphi,v,D}(\varphi(z)) \ = \ 
d \cdot \hhat_{\varphi,v,D}(z) - \log |f|_v
\end{equation*}
for all 
$z \in \PP^1(\CC_v) \setminus ({\rm supp}(D) \cup \supp(\varphi^{*}(D)))$.
It is proved in \cite{CS} that a canonical local height function
exists for every divisor $D$, and is unique up to an additive constant.

\medskip
Choose coordinates on $\PP^1(\CC_v)$ in such a way that $\infty$
corresponds to the point $[0:1]$ and $0$ corresponds to $[1:0]$.
Let $V_{\infty} = \PP^1(\CC_v) \setminus \{ \infty \}$,
so that every $z=(z_0:z_1) \in V_{\infty}$
can be expressed uniquely as $[1:T(z)]$ with $T(z) = z_1/z_0 \in \CC_v$.

Define
$\hhat_{F,v,(\infty)} : V_\infty \to \RR$ by
\[
\hhat_{F,v,(\infty)}(z) \ = \ \Hhat_{F,v}(1,T(z)),
\]
and note that for $z \in V_{\infty}$, we have $(1,T(z)) \wedge (0,1) = 1$, so
\begin{equation}
\label{CoordinateGreen1}
g_{\varphi,v}(z,w) = \left\{
\begin{array}{ll}
-\log |T(z)-T(w)| + \hhat_{F,v,(\infty)}(z)
        + \hhat_{F,v,(\infty)}(w) + \log c_v(F) & w \neq \infty \\
\hhat_{F,v,(\infty)}(z) + \Hhat_{F,v}((0,1)) + \log c_v(F) & w = \infty. \\
\end{array}
\right.
\end{equation} 

For
$z \in \PP^1(\CC_v) \setminus ( \{ \infty \} \cup \varphi^{-1}(\infty))$,
we have the identity
\begin{equation}
\label{IterateIdentity}
F^{(n-1)}(1,T(\varphi(z))) \ = \ F^{(n-1)}(1,\frac{F_2(1,T(z))}{F_1(1,T(z))})
\ = \ \frac{F^{(n)}(1,T(z))}{F_1(1,T(z))^{d^{n-1}}}.
\end{equation}
Taking logarithms in (\ref{IterateIdentity}) and letting $n\to\infty$
gives the functional equation
\begin{equation}
\label{FunctionalEquation1}
\hhat_{F,v,(\infty)}(\varphi(z)) \ = \
d \cdot \hhat_{F,v,(\infty)}(z) - \log |F_1(1,T(z))|_v,
\end{equation}
which is valid for all
$z \in \PP^1(\CC_v) \setminus ( \{ \infty \} \cup \varphi^{-1}(\infty))$.

By Remark~\ref{BoundedDifferenceRemark}, there exists a constant
$C>0$ such that
\begin{equation}
\label{BoundednessEquation1}
|\hhat_{F,v,(\infty)}(z) - \log \max (1, |T(z)|_v) | \ \leq \ C
\end{equation}
for all $z \in V_{\infty}$.

Equations (\ref{FunctionalEquation1}) and (\ref{BoundednessEquation1})
show that $\hhat_{F,v,(\infty)}$ is a Call-Silverman
canonical local height
function on $\PP^1(\CC_v)$ for $\varphi$ relative to the divisor
$D = (\infty)$.   By (\ref{CoordinateGreen1}), the function
$g_{\varphi,v}(z,\infty)$
is also a Call-Silverman canonical local height for $\varphi$
relative to $(\infty)$, since $g(z,\infty)$ and
$\hhat_{F,v,(\infty)}(z)$ differ by a constant.
More generally, using the fact that $\SU(2)$ (resp. $\SL(2,\hcO_v)$)
acts transitively on $\PP^1(\CC_v)$,
it follows that for each $w \in \PP^1(\CC_v)$,
the function $g_{\varphi,v}(z,w)$
is a Call-Silverman local height for $\varphi$ relative to the
divisor $D = (w)$.  Thus $g_{\varphi,v}(z,w)$ is a continuously varying 
one-parameter family of local heights.

\smallskip

As a concrete example, take $w=0$.  In the affine patch
$V_0 = \PP^1(\CC_v) \setminus \{ 0 \}$, every point $z \in V_2$ can be
represented uniquely as $[U(z):1]$
with $U(z) = z_0/z_1 \in \CC_v$.

Define
$\hhat_{F,v,(0)}(z) : V_0 \to \RR$ by
\[
\hhat_{F,v,(0)}(z) \ = \ \Hhat_{F,v}(U(z),1)\ ,
\]
so that for $z \in V_0$, we have
\begin{equation}
\label{CoordinateGreen2}
g_{\varphi,v}(z,w) = \left\{
\begin{array}{ll}
-\log |U(z)-U(w)| + \hhat_{F,v,(0)}(z) + \hhat_{F,v,(0)}(w) + \log c_v(F)
                   & w \neq 0 \\
\hhat_{F,v,(0)}(z) + \Hhat_{F,v}((0,1)) + \log c_v(F) & w = 0. \\
\end{array}
\right.
\end{equation}
Then for
$z \in \PP^1(\CC_v) \setminus ( \{ 0 \} \cup \varphi^{-1}(0))$,
the identity
\begin{equation*}
F^{(n-1)}(U(\varphi(z)),1) \ = \ F^{(n-1)}(\frac{F_1(U(z),1)}{F_2(U(z),1)},1)
\ = \ \frac{F^{(n)}(U(z),1)}{F_2(U(z),1)^{d^{n-1}}},
\end{equation*}
gives the functional equation
\begin{equation}
\label{FunctionalEquation2}
\hhat_{F,v,(0)}(\varphi(z))
\ = \ d \cdot \hhat_{F,v,(0)}(z) - \log |F_2(U(z),1)|_v,
\end{equation}
valid for all
$z \in \PP^1(\CC_v) \setminus ( \{ 0 \} \cup \varphi^{-1}(0))$.
We also have
\begin{equation*}
|\hhat_{F,v,(0)}(z) - \log \max (|U(z)|_v,1) | \ \leq \ C
\end{equation*}
for all $z \in V_0$.

Finally, note that letting $n$ tend to infinity in the identity
\[
F^{(n)}(U(z),1) \ = \ F^{(n)}(1,T(z)) / T(z)^{d^n}
\]
and taking logarithms gives 
\begin{equation}
\label{HeightConversionEquation}
\hhat_{F,v,(0)}(z) \ = \ \hhat_{F,v,(\infty)}(z) - \log |T(z)|_v
\end{equation}
for all $z \in V_0 \cap V_{\infty} = \PP^1(\CC_v) \setminus \{ 0, \infty \}$.

\medskip

\subsection{Arakelov Green's functions and the canonical measure.} { \ \ \ }

In this subsection, we will show that 
$\frac{1}{\log(q_v)} g_{\varphi,v}(z,w)$ is in fact an Arakelov Green's function.  
This means showing that for each $w$,
the Laplacian of $g_{\varphi,v}(z,w)$ satisfies
\begin{equation*}
\frac{1}{\log(q_v)} \Delta(g_{\varphi,v}(z,w)) 
        \ = \ \delta_w(z) - \mu_{\varphi,v}(z), 
\end{equation*}
where $\mu_{\varphi,v}$ is a probability measure, independent of $w$.  
As will be explained below, in the nonarchimedean case the Laplacian 
is taken on the Berkovich space $\PP^1_{\Berk,v}$.
In the archimedean case,  $\mu_{\varphi,v}$
turns out to be the canonical measure supported on the Julia
set of $\varphi$ (see \S\ref{RationalFunctionSection}).  Thus   
``the canonical measure is the minus Laplacian of the local height''.

The measures $\mu_{\varphi,v}$ play a central role in our theory:  
they are the target measures in our main equidistribution theorem,
Theorem \ref{DynamicalEquidistributionTheorem}.  

\medskip 

In the archimedean case, for any Riemann surface $X/\CC$,
we define an Arakelov Green's function to be a function
$g(z,w) : X(\CC) \times X(\CC) \to \RR \cup \{\infty\}$
which satisfies the following two conditions:

\begin{itemize}
\item[(RS1)] (Continuity) The function $g(z,w)$
is a continuous as a function from $X(\CC) \times X(\CC)$ 
to the extended reals, and is finite off the diagonal.

\item[(RS2)] (Differential equation)
There is a probability measure $\mu$ on $X(\CC)$
such that for each fixed $w$,
$g(z,w)$ satisfies the distributional identity
\[
\Delta_z g(z,w) \ = \ \delta_w(z) - \mu(z).
\]
\end{itemize}

Conditions (RS1) and (RS2) imply that $g(z,w)$ is symmetric 
and bounded below, with a logarithmic singularity along the diagonal.  
These two conditions determine the
function $g(z,w)$ up to an additive constant.  There is a canonical
way to normalize it:  if 
\begin{itemize}
\item[(RS3)] (Normalization) \quad $\iint g(z,w) \, d\mu(z) d\mu(w) = 0$,
\end{itemize}
we will say $g(z,w)$ is a {\em normalized} Arakelov Green's function.
In any case, a non-normalized Arakelov Green's function still satisfies
\begin{equation}
\label{eq:LF'}
\iint g(z,w) \, d\mu(z) d\mu(w) \ < \ \infty \ .
\end{equation}
As noted in \cite{CR} (see also Lemma~\ref{LaplacianCalculationLemma}
below), if $g(z,w)$ satisfies (RS1) and (RS2), then
differentiating $\lambda(z) := \int g(z,w) \, d\mu(w)$ under the
integral sign shows that $\lambda(z)$ is harmonic everywhere and therefore
constant.  Thus (RS3) is equivalent to the following apparently 
stronger condition:
\begin{itemize}
\item[${\rm (RS3)}^{\prime}$] 
(Strong Normalization) \quad $\int g(z,w) \, d\mu(w) \equiv 0$ \ .
\end{itemize}

\begin{remark}
\label{LaplacianRemark}
This definition of an Arakelov Green's function, 
taken from \cite{RumelyBook} (see also \cite{R3} and \cite{Mai}),
is slightly looser than the one commonly used
in the literature, where the measure $\mu$ is required to be a smooth
positive $(1,1)$-form $\omega$ with total mass $1$.
\end{remark}

The operator $\Delta = -dd^c$ on $X(\CC)$ is to be considered
in the distributional sense.  In local coordinates, if $f$ is $\cC^2$
then in terms of the standard real Laplacian we have
\[
\Delta f = -\frac{1}{2\pi} \left( \frac{\partial^2 f}{\partial x^2} +
  \frac{\partial^2 f}{\partial y^2} \right)  dx \wedge dy \ .
\]


The fact that the distributional Laplacian of $g(z,w)$ is a negative
measure on $X(\CC) \backslash \{ w \}$
means that the restriction of $g(z,w)$ to $X(\CC) \backslash \{w\}$
is subharmonic.  
It should also be noted that the continuity of $g(z,w)$
imposes conditions on the measure $\mu$;
in particular, $\mu$ cannot have any point masses.
This follows from the Riesz Decomposition theorem 
(\cite{Ts}, Theorem II.24, p.45): on any local coordinate patch
$V \subset X(\CC) \backslash \{w\}$, there is a harmonic function $h_V(z)$
such that for all $z \in V$
\[
g(z,w) \ = \ h_V(z) + \int_V \log |z-x| \, d\mu(x) 
\] 
If $\mu$ had a point mass at some $p \in V$, then we would have 
$g(p,w) = -\infty$, a contradiction.  

The continuity of $g(z,w)$ also shows that 
$u_V(z) = \int_V \log |z-x| \, d\mu(x)$ 
is a continuous function of $z$ on $V$.
Therefore $\mu$ must be log-continuous, in the following sense:

\begin{defn} \label{LogContDef}
A bounded Borel measure $\mu$ is {\em log-continuous} if for
each $p \in X(\CC)$ there is a neighborhood $V$ of $p$ such that 
\[
u_V(z) \ := \ \int_V \log(|z-x|) \, d\mu(x) 
\]
is continuous on $V$.   
\end{defn}

In \cite{RumelyBook} this concept was called {\em log-finiteness}.  
However, the terminology log-continuous used here seems more appropriate.  

\medskip

In the nonarchimedean case, there is also a notion of an Arakelov
Green's function.  
In theory, one could define Arakelov Green's functions 
on an arbitrary Berkovich curve over $\CC_v$, 
but we restrict ourselves here to the case $X = \PP^1_{\Berk,v}$.
We only sketch the basic framework; for further details, 
see \cite[\S 4 --\S 7]{RumelyNotes}.

There is a class of functions
on $\PP^1_{\Berk,v}$, called {\em functions of bounded differential variation},
for which it is possible to define a measure-valued Laplacian;
this class is denoted $\BDV(\PP^1_{\Berk,v})$
(see \cite{RumelyNotes}, \S5.3).  The Laplacian is defined
first for functions on finitely branched subgraphs of 
$\PP^1_{\Berk,v} \backslash \PP^1(\CC_v)$ via the construction in 
(\cite{BR}, \S4) 
which generalizes the approaches of \cite{CR} and \cite{ZhangAP}.
It is then extended by a limiting process to functions on open 
subdomains of $\PP^1_{\Berk,v}$, 
using the Riesz Representation theorem.  
There are analogues of harmonic functions and subharmonic functions 
on $\PP^1_{\Berk,v}$.  Harmonic functions satisfy a maximum principle
(\cite{RumelyNotes}, Proposition 5.14), a Poisson formula 
(\cite{RumelyNotes}, Proposition 5.18), and Harnack's principle
(\cite{RumelyNotes}, Proposition 5.24).  Subharmonic functions are 
functions which locally belong to $\BDV(\PP^1_{\Berk,v})$ and have 
non-negative Laplacian (\cite{RumelyNotes}, Proposition 6.1).  They
have stability properties similar to classical subharmonic functions
(\cite{RumelyNotes}, Proposition 6.11) and  
satisfy a maximum principle (\cite{RumelyNotes}, Proposition 6.15), 
a comparison theorem (\cite{RumelyNotes}, Proposition 6.16), and a 
Riesz Decomposition Theorem (\cite{RumelyNotes}, Proposition 6.19).  
The pullback of a subharmonic function by a rational map is subharmonic
(\cite{RumelyNotes}, Proposition 7.13).  
In brief, \cite{RumelyNotes} provides all the tools necessary 
carry through arguments of 
classical potential theory on $\PP^1_{\Berk,v}$.  

For a rational function $\varphi$ acting on $\PP^1_{\Berk,v}$,
there is a theory of multiplicities at points of $\PP^1_{\Berk,v}$,
extending the usual algebraic multiplicities on $\PP^1(\CC_v)$
(see \cite{RumelyNotes}, Proposition 7.2).  
Given a Borel measure $\mu$ on $\PP^1_{\Berk,v}$,
this makes it possible to define pushforward and pullback measures
$\varphi_* \mu$ and $\varphi^* \mu$ with the usual formal properties
(\cite{RumelyNotes}, \S 7.3).

A (Berkovich) Arakelov Green's function is a function
$g(z,w) : \PP^1_{\Berk,v} \times \PP^1_{\Berk,v} 
\rightarrow \RR \cup \{\infty\}$
such that

\begin{itemize}
\item[(B1)] (Semicontinuity)
The function $g(z,w)$ 
is finite and continuous off the diagonal, 
and is strongly lower-semicontinuous on the diagonal, in the sense that 
for each $z \in \PP^1_{\Berk,v}$
\[
g(z,z) \ = \ \liminf \begin{Sb} (x,y) \rightarrow (z,z) \\ x \ne y \end{Sb}
                    g(x,y)  \ .
\]

\item[(B2)] (Differential equation)  For each $w \in \PP^1_{\Berk}$,
$g(z,w)$ belongs to $\BDV(\PP^1_{\Berk,v})$.  Furthermore,
there is a probability measure $\mu$ on $\PP^1_{\Berk,v}$
such that for each $w$, $g(z,w)$ satisfies the identity
\[
\Delta_z g(z,w) \ = \ \delta_w(z) - \mu(z).
\]
\end{itemize}

\medskip

As in the archimedean case, conditions (B1) and (B2) 
imply that $g(z,w)$ is symmetric and bounded below.  
The semicontinuity along the diagonal
is a technical condition which arises naturally from properties of the space 
$\PP^1_{\Berk,v}$ (see \cite{RumelyNotes}, Proposition 3.1).  
Together, (B1) and (B2) determine 
$g(z,w)$ up to an additive constant by the maximum principle
(\cite{RumelyNotes}, Proposition 5.14).  If in addition 
\begin{itemize}
\item[(B3)] (Normalization) \quad $\iint g(z,w) \, d\mu(z) d\mu(w) \ = \ 0$,
\end{itemize}
we will say $g(z,w)$
is a {\em normalized} Berkovich Arakelov Green's function.

Again, our assumption that $g(z,w)$ is continuous off the diagonal 
means that $\mu$ is log-continuous (the precise definition, and proof,
are given in \cite{RumelyNotes}, Proposition 7.15).  
And as in the archimedean case, 
log-continuity implies that $\mu$ has no point masses on $\PP^1(\CC_v)$. 
However, it can have point masses on $\PP^1_{\Berk,v} \backslash
\PP^1(\CC_v)$ (see Example~\ref{GoodReductionMeas} below).

\medskip
The function $g_{\varphi,v}(z,w)$ has a natural extension
`by continuity' to $\PP^1_{\Berk,v}$;   for details, see
(\cite{RumelyNotes}, \S 7.5).  We will write $g_{\varphi,v}(z,w)$ for both
the function on $\PP^1(\CC_v)$ constructed above, and its extension
to $\PP^1_{\Berk,v}$.

\medskip 
Recall that $q_v$ is the order of the residue field of $k_v$.  
We will now show that for each $v$, the function 
$\frac{1}{\log(q_v)} g_{\varphi,v}(z,w)$
is an Arakelov Green's function.  
The probability measure $\mu_{\varphi,v}$ associated to
$g_{\varphi,v}(z,w)$ (i.e., the measure occurring in (RS2) or (B2))
plays a key role in our theory.
In the nonarchimedean case, $\mu_{\varphi,v}$ 
is a measure on $\PP^1_{\Berk,v}$;  
in the archimedean case, it is a measure on $\PP^1(\CC)$.  
As noted earlier, $\PP^1_{\Berk}/\CC \cong \PP^1(\CC)$, 
so in fact we can view
$\mu_{\varphi,v}$ as a measure on $\PP^1_{\Berk,v}$ for all $v$.

\begin{prop}
\label{LyubichMeasureProp} { \ \ \ }

$A)$  For each $v \in M_k$, the function
$\frac{1}{\log(q_v)} g_{\varphi,v}(z,w)$ is an Arakelov Green's function associated to
a log-continuous probability measure $\mu_{\varphi,v}$ on $\PP^1_{\Berk,v}$.
For each $w \in \PP^1(\CC_v)$, the measure $\mu_{\varphi,v}$ 
is given locally on $V_w := \PP^1_{\Berk,v} \backslash \{w\}$ by
\begin{equation}
\label{DefinitionOfMuPhi}
\mu_{\varphi,v} |_{V_w} \ = \
       -\frac{1}{\log(q_v)} \Delta g_{\varphi,v}(z,w).
\end{equation}
Furthermore, $\varphi^* \mu_{\varphi,v} = d \cdot \mu_{\varphi,v}$ and
$\varphi_* \mu_{\varphi,v} = \mu_{\varphi,v}$.

$B)$ If $v \in M_k$ is archimedean, then $\mu_{\varphi,v}$ coincides 
with the canonical measure on $\PP^1(\CC)$ associated 
to $\varphi$ by Lyubich and Freire-Lopes-Ma{\~n}{\'e}.  

\end{prop}

\begin{proof}
For $v$ nonarchimedean, this is \cite[Theorem 7.14]{RumelyNotes}.
Henceforth assume $v$ is archimedean. For part A), 
note first that both $\hhat_{F,v,(0)}(z)$ and
$\hhat_{F,v,(\infty)}(z)$ are uniform limits of subharmonic
functions, and are therefore subharmonic.
Thus both
$-\Delta \hhat_{F,v,(0)}(z)$ and
$-\Delta \hhat_{F,v,(\infty)}(z)$ are nonnegative measures.
Using the relation (\ref{HeightConversionEquation}), we see that 
$\hhat_{F,v,(0)}(z)$ and
$\hhat_{F,v,(\infty)}(z)$ differ by a harmonic function on $V_0 \cap
V_\infty$, and therefore
$\Delta \hhat_{F,v,(\infty)}(z) = \Delta \hhat_{F,v,(0)}(z)$ on
$V_0 \cap V_\infty$.
It follows that there is a non-negative measure
$\mu_{\varphi,v}$ given locally by (\ref{DefinitionOfMuPhi}).
By (\ref{CoordinateGreen1}) and (\ref{CoordinateGreen2}), for any $w$
we have $\frac{1}{\log(q_v)} \Delta_z g_{\varphi,v}(z,w) = \delta_w(z) - \mu_{\varphi,v}$
on both $V_0$ and $V_\infty$,
and hence on $V_0 \cup V_\infty = \PP^1(\CC)$, as desired.

The fact that $\mu_{\varphi,v}$ is a {\em probability} measure (i.e., that
$\mu_{\varphi,v}(\PP^1(\CC)) = 1$) follows immediately from the identity
$\frac{1}{\log(q_v)} \Delta_z g_{\varphi,v}(z,w) = \delta_w(z) - \mu_{\varphi,v}$,
since the distributional Laplacian of a function on $\PP^1(\CC)$
always has total mass zero.

To see that $\varphi^*(\mu_{\varphi,v}) = d \cdot \mu_{\varphi,v}$,
combine (\ref{FunctionalEquation1}) and (\ref{FunctionalEquation2}),
using the fact that $F_1$ and $F_2$ have no common zeros in
$\CC^2$ by assumption.  Finally, the relation
$\varphi_*(\mu_{\varphi,v}) = \mu_{\varphi,v}$ follows formally from
$\varphi^*(\mu_{\varphi,v}) = d \cdot \mu_{\varphi,v}$ using the
fact that $\varphi_*(\varphi^*(\mu)) = d \cdot \mu$
for all measures $\mu$ on $\PP^1(\CC)$.

For part B), recall from Theorem~\ref{LyubichTheorem}
that the canonical measure is the unique probability
measure $\mu$ on $\PP^1(\CC)$ with no point masses such that
$\varphi^*(\mu) = d \cdot \mu$.
As noted above, the continuity of $g_{\varphi,v}(z,w)$ off the
diagonal implies that $\mu_{\varphi,v}$ has no point masses, and the
functional equation $\varphi^*(\mu_{\varphi,v}) = d \cdot \mu_{\varphi,v}$
has been established in A).  Hence $\mu_{\varphi,v}$ coincides with
the canonical measure.
\end{proof}


\begin{remark}
\label{NormalizedGreenRemark}
We will see in Corollary \ref{NormConstantCor}, as a consequence
of global considerations, that $\frac{1}{\log(q_v)} g_{\varphi,v}(z,w)$ is in fact a 
{\em normalized} Arakelov Green's function for each $v$.
\end{remark}


\medskip
Assuming Remark~\ref{NormalizedGreenRemark}, we can establish the following 
invariance property of $g_{\varphi,v}(z,w)$:

\begin{cor} \label{InvarianceCorollary}
Given $w \in \PP^1_{\Berk,v}$, write $\varphi^*((w)) = \sum_{i=1}^r m_i(w_i)$.  
Then for all $z \in \PP^1_{\Berk,v}$,
\[
g_{\varphi,v}(\varphi(z),w) \ = \ \sum_{i=1}^r m_i \, g_{\varphi,v}(z,w_i) \ .
\]
\end{cor}

\begin{proof}  We will only give the proof in the archimedean case;  
in the nonarchimedean case the proof is formally identical, 
using properties of the Berkovich Laplacian.

Since $\SU(2)$ acts transitively on $\PP^1(\CC)$, we can assume without loss
that $w = \infty$.  By formula (\ref{CoordinateGreen1}),
there is a constant $C_1$ such that
\begin{equation} \label{FAA1}
g_{\varphi,v}(\varphi(z),\infty) 
\ = \ \hhat_{F,v,(\infty)}(\varphi(z)) + C_1 \ .
\end{equation}
By the functional equation (\ref{FunctionalEquation1}) 
of the Call-Silverman local height $\hhat_{F,v,(\infty)}$,  
\begin{equation} \label{FAA2}
\hhat_{F,v,(\infty)}(\varphi(z)) \ = \
d \cdot \hhat_{F,v,(\infty)}(z) - \log |F_1(1,T(z))|_v \ .
\end{equation}
Here $f(z) := F(1,T(z))$ is a polynomial with divisor 
$\div(f) = \sum m_i(w_i) - d \cdot (\infty)$, 
where $\sum m_i(w_i) = \varphi^*(\infty)$.  We claim
there is a constant $C_2$ such that
\begin{equation} \label{FAA3}
-\log(|f(z)|_v) \ = \ 
\sum m_i g_{\varphi,v}(z,w_i) - d \cdot g_{\varphi,v}(z,\infty) + C_2.
\end{equation}
To see this, note that both sides have Laplacian equal to 
$\log(q_v)$ times $\sum m_i \delta_{w_i}(z) - d \cdot
  \delta_{\infty}(z)$;  hence their
difference is a function which is harmonic everywhere, thus constant.  
(In the nonarchimedean case this argument is justified by 
\cite{RumelyNotes}, Proposition 5.14 and Lemma 5.12).    
Combining (\ref{FAA1}), (\ref{FAA2}) and (\ref{FAA3}) shows that
\[
g_{\varphi,v}(\varphi(z),\infty) 
\ = \ \sum_{i=1}^r m_i g_{\varphi,v}(z,w_i) + C
\]
for some constant $C$.  Integrating the left-hand side against
$\mu_{\varphi,v}$ and using the invariance property of $\mu_{\varphi,v}$,
the fact that $\frac{1}{\log(q_v)} g_{\varphi,v}(z,w)$ is normalized, and ${\rm (RS3)}^{\prime}$,
we have 
\[
\int g_{\varphi,v}(\varphi(z),\infty) \, d\mu_{\varphi,v}(z)  
= \int g_{\varphi,v}(z,\infty) \, d(\varphi_*\mu_{\varphi,v})(z) 
= \int g_{\varphi,v}(z,\infty) \, d\mu_{\varphi,v}(z) =  0 \ .
\]
Computing the integral of the right-hand side, we get $C$.  
Therefore $C = 0$ as desired.
\end{proof}  

\medskip

\begin{example} \label{GoodReductionMeas}
Recall from \cite{CS} that $\varphi$ is said to have {\em good reduction} 
at a place $v$ if it can be written as $\varphi(T) = G_2(T)/G_1(T)$
where $G_1, G_2 \in \cO_v(T)$ are such that the reduced polynomials
$g_1 = G_1 \mod(m_v)$, $g_2 = G_2 \mod(m_v)$ are nonzero and coprime,
with $\max(\deg(g_1),\deg(g_2)) = d$.  

If $\varphi$ has good reduction at $v$, then by Example 7.2 of
\cite{RumelyNotes}, $\mu_{\varphi,v}$ is the discrete measure 
supported at the {\em Gauss point} $\zeta_0$ of $\PP^1_{\Berk,v}$, and 

\begin{equation}
\label{IntersectionEquation}
g_{\varphi,v}(z,w) = \left\{ \begin{array}{ll} -\log |z-w|_v + \log^+|z|_v +
   \log^+|w|_v & z,w \neq \infty \\
\log^+|z|_v & w=\infty \\ 
\log^+|w|_v & z=\infty. \\ 
\end{array} \right.
\end{equation}
\end{example}

\medskip

\subsection{The Energy Minimization Principle.}

If $g(z,w)$ is an Arakelov Green's function on $\PP^1_{\Berk,v}$
associated to the measure $\mu$, we will often write $g_{\mu}(z,w)$ 
instead of $g(z,w)$.  
With this notation, $\frac{1}{\log(q_v)} g_{\varphi,v}(z,w) = g_{\mu_{\varphi,v}}(z,w)$.

Arakelov Green's functions on $\PP^1_{\Berk,v}$
satisfy the following important energy minimization principle:

\begin{theorem} \label{EnergyMinimizationTheorem} 
Let $v$ be a place of $k$, and let $g_\mu(z,w)$ be
an Arakelov Green's function on $\PP^1_{\Berk,v}$ 
whose associated measure $\mu$ is log-continuous.  Define the
``energy functional'' $I_\mu(\nu)$ on the space $\PP$ of probability
measures on $\PP^1_{\Berk,v}$ by the formula
\[
I_\mu(\nu) \ := \ 
\iint_{\PP^1_{\Berk,v} \times \PP^1_{\Berk,v}} g_\mu(z,w) \, d\nu(z) d\nu(w).
\]
Then $I_\mu(\nu) \geq I_\mu(\mu)$ for all $\nu \in \PP$,
with equality if and only if $\nu = \mu$.
\end{theorem}
 
In the archimedean case, Theorem \ref{EnergyMinimizationTheorem} will be
proved in \S\ref{RiemannSurfaceSection}
as a consequence of the more general Theorem \ref{RiemannSurfaceTheorem};  in the nonarchimedean case, 
it is proved in \cite[Theorem 7.20]{RumelyNotes}. 

\medskip

\subsection{Discrete approximations to the energy integral.}

In this section, $v$ denotes an arbitrary place of $k$, and we work
on the Berkovich space $\PP^1_{\Berk,v}$.  Recall that if $v$
is archimedean, the space $\PP^1_{\Berk,v}$ is just $\PP^1(\CC)$.


The following lemma enables us to
apply the energy minimization principle in a useful way
to discrete measures.  We state it abstractly because it was 
also used in Lemma \ref{EqualityLemma}.

\begin{lemma}
\label{LimInfLemma}
Let $(X,\nu)$ be a measure space, with $\nu$ a probability measure.
Let $\{S_n\}_{n \ge 1}$ be a sequence of finite subsets of $X$,
and for each $n$ let $\delta_n$ be the discrete
probability measure supported equally at all elements of $S_n$.
Suppose the measures $\delta_n$ converge weakly to $\nu$.
Let $g : X \times X \to \RR \cup \{\infty\}$
be a function which is finite, continuous, and bounded from below on
$X \times X \backslash {\rm (Diag)}$.
Then 
\[
\liminf_{n\to\infty} 
\iint_{X \times X \backslash {\rm (Diag)}} g(z,w) \,
        d\delta_n(z) d\delta_n(w) \ \geq \
\iint_{X \times X} g(z,w) \, d\nu(z) d\nu(w).
\]
\end{lemma}

\begin{proof}

Define $N_n := \# S_n$.
For any fixed real number $M>0$, we have
\begin{equation}
\label{eqnsmalldiag}
\iint_{{\rm (Diag)}} \min \{ M, g(x,y) \} \, d\delta_n(x) d\delta_n(y) 
\ = \frac{1}{N_n} \cdot M 
\end{equation}
by the definition of $\delta_n$, and therefore
{\allowdisplaybreaks
\begin{equation*}
\begin{aligned}
& \liminf_{n\to\infty}
\iint_{X \times X \backslash {\rm (Diag)}} g(x,y) \,
                d\delta_n(x) d\delta_n(y) \notag\\
&\ \ge    \lim_{M\to \infty}
\liminf_{n\to\infty} \iint_{X \times X \backslash {\rm (Diag)}}
       \min\{M, g(x,y) \} \, d\delta_n(x) d\delta_n(y)
\qquad\text{(since $(*) \ge \min\{M, (*) \}$)}
\notag\\
&\ =  \lim_{M\to \infty}
\liminf_{n\to\infty} \iint_{X \times X }
         \min\{M, g(x,y) \} \, d\delta_n(x) d\delta_n(y)
\qquad\text{(by (\ref{eqnsmalldiag}))} 
\notag\\
&\ = \lim_{M\to \infty}  \iint_{X \times X}
        \min\{M, g(x,y) \} \, d\nu (x) d\nu(y)
\qquad\text{($\delta_n \to \nu$ weakly)} 
\notag\\
&\ = \iint_{X \times X} g(x,y) \, d\nu (x) d\nu(y)
\qquad\text{(monotone convergence theorem).} 
\notag\\
\end{aligned}
\end{equation*}
}
\end{proof}

\medskip

Let $g_{\mu}(z,w)$ be an Arakelov Green's function on $\PP^1_{\Berk,v}$ with 
associated log-continuous measure $\mu$.  
We now introduce a quantity $D_\infty(\mu)$ analogous to the 
(negative logarithm of the) classical transfinite diameter.   
For $n\geq 2$, define
\[
D_n(\mu) \ := \ 
\inf_{z_1,\ldots,z_n \in \PP^1(\CC_v)} \frac{1}{n(n-1)} \sum_{i \neq j} 
        g_{\mu}(z_i,z_j).
\]
Since $g_\mu(z,w)$ is bounded below and is 
finite off the diagonal, each $D_n(\mu)$ is a well-defined real number.  

The proof of the following lemma is similar to that of 
Lemma~\ref{HomogeneousMonotonicLemma}.
%

\begin{lemma}
\label{MonotonicLemma}
The sequence $D_n(\mu)$ is non-decreasing.
\end{lemma}

\begin{proof}
Take $n\geq 2$, 
fix $\varepsilon > 0$, and choose $w_1,\ldots,w_{n+1}$ such that 
\[
\sum_{i \neq j} g_{\mu}(w_i,w_j) \ \le \ n(n+1) (D_{n+1}(\mu) + \varepsilon) \ .
\]
By the definition of $D_n(\mu)$, we have (for each $1\leq m \leq n+1$)
\[
n(n-1) D_{n}(\mu) \ \leq \ \sum\begin{Sb} i,j \neq m \\ i \neq j \end{Sb}
g_{\mu}(w_i,w_j).
\]
Adding together these $n+1$ inequalities gives
\[
(n+1) n(n-1) D_n(\mu) \ \leq \ 
(n-1) \sum_{i \neq j} g(w_i,w_j) 
\ \leq \ (n-1)n(n+1) (D_{n+1}(\mu) + \varepsilon) \ .
\]
Since $\varepsilon > 0$ is arbitrary, $D_n(\mu) \leq D_{n+1}(\mu)$ as desired.
\end{proof}

\medskip

Define 
\[
D_\infty(\mu) \ = \ \lim_{n\to\infty} D_n(\mu) \ .
\]
The following
result is analogous to the equality of the transfinite diameter and
the capacity of a compact set in classical complex potential theory.

\begin{theorem}
\label{DiscreteToContinuousTheorem}
$D_\infty(\mu) = I_\mu(\mu).$
\end{theorem}

\begin{proof}
Let $n\geq 2$.  We first claim that for 
all $z_1,\ldots,z_n \in \PP^1_{\Berk,v}$, 
\begin{equation} \label{FBB2}
\frac{1}{n(n-1)} \sum_{i \neq j} g_\mu(z_i,z_j) \ \geq \ D_n(\mu) \ .
\end{equation}
In the archimedean case this is immediate, since 
$\PP^1_{\Berk,v} = \PP^1(\CC)$.  
To see it in the nonarchimedean case, first suppose the $z_i$ are distinct,
and note that by the continuity of $g_{\mu}(z,w)$ 
off the diagonal and the fact that $\PP^1(\CC_v)$ 
is dense in $\PP^1_{\Berk,v}$,
for any $\varepsilon > 0$ 
there are points $x_1, \ldots, x_n \in \PP^1(\CC_v)$ with 
\[
|g_{\mu}(z_i,z_j)-g_{\mu}(x_i,x_j)| \ < \ \varepsilon \ .
\]
By definition we have 
$\frac{1}{n(n-1)} \sum_{i \neq j} g_\mu(x_i,x_j) \geq D_n(\mu)$, 
so letting $\varepsilon \rightarrow 0$ gives (\ref{FBB2}).  The
general case follows by the strong lower semicontinuity of $g_{\mu}(z,w)$
(see axiom (B1) for Berkovich Arakelov Green's functions).   

Integrating (\ref{FBB2}) against $d\mu(z_1)\cdots d\mu(z_n)$, we
see that
\[
\frac{1}{n(n-1)} \sum_{i \neq j} \iint g_\mu(z_i,z_j) d\mu(z_i) d\mu(z_j)
\ \geq \ D_n(\mu),
\]
for all $n$, and therefore $I_\mu(\mu) \geq D_\infty(\mu)$.

\medskip

For the other direction, for each $n$ 
choose $w_1,\ldots,w_n \in \PP^1(\CC_v)$ such that 
\[
\frac{1}{n(n-1)} \sum_{i \neq j} g_\mu(w_i,w_j) 
    \ \le \ D_n(\mu) + \frac{1}{n}
\]
and let $\nu_n$ be the discrete measure supported
equally on each of the points $w_i$, i.e.,
\[
\nu_n := \frac{1}{n} \sum_i \delta_{w_i}.
\]

By passing to a subsequence if necessary, we may assume that the $\nu_n$
converge weakly to some measure $\nu$ on $\PP^1_{\Berk,v}$.
Noting that
\[
\frac{n-1}{n} (D_n(\mu) + \frac{1}{n}) 
\ \ge \ \frac{1}{n^2} \sum_{i \neq j} g_\mu(w_i,w_j) 
\ = \ \iint_{\PP^1_{\Berk,v} \times \PP^1_{\Berk,v} \backslash {\rm (Diag)}}
g_\mu(w,z) \, d\nu_n d\nu_n\ ,
\]
it follows from Theorem~\ref{EnergyMinimizationTheorem} 
and Lemma~\ref{LimInfLemma} that
\[
D_\infty(\mu) \ = \ \liminf_{n \rightarrow \infty} D_n(\mu)
\ \geq \ I_\mu(\nu) \ \geq \ I_\mu(\mu)
\]
as desired.
\end{proof}

\begin{remark}
By the exact same arguments, one sees that Lemma~\ref{MonotonicLemma}
and Theorem~\ref{DiscreteToContinuousTheorem}
remain true in the archimedean case for an arbitrary Riemann surface,
using Theorem~\ref{RiemannSurfaceTheorem} instead of 
Theorem~\ref{EnergyMinimizationTheorem}.
\end{remark}



\section{Proof of the main equidistribution theorem}
\label{EquidistributionSection}

\medskip

We now turn to the proof of Theorem \ref{DynamicalEquidistributionTheorem}.
As will be seen, the theorem follows rather formally from the machinery
developed above.
Before giving the argument, we deal with some technical preliminaries.

\subsection{Base change lemmas}

\medskip

In this subsection, we formulate a lemma which relates local
Arakelov Green's functions over different base fields.

\medskip

Let $\varphi$ be a rational function of degree $d \geq 2$ defined over
the number field $k$.  For $v \in M_k$,
define $g_{\varphi,v}(z,w) := G_{F,v}(z,w)$ for some lift $F$ of $\varphi$
to $k[x,y] \times k[x,y]$, i.e., given $z, w \in \PP^1(\kbar)$,
take lifts $\tz, \tw$ of $z$ and $w$ to $\kbar^2$; then
\[
g_{\varphi,v}(z,w) \ = \ -\log|\tz \wedge \tw|_v +
            \Hhat_{F,v}(\tz) + \Hhat_{F,v}(\tw) + \log c_v(F),
\]
where $c_v(F) = |\Res(F)|^{-\frac{1}{d(d-1)}}$ as before.
If $k^{\prime}/k$ is a finite extension, we can in a similar way define
$g_{\varphi,v^{\prime}}(z,w)$ for $v^{\prime} \in M_{k^{\prime}}$.  We have:

\begin{lemma}
\label{BaseChangeLemma}

$A)$ The expression 
\[
g_\varphi(z,w)
\ := \ \frac{1}{[k^{\prime}:\QQ]} \sum_{v^{\prime} \in M_{k^{\prime}}}
   g_{\varphi,v^{\prime}}(z,w)
\]
is independent of the choice of a number field $k^{\prime}$
containing $z$ and $w$, and therefore gives a well-defined 
function on $\kbar \times \kbar \backslash {\rm (Diag)}$.

$B)$ 
For all $z,w \in \kbar$, $z \neq w$, 
\[
g_\varphi(z,w) \ = \ \hhat_{\varphi}(z) + \hhat_{\varphi}(w).
\]

$C)$
Let $k^{\prime}$ be a finite extension of $k$.  Take $v \in M_k$,
and let $v^{\prime}$ be a place of $k^{\prime}$ with $v^{\prime} \mid v$.
If $S$ is a finite $\Gal(k^{\prime}/k)$-invariant subset of $k^{\prime}$,
then for all $z,w \in S$, $z \neq w$, the expression
$\sum_{z\neq w \in S} g_{v^{\prime}}(z,w)$
is independent of the place $v^{\prime}$, and
\[
\frac{1}{[k:\QQ]} \left( \sum_{z\neq w \in S} g_{\varphi,v}(z,w) \right)
\ = \ \frac{1}{[k^{\prime}:\QQ]} \sum_{v^{\prime} \mid v}
\left( \sum_{z\neq w \in S} g_{\varphi,v^{\prime}}(z,w) \right).
\]

$D)$
Let $z_1,\ldots,z_N$ be the Galois conjugates of an element
$z \in \PP^1(\QQbar) \backslash \PP^1(\QQ)$.  Then
\begin{equation}
\label{eq:BaseChange}
\frac{1}{[k:\QQ]} \sum_{v \in M_k}
\left( \frac{1}{N(N - 1)} \sum_{i \neq j} g_{\varphi,v}(z_i,z_j) \right)
\ = \ 2\hhat_\varphi(z).
\end{equation}

\end{lemma}

\begin{proof}
The proofs of A) and C) are straightforward consequences
of our choice of normalizations for absolute values.
B) follows from A) by the product formula (applied twice):
if $k^{\prime}/k$ is a finite extension
such that the lifts $\tz$ and $\tw$ are rational over $k^{\prime}$, then
\begin{eqnarray*}
g_\varphi(z,w)
& = & \frac{1}{[k^{\prime}:\QQ]} \sum_{v^{\prime} \in M_{k^{\prime}}}
   g_{\varphi,v^{\prime}}(z,w) \\
& = & \frac{1}{[k^{\prime}:\QQ]} \sum_{v^{\prime} \in M_{k^{\prime}}}
  \left( -\log|\tz \wedge \tw|_{v^\prime} + \Hhat_{F,v^{\prime}}(\tz)
      + \Hhat_{F,v^{\prime}}(\tw)  + \log c_{v^\prime}(F) \right) \\
& = & \hhat_{\varphi}(z) + \hhat_{\varphi}(w)
\end{eqnarray*}
since $\tz \wedge \tw = \tz_0 \tw_1 - \tz_1 \tw_0 \in (k^{\prime})^*$ and
$c_{v^{\prime}}(F) = |\Res(F)|_{v^\prime}^{-1/d(d-1)}$.
Finally, D) follows from C) by summing both sides over all places $v$ of $k$.
\end{proof}

\subsection{Lemmas on double sums.}

\medskip

Before turning to Theorem~\ref{DynamicalEquidistributionTheorem} and its
proof, we need two lemmas on doubly-indexed sums.
The first is a discrete analogue of Fatou's lemma from real analysis:

\begin{lemma}
\label{FatouLemma}
Suppose $a_n^{(j)}$ is a doubly-indexed sequence of real numbers 
which satisfy the following two properties:
\begin{itemize}
\item[(F1)] For each $n$, $\sum_j a_n^{(j)}$ converges.
\item[(F2)] There is a collection $\{ M_j \}$ of real numbers,
  almost all zero, such that $a_n^{(j)} \geq -M_j$
for all $j,n$.
\end{itemize}
Then
\begin{equation}
\label{FatouEquation}
\sum_j \liminf_{n \to\infty} a_n^{(j)} 
\ \leq \
\liminf_{n \to\infty} \sum_j a_n^{(j)}
\end{equation}
as extended real numbers.
\end{lemma}

\begin{proof}
Replacing $a_n^{(j)}$ by $a_n^{(j)} + M_j$ if necessary, we may assume
without loss of generality that $a_n^{(j)} \geq 0$ for all $n,j$.  The
result now follows immediately from the usual version of Fatou's lemma
(see \cite{Royden}, Theorem 4.3.9)
applied to the sequence $f_n$ of locally
constant functions defined by
$f_n(x) = a_n^{(j)} {\rm \; if \;} x \in [j,j+1)$.
\end{proof}

\medskip

The next lemma is a simple application of Lemma~\ref{FatouLemma}.

\begin{lemma}
\label{TechnicalLemma}
Suppose $a_n^{(j)}$ is a doubly-indexed sequence of real numbers which
satisfy properties $(F1)$ and $(F2)$.  Consider the following conditions,
where $L, L_j \in \RR$.
\begin{itemize}
\item[$(L1)$] \qquad
$\displaystyle{ \limsup_{n\to\infty} \sum_j a_n^{(j)} \ \leq \ L \ ,}$

\item[$(L2)$] \qquad 
$\displaystyle{ \sum_j \liminf_{n\to\infty} a_n^{(j)} \ \geq \ L \ ,}$

\item[$\ (L2)^{\prime}$] \qquad
$\displaystyle{\liminf_{n\to\infty} a_n^{(j)} \ \geq \ L_j \ .}$
\end{itemize} 

\noindent{Then:}

$A)$  If $(L1)$ and $(L2)$ hold for some $L$, then
$\displaystyle{\lim_{n\to\infty} a_n^{(j)}}$ exists for all $j$.

$B)$  If $(L1)$ holds, and there are numbers $L_j$ with
with $\sum_j L_j \geq L$ such that $(L2)^{\prime}$ holds for all $j$,
then $\displaystyle{\lim_{n\to\infty} a_n^{(j)} = L_j}$ for all $j$.
\end{lemma}

\begin{proof}
For any sequences $a_n,b_n$ of real numbers which are bounded from below,
it is easy to see that
$\limsup (a_n + b_n) \geq \limsup (a_n) + \liminf(b_n)$.

For any index $i$, one therefore sees from Lemma~\ref{FatouLemma} that
\begin{eqnarray*}
L &\geq& \limsup \sum_j a_n^{(j)} 
     \ \geq \ \limsup a_n^{(i)} + \liminf \sum_{j \neq i} a_n^{(j)} \\
&\geq& \limsup a_n^{(i)} + \sum_{j\neq i} \liminf a_n^{(j)} 
     \ \geq \ \limsup a_n^{(i)} - \liminf a_n^{(i)} + L \ ,
\end{eqnarray*}
which implies that 
$\limsup a_n^{(i)} \leq \liminf a_n^{(i)}$.
Therefore 
$\lim_{n \to \infty} a_n^{(i)}$ exists for all $i$, which proves A).

For B), note that
\begin{eqnarray*}
L &\geq& \limsup \sum_j a_n^{(j)} \ \geq \ \liminf \sum_j a_n^{(j)} \\
& \geq & \sum_j \liminf a_n^{(j)} \ \geq \ \sum L_j  \ \ge \ L.
\end{eqnarray*}
Therefore equality holds everywhere, so
$\liminf a_n^{(j)} = L_j$ for all $j$.  
By A), we conclude that
$\lim a_n^{(j)} = L_j$ for all $j$, which establishes B).
\end{proof}

\medskip

In applying Lemma~\ref{TechnicalLemma}, we will use the following
easily verified properties of the collection of functions $g_{\varphi,v}(z,w)$:
\begin{itemize}
\item[(G1)] For fixed $z,w \in \PP^1(\kbar)$ with $z \neq
  w$, we have $g_{\varphi,v}(z,w) = 0$ for almost all $v$.
\item[(G2)] For almost all $v$, we have
$g_{\varphi,v}(z,w)\geq 0$ for all $z,w \in \PP^1(\CC_v)$.
\end{itemize}

\medskip

\subsection{Pseudo-equidistribution.} 

{ \ \ \ }

Let $v$ be a place of $k$.
If $S$ is a finite subset of $\PP^1(\CC_v)$ of cardinality $N$, we define a 
discrete probability measure $\delta(S)$ on $\PP^1(\CC_v)$ by
\[
\delta(S) \ := \ \frac{1}{N} \sum_{z \in S} \delta_z\ .
\]
Note that if $S$ is a subset of $\PP^1(\kbar)$, we can consider $S$ as
a subset of $\PP^1(\CC_v)$ for each $v \in M_k$, 
since we have fixed an embedding $\kbar \into
\CC_v$ for each $v$.  If $S$ is $\Gal(\kbar/k)$-stable, the
resulting subset of $\PP^1(\CC_v)$ is independent of the choice of embedding.


\begin{defn} 
A sequence of finite subsets $\{ S_n \}_{n \ge 1}$ of $\PP^1(\CC_v)$ 
is {\em pseudo-equidistributed} with respect to $g_{\varphi,v}$ if
$N_n = \#(S_n) \to \infty$ and
\begin{equation}
\label{PseudoEquidistributionEquation}
\lim_{n\to\infty} \frac{1}{N_n(N_n-1)} \sum \begin{Sb} z,w \in S_n \\
  z\neq w \end{Sb} g_{\varphi,v}(z,w) \ = \ 0\ .
\end{equation}
\end{defn}
Recall that by Lemma \ref{DynamicalPositivityLemma} the minimal possible
value for the left-hand side of (\ref{PseudoEquidistributionEquation})
is $0$.  Thus, the sequence $\{S_n\}_{n \ge 1}$ is
pseudo-equidistributed if and only if it achieves this minimum value.

\begin{remark} 
This definition is a bit different from the definition of
pseudo-equidistribution in (\cite{BakerHsia}).  It anticipates the
fact, shown in Corollary \ref{NormConstantCor} below, 
that $\frac{1}{\log(q_v)} g_{\varphi,v}(z,w)$ is a normalized Arakelov Green's function.  
\end{remark}


We now prove the following adelic pseudo-equidistribution result:

\begin{theorem}
\label{PseudoEquidistributionTheorem}
Let $z_n$ be a sequence of distinct points of $\PP^1(\kbar)$ such that
$\hhat_\varphi(z_n) \to 0$.  Let $S_n$ denote
the set of Galois conjugates $($over $k$$)$ of $z_n$.  
Then the sequence $\{ S_n \}$ is pseudo-equidistributed with respect
to $g_{\varphi,v}$ for all $v \in M_k$.
\end{theorem}


%


\begin{proof}

Let $N_n$ be the cardinality of $S_n$.  By Northcott's finiteness theorem, 
the hypothesis $\hhat_\varphi(z_n) \to 0$ (and the fact that the points
$z_n$ are all distinct) implies that $N_n \to \infty$ as $n\to \infty$.  

For $v \in M_k$ and $n\geq 1$, set
\[
g_{v,n} \ := \ \frac{1}{N_n(N_n - 1)} \sum \begin{Sb} z,w \in S_n \\
  z\neq w \end{Sb} g_{\varphi,v}(z,w)\ .
\]
By (\ref{eq:BaseChange}), we have
\begin{equation}
\label{gnEstimate}
g_n := \frac{1}{[k:\QQ]} \sum_{v \in M_k} g_{v,n}
= 
\frac{1}{N_n(N_n - 1)} \frac{1}{[k:\QQ]} 
\sum \begin{Sb} z,w \in S_n \\ z\neq w \end{Sb} \left( \sum_{v \in M_k} g_{\varphi,v}(z,w) \right)
= 2\hhat_\varphi(z_n) \to 0\ .
\end{equation}
In particular, 
\begin{equation}
\label{limsupequation}
\limsup_{n\to\infty} \sum_{v \in M_k} g_{v,n} \ \leq \ 0 \ .
\end{equation}

Now let
$t_{v,N} := \inf_{z_1,\ldots,z_N \in \PP^1(\CC_v)} 
\frac{1}{N(N - 1)} \sum_{i \neq j} g_{\varphi,v}(z_i,z_j),$
so that $\Gamma_{\varphi,v} = \liminf_{N\to\infty} t_{v,N} \geq 0$ by
Lemma~\ref{DynamicalPositivityLemma}.
Note that
$g_{v,n} \geq t_{v,N_n}$ and
$\liminf_{n \to \infty} t_{v,N_n} \geq \liminf_{N \to\infty} t_{v,N},$.  
Thus
\begin{equation}
\label{liminfequation}
\liminf_{n\to\infty} g_{v,n} \ \geq \ 0
\end{equation}
for all $v$.


Finally, we apply Lemma~\ref{TechnicalLemma} to $a_n^{(v)} := g_{v,n}$.  
The hypotheses (F1) and (F2) in that lemma are satisfied
because the functions $g_{\varphi,v}$ satisfy (G1) and (G2),
and conditions ${\rm (L1)}$ and ${\rm (L2)}^{\prime}$
are satisfied because of (\ref{liminfequation}) and
(\ref{limsupequation}), respectively.
We conclude that $\lim_{n\to\infty} g_{v,n} = 0$
for each $v \in M_k$, as desired.

%
%
\end{proof}

\medskip
As a consequence of this result, and the fact that there are 
infinitely many pre-periodic points, we obtain 

\begin{cor}  \label{NormConstantCor}
Let $\varphi \in k(T)$ be a rational function of degree $d \ge 2$.  
For each place $v$ of $k$:

$A)$  The constants $\Gamma_{\varphi,v}$ 
and $I_{\mu_{\varphi,v}}(\mu_{\varphi,v})$ are equal to zero.

$B)$  The Arakelov Green's function $\frac{1}{\log(q_v)} g_{\varphi,v}(z,w)$ is normalized.
\end{cor}

\begin{proof} 
For A), choose an infinite sequence $\{x_n\}$ of distinct pre-periodic points; 
then $\hhat_{\varphi}(x_n) = 0$ for each $n$. Let $S_n$ be the set of 
Galois conjugates of $x_n$, and put $N_n = \#(S_n)$.  Applying Theorem 
\ref{PseudoEquidistributionTheorem}, we see that 
\[
\liminf \frac{1}{N_n(N_n - 1)} \sum \begin{Sb} z,w \in S_n \\
  z\neq w \end{Sb} g_{\varphi,v}(z,w) \ = \ 0 .
\]
Hence $\Gamma_{\varphi,v}$, defined in (\ref{GammaDef}), is $\le 0$.
Combined with the inequality $\Gamma_{\varphi,v} \ge 0$ 
proved in Lemma \ref{DynamicalPositivityLemma}, 
this gives $\Gamma_{\varphi,v} = 0$.  
Write $\mu = \mu_{\varphi,v}$.  
Since $D_{\infty}(\mu) = \Gamma_{\varphi,v}$ by the definitions,
it follows from Theorem \ref{DiscreteToContinuousTheorem} 
that $I_{\mu}(\mu) = 0$.   

For B), it is only necessary to show that axiom (RS3) (resp (B3))
is satisfied, i.e. we must show that  
\[
\frac{1}{\log(q_v)} \iint g_{\varphi,v}(z,w) \, d\mu(z) d\mu(w) \ = \ 0 \ .
\] 
However, this is exactly the assertion that $I_{\mu}(\mu) = 0$.  
\end{proof}
 
\subsection{The equidistribution theorem for 
                   dynamical systems on $\PP^1$.} 

{ \ \ \ }

In this subsection we will show that pseudo-equidistribution, combined
with the energy minimization principle, implies equidistribution.
 
\begin{defn}
If $S_n$ is a finite subset of $\PP^1(\CC_v)$ 
for each $n \geq 1$, we say that the sequence
$\{ S_n \}_{n \ge 1}$ is {\em equidistributed} with respect to a probability
measure $\mu$ on $\PP^1_{\Berk,v}$ over $\CC_v$ if the sequence of
measures $\delta_n = \delta(S_n)$ converges
weakly to $\mu$ on $\PP^1_{\Berk,v}$.
\end{defn}

For each $v \in M_k$, pseudo-equidistribution implies
equidistribution, in the following precise sense:

\begin{theorem}
\label{PseudoEquiImpliesEquiTheorem}
Let $\{S_n\}_{n \ge 1}$ be a sequence
of finite subsets of $\PP^1(\CC_v)$
which is pseudo-equidistributed with
respect to $g_{\varphi,v}$.
Then $\{S_n\}_{n \ge 1}$ is equidistributed with respect to $\mu_{\varphi,v}$
on $\PP^1_{\Berk,v}$.
\end{theorem}

\begin{proof}  Write $\mu = \mu_{\varphi,v}$.  
Since $\PP^1_{\Berk,v}$ is compact, it follows from
Prohorov's theorem that $\delta_n$ has a weakly convergent subsequence.
If $\nu$ is any weak limit of a subsequence of $\delta_n$,
then passing to that subsequence
\begin{eqnarray*}
0 & = & \lim_{n\to\infty} \frac{1}{N_n(N_n-1)} \frac{1}{\log(q_v)} \sum 
             \begin{Sb} z,w \in S_n \\ z\neq w \end{Sb} g_{\varphi,v}(z,w) 
             \qquad \text{by pseudo-equidistribution} \\
 & = & \lim_{n\to\infty}
     \iint_{\PP^1_{\Berk,v} \times \PP^1_{\Berk,v} \backslash {\rm (Diag)}}
           g_{\mu}(w,z)\, d\delta_n(w)d\delta_n(z) \\
& \geq & \iint_{\PP^1_{\Berk,v} \times \PP^1_{\Berk,v}} 
     g_{\mu}(w,z) \, d\nu(w)d\nu(z) \qquad \text{by Lemma \ref{LimInfLemma}} \\ 
& = & I_{\mu}(\nu) \ \ge \ I_{\mu}(\mu) \qquad 
             \text{by Theorem \ref{EnergyMinimizationTheorem}.} 
\end{eqnarray*}
Since $I_{\mu}(\mu) = 0$ by Corollary \ref{NormConstantCor},
it follows that $I_{\mu}(\nu) = I_{\mu}(\mu)$, 
so Theorem \ref{EnergyMinimizationTheorem} gives $\nu = \mu$. 
\end{proof}

\begin{remark}
For archimedean places $v$ one can give an alternative proof of 
Theorem~\ref{PseudoEquiImpliesEquiTheorem} 
using a theorem of DeMarco \cite[Theorem 1.3]{DeMarco} instead
of Theorem~\ref{EnergyMinimizationTheorem}, 
and working on $\CC^2$ rather than on $\PP^1(\CC)$.
\end{remark}

\begin{remark}
\label{PseudoEquiImpliesEquiRemark}
For archimedean $v$, the same proof shows that Theorem~\ref{PseudoEquiImpliesEquiTheorem}
remains valid if we replace $\PP^1_{\Berk,v}$ by an arbitrary compact Riemann
surface $X / \CC$ and $g_{\varphi,v}$ by any normalized Arakelov
Green's function on $X(\CC)$. 
\end{remark}

\medskip

Combining Theorem~\ref{PseudoEquidistributionTheorem} and
Theorem~\ref{PseudoEquiImpliesEquiTheorem}, we have finally proved:

\medskip

\noindent {\bf Theorem~\ref{DynamicalEquidistributionTheorem}} (Main Theorem).
{\it Let $z_n$ be a sequence of distinct points of $\PP^1(\kbar)$ such that
$\hhat_\varphi(z_n) \to 0$.  Let $S_n$ denote
the set of Galois conjugates $($over $k$$)$ of $z_n$.  
Then the sequence $\{ S_n \}_{n \ge 1}$ is equidistributed with respect
to $\mu_{\varphi,v}$ on $\PP^1_{\Berk,v} / \CC_v$ for all $v \in M_k$.}




\medskip

\section{Potential theory on Riemann surfaces}
\label{RiemannSurfaceSection}

\medskip

The goal of this section is to prove the Energy Minimization Principle
(Theorem~\ref{RiemannSurfaceTheorem}) for Arakelov Green's functions
on a compact Riemann surface.  This result was used in
\S\ref{AdelicDynamicsSection}, and is needed for the proof of our
main equidistribution theorem (Theorem~\ref{DynamicalEquidistributionTheorem}).

\medskip

\subsection{Arakelov Green's functions.}

\medskip

Let $X / \CC$ be a compact Riemann surface of genus $g$.  
Arakelov discovered that by fixing a volume form $\omega$ on $X(\CC)$,
one could define an extension of N{\'e}ron's archimedean local height
pairing from divisors of degree zero with disjoint support to
arbitrary divisors with disjoint support.  When $g \geq 1$, 
Arakelov defined a canonical volume form $\omega_{\rm can}$ 
(the pullback of the flat metric on the Jacobian
of $X$ under an Albanese embedding)
which plays a distinguished role in his theory.

The extension of N{\'e}ron's pairing arises via 
{\it Arakelov Green's functions}.  We proceed slightly more generally
than Arakelov did, using positive measures rather
than smooth $(1,1)$-forms.

Recall (Definition \ref{LogContDef}) that a 
measure $\mu$ on $X(\CC)$ is called {\it log-continuous} 
if in every coordinate patch $U \subset X(\CC)$, the function 
$\int_U \log |z-w| d\mu(w)$ is finite and continuous for all $z \in U$.
For example, any measure $\mu$ which locally has the form 
$\mu = f(z) dx \wedge dy$, where $f(z)$ is continuous and 
$dx \wedge dy$ is Lebesgue measure, is log-continuous.    


We have seen that given an Arakelov Green's function $g_{\mu}(z,w)$,
the associated measure $\mu$ is log-continuous. 
Conversely, given a log-continuous probability measure 
(i.e., a positive measure of total mass 1) $\mu$ on $X(\CC)$,
it follows from \cite[\S 2.3]{RumelyBook} or \cite[\S 4.1]{R3}
that there is a unique
pairing $\( z,w \)_\mu$ on $X(\CC)\times X(\CC)$ and a corresponding
normalized Arakelov Green's function $g_\mu(z,w) : X(\CC)\times X(\CC)
\backslash ({\rm Diag}) \to \RR$
defined by $g_\mu(z,w) = -\log \(z,w\)_\mu$ such that
axioms (RS1), (RS2) and (RS3) hold.  

\medskip

One way to prove the existence of an Arakelov Green's function $g_\mu(z,w)$ 
attached to $\mu$
is by utilizing a continuously varying family of 
{\it canonical distance functions} $[z,w]_\zeta$ on
$X(\CC)$, whose existence is proved in \cite[Theorem~2.1.1]{RumelyBook}.
It is shown in \cite[Theorem~2.3.4]{RumelyBook} that the integral
\begin{equation}
\label{CanonicalDistanceToGreen}
\int_{X(\CC)} -\log [z,w]_\zeta \, d\mu(\zeta) 
\end{equation}
satisfies properties (RS1) and (RS2) above.
Therefore $g_\mu(z,w) := \int_{X(\CC)} -\log [z,w]_\zeta d\mu(\zeta)$
is an Arakelov Green's function for $\mu$.
Furthermore, there is a unique choice of $B$ such that 
$\int_{X(\CC)} -\log [z,w]_\zeta d\mu(\zeta) + B$
satisfies condition (RS3) above and yields a normalized Arakelov
Green's function.

Conversely, given an Arakelov Green's function $-\log\( z,w \)_\mu$, 
it is shown in \cite[Theorem~2.3.3]{RumelyBook} that
one can construct a continuously varying family of canonical distance 
functions via
\begin{equation}
\label{GreenToCanonicalDistance}
[z,w]_\zeta \ := \ 
\frac{\( z,w \)_\mu}{ \( z,\zeta \)_\mu \( w,\zeta \)_\mu}.
\end{equation}
One deduces formula (\ref{GreenToCanonicalDistance}) from the relation
\begin{equation}
\label{CanonicalDistanceTransferFormula}
-\log [z,w]_\zeta = -\log [z,w]_p + \log [z,\zeta]_p + \log
 [w,\zeta]_p + C(p),
\end{equation}
which is valid for all $p,\zeta,z,w \in X(\CC)$ with $z \neq w$ (see
\cite[Corollary 2.1.5]{RumelyBook}).  Here $C(p)$ is a constant
depending only on $p$, and the right-hand side must be suitably
interpreted when $z=p$ or $w=p$.


\subsection{Examples of archimedean Arakelov Green's functions.}

\medskip

\begin{example}
\label{ArchGreenExample1}
{\bf Arakelov Green's functions on $\PP^1$}

Suppose $X = \PP^1$ and $\mu = \mu_{S^1}$ is the uniform 
probability measure on the
unit circle in $\CC = \PP^1(\CC) \backslash \{ \infty \}$.  
Then an Arakelov Green's function associated to $\mu$ is
\begin{equation}
\label{ArchGreenEquation1}
g_\mu(z,w) = \left\{ \begin{array}{ll} -\log |z-w| + \log^+|z| +
    \log^+|w| & z,w \neq \infty \\
\log^+|z| & w=\infty \\ 
\log^+|w| & z=\infty. \\ \end{array} 
\right.
\end{equation}
Note that the function $g_\mu(z,\infty) = \log^+|z|$ is the Green's
function for the unit circle in $\CC$ relative to the point at
infinity, and is also the archimedean contribution to the logarithmic 
Weil height on $\QQbar = \PP^1(\QQbar) \backslash \{ \infty \}$.

If we write (\ref{ArchGreenEquation1}) in terms of a choice of homogeneous
coordinates $z = (z_1 : z_2)$, $w = (w_1 : w_2)$, we obtain
\begin{equation}
\label{ArchGreenEquation2}
g_\mu(z,w) = -\log |z \wedge w| + \log ||z|| + \log ||w||,
\end{equation}
where $z \wedge w = z_1 w_2 - z_2 w_1$ and
$||z|| = \max \{ |z_1|,|z_2| \}$.
\end{example}

\begin{example}
\label{ArchGreenExample2}
{\bf Arakelov Green's functions on elliptic curves}

If $X=E$ is an elliptic curve over $\CC$ and $\mu = \mu_{\Haar}$ is the
normalized Haar measure on $E$, then we can take 
$g_\mu(z,w) = \lambda_\infty (z-w)$, where $\lambda_\infty$ is a
(suitably normalized) archimedean {N{\'e}ron} local height
function on $E(\CC)$ (see \cite[Section 7]{Faltings}). 
One can explicitly describe the function $\lambda_\infty$ in terms of the
Weierstrass $\sigma$-function and the quasi-period homomorphism
$\eta$ (see \cite[Chapter VI]{SilvermanII}).  

\end{example}

\subsection{Statement and discussion of Theorem~\ref{RiemannSurfaceTheorem}.}
\label{RiemannSurfaceTheoremStatementSection}

\medskip

The main result of this section is the following energy minimization
principle for Arakelov Green's functions:

\begin{theorem}
\label{RiemannSurfaceTheorem}
Let $X$ be a compact Riemann surface, let $\mu$ be a log-continuous
probability measure on $X(\CC)$, and let $g_\mu(z,w)$ be an Arakelov
Green's function for $\mu$.  Define the ``energy functional'' $I_\mu$
on the space $\PP$ of probability measures on $X(\CC)$ by the formula
\[
I_\mu(\nu) \ := \ \iint_{X(\CC) \times X(\CC)} g_\mu(z,w) \, d\nu(z) d\nu(w).
\]
Then $I_\mu(\nu) \geq I_\mu(\mu)$ for all probability measures $\nu \in \PP$,
with equality if and only if $\nu = \mu$.
\end{theorem}

In other words, $\mu$ is the unique probability measure minimizing the
energy functional $I_\mu$.  Note that by definition, $g_\mu(z,w)$ is
normalized if and only if $I_\mu(\mu) = 0$.

\medskip

The most important difference between Theorem~\ref{RiemannSurfaceTheorem}
and previous energy minimization results on Riemann surfaces (e.g. \cite[Theorem 3.1.12]{RumelyBook}) 
is that we consider the space $\PP$ of probability measures supported on all of $X(\CC)$, whereas 
in classical potential theory, one restricts attention to probability measures supported on
a compact set $E \subset X(\CC) \backslash \{ \zeta \}$ for a fixed reference point $\zeta$.
For the applications in the present paper, it is crucial to allow $\nu$ to vary over all of $\PP$,
since it is well-known that the canonical measure attached to a rational map can have
support equal to all of $\PP^1(\CC)$ (this happens, for example, with the degree 4 Latt{\`e}s maps associated to 
multiplication by 2 on an elliptic curve).  For polynomial maps, where the filled Julia set stays bounded away
from the point at infinity, one can get by with more classical results from potential theory (see \cite{BakerHsia}).

\medskip

We recall that the {\em capacity} of a compact set $F \subset \CC$ is
defined as $c(F) = e^{-V(F)}$, where $V(F)$ (the``Robin's constant''
of $F$) is the infimum (which may be a real number or $+\infty$) over
all probability measures $\nu$ supported on $F$ of expression
(\ref{capacityequation}) below.  Theorem~\ref{RiemannSurfaceTheorem} can be
viewed as a generalization
of the following fundamental result from capacity theory (see \cite{Ransford}):

\begin{theorem}
\label{FTPT}
Let $F$ be a compact subset of $\CC$ having positive 
capacity.  Then there exists a unique probability measure $\mu_F$
supported on $F$ $($called the {\it equilibrium measure} for $F$$)$ 
which minimizes the energy functional
\begin{equation}
\label{capacityequation}
I(\nu) \ = \ \iint_{F \times F} -\log |z-w|\, d\nu(z) d\nu(w).
\end{equation}
\end{theorem}

We claim that Theorem~\ref{RiemannSurfaceTheorem} implies
Theorem~\ref{FTPT} for all compact sets $F$ such that each $x \in F$
is regular for the Dirichlet problem. (For example, this holds 
if each connected component of $F$ is a continuuum;  see \cite{Ts}, 
Theorem I.11, p.7).  Such a set necessarily has positive capacity 
(\cite{Ts}, Theorem III.5, p.56). 
If $g_F(z)$ is a Green's function for $\CC \backslash F$ relative
to $\infty$, then $g_F(z)$ is continuous (\cite{Ts}, Theorem III.36, p.82).
One verifies easily that if $\mu_F = -\Delta g_F(z)$, then the function
\[
-\log |z-w| + g_F(z) + g_F(w)
\]
satisfies conditions (RS1) and (RS2) above, 
and therefore the normalized Arakelov Green's function for $\mu_F$ is
given for $z,w \in \CC, z\neq w$ by 
\[
g_{\mu_F}(z,w) \ = \ -\log |z-w| + g_F(z) + g_F(w) - C,
\]
where the constant $C$ is chosen so that (RS3) is satisfied.  (Compare
with Example~\ref{ArchGreenExample1}).
Since $\mu_F$ is supported on $F$ and $g_F \equiv 0$ on $F$, we must
in fact have
\[
C \ = \ \iint_{F\times F} -\log |z-w| \, d\mu_F(z) d\mu_F(w).
\]
Also, if $\nu$ is any probability measure supported on $F$, then
Theorem~\ref{RiemannSurfaceTheorem} implies that 
\[
I_{\mu_F}(\nu) \ = \ \iint_{F\times F} -\log |z-w| \, d\nu(z) d\nu(w) - C 
\ \geq \ 0.
\]
Therefore $\mu_F$ is the unique probability measure supported on $F$
which minimizes the energy functional $I$ (i.e., $\mu_F$ is the
equilibrium measure for $F$).  By the definition of capacity, we also
see that the constant $C$ is just the Robin's constant $V(F)$ 
of $F$, so that
\[
g_{\mu_F}(z,w) \ = \ -\log |z-w| + g_F(z) + g_F(w) + \log c(F).
\]


\begin{remark}
P.~Autissier has obtained essentially the same result as
Theorem~\ref{RiemannSurfaceTheorem} in the case where $\mu$ is the
equilibrium measure of a compact set $F \subset \CC$.
\end{remark}


\medskip

We have already discussed the relationship between 
Theorem~\ref{RiemannSurfaceTheorem} and classical potential theory in
$\CC$.  
There is also a relationship between
Theorem~\ref{RiemannSurfaceTheorem} and a result which is used in
Arakelov theory in order to establish an analogue of the Riemann-Roch theorem
for arithmetic surfaces (see \cite{LangAT}).
Indeed, when $g \ge 1$ and $\mu = \omega_{\rm can}$,
the nonnegativity of $I_\mu(\nu)$ in 
Theorem~\ref{RiemannSurfaceTheorem}
is the continuous analogue of the following theorem of Faltings 
\cite{Faltings},
originally proved using the spectral theory of the Laplacian:

\begin{theorem}[Faltings]
\label{FaltingsTheorem}
For each integer $n\geq 1$, choose an $n$-tuple $z_1,\ldots,z_n$ of
distinct points in $X(\CC)$. Then
\[
\liminf_{n\to\infty} 
\frac{1}{n(n-1)} \sum_{i\neq j} g_{\omega_{\rm can}}(z_i,z_j) \ \geq \ 0\ .
\]
\end{theorem}

It is easy to see that Theorem~\ref{RiemannSurfaceTheorem} implies
Theorem~\ref{FaltingsTheorem}.
Indeed, if $\delta_n$ is the discrete probability measure supported
equally at $z_1,\ldots,z_n$ and if $\nu$ is any weak limit of a
subsequence of $\delta_n$, then 
\[
\frac{1}{n(n-1)} \sum_{i\neq j} g_{\omega_{\rm can}}(z_i,z_j) \ = \ 
\iint_{X(\CC) \times X(\CC) \backslash {\rm (Diag)}} g_{\omega_{\rm can}}(z,w) \, 
d\delta_n(z) d\delta_n(w),
\]
and it follows from Lemma~\ref{LimInfLemma} that
\[
\liminf_{n\to\infty} 
\iint_{X(\CC) \times X(\CC) \backslash {\rm (Diag)}} g_{\omega_{\rm can}}(z,w) \,
d\delta_n(z) d\delta_n(w) \geq 
\iint_{X(\CC) \times X(\CC)} g_{\omega_{\rm can}}(z,w) \, d\nu(z) d\nu(w).
\]
This last quantity is nonnegative by
Theorem~\ref{RiemannSurfaceTheorem}, proving the claim.




\subsection{Proof of Theorem~\ref{RiemannSurfaceTheorem}.}

\medskip

We now turn to the proof of Theorem~\ref{RiemannSurfaceTheorem}.
The proof uses the representation of $g_{\mu}(z,w)$ 
in terms of the canonical distance function,
and is similar to the classical proof of 
Theorem~\ref{FTPT}.
Namely, our plan is to prove analogues of
Maria's theorem and Frostman's theorem, and to deduce 
Theorem~\ref{RiemannSurfaceTheorem} from those results.
As discussed in \S \ref{RiemannSurfaceTheoremStatementSection}, 
a key difference between Theorem~\ref{RiemannSurfaceTheorem} and
Theorem~\ref{FTPT} is the presence in Theorem~\ref{FTPT} 
of a fixed reference point at infinity.
It is the assumption that $\mu$ is log-continuous which ultimately allows us
to apply techniques from classical potential theory to the present situation.

\medskip

We fix a (possibly non-normalized) Arakelov Green's function $g_\mu(z,w)$ for
$\mu$, and for each $\zeta \in X(\CC)$ we define a canonical distance
function $[z,w]_\zeta$ by (\ref{GreenToCanonicalDistance}).

Using formula (\ref{CanonicalDistanceTransferFormula}),
we see that if $\zeta \in X(\CC)$ then for all $z,w \ne \zeta$ 
\begin{eqnarray}
g_\mu(z,w) 
& = & \int -\log [z,w]_{p} \, d\mu(p) \ \notag \\
& = & -\log [z,w]_\zeta - u_\mu(z,\zeta) - u_\mu(w,\zeta) + C_\zeta \ ,
                      \label{GMuZetaEquation}
\end{eqnarray}
where 
\[
u_\mu(z,\zeta) := \int_{X(\CC)} -\log [z,w]_\zeta \, d\mu(w)\ .
\]

By Proposition 2.1.3 of \cite{RumelyBook}, for fixed $\zeta
\in X(\CC)$, the two-variable function $-\log [z,w]_\zeta$ can be
expressed locally on $X(\CC) \times X(\CC)$ as a linear combination of
$\log |z-w|, \log |z-\zeta|$, $\log |w - \zeta|$, and a continuous
function which is harmonic in $z$ and $w$ separately.  Since $\mu$ is
log-continuous, it follows that for fixed $\zeta$, the function
$u_\mu(z,\zeta)$ is continuous on  $X(\CC) \backslash \{\zeta\}$.



\medskip

If $\nu$ is any probability measure on $X(\CC)$, we define the
generalized potential function 
$u_\mu(z,\nu) : X(\CC) \to \RR \cup \{ \infty \}$ to be
\[
u_\mu(z,\nu) \ := \ \int_{X(\CC)} g_\mu(z,w)\, d\nu(w)\ .
\]

\begin{lemma}
\label{LaplacianCalculationLemma}
For any probability measure $\nu$ on $X(\CC)$, we have
$\Delta u_\nu(z,\mu) = \nu - \mu$ as distributions.
\end{lemma}

\begin{proof}
For any test function $\psi$, we have
\begin{equation*}
\begin{aligned}
\int \psi(z) \Delta u_\nu(z,\mu) &= \int u_\nu(z,\mu) \Delta \psi(z) \\
&= \int \left( \int g_\mu(z,w)\, d\nu(w) \right) \Delta \psi(z) \\
&= \int \left( \int g_\mu(z,w)\, \Delta \psi(z) \right) d\nu(w) 
\qquad\text{(Fubini's theorem)} \notag\\
&= \int \left( \int \psi(z) \Delta g_\mu(z,w)\, \right) d\nu(w) \\
&= \int \psi(w) - \left( \int \psi(z)\, d\mu(z) \right) d\nu(w) \\
&= \int \psi(w) d\nu(w) - \int \psi(z) d\mu(z). \\
\end{aligned}
\end{equation*}
The use of Fubini's theorem at the third step is justified by the
fact that the measure $\Delta \psi(z)$  locally has
the form $f(z) \, dx \wedge dy$ for a continuous function $f(z)$.  
It follows that positive and negative parts 
$\Delta \psi^+$ and $\Delta \psi^-$ in the Jordan decomposition of 
$\Delta \psi$ are log-continuous.  Now Fubini's theorem,
in the form given in (\cite{Rudin}, Theorem 7.8(b), p.150), says that
$\int \left( \int g_\mu(z,w)\, \Delta \psi^{\pm}(z) \right) d\nu(w) 
= \int \left( \int g_\mu(z,w)\, d\nu(w) \right) \Delta \psi^{\pm}(z)$. 
\end{proof}

\begin{lemma}
\label{PotentialFunctionLemma}
For any probability measure $\nu$ on $X(\CC)$,

$A)$ $u_\nu(z,\mu)$ is a lower semicontinuous function on $X(\CC)$.

$B)$ $u_\nu(z,\mu)$ is continuous and subharmonic outside $\supp(\nu)$.

\end{lemma}

\begin{proof}
For each $M \in \RR$ define $g_\mu^{(M)}(z,w) := \min \{ M, g_\mu(z,w)
\}$, with $g_\mu^{(M)}(z,z) := M$ for all $z$.
Then $g_\mu^{(M)}(z,w)$ is a continuous function on $X(\CC) \times
X(\CC)$, and
\[
u_\nu(z,\mu) = \lim_{M\to\infty} \int_{X(\CC)} g_\mu^{(M)}(z,w) d\nu(w).
\]
Therefore $u_\nu(z,\mu)$ is lower semicontinuous, being an increasing
limit of continuous functions.  This proves part A) of the lemma.

\medskip

To prove part B), note that if $z_0 \not\in \supp(\nu)$, then by
definition there exists an open neighborhood $U$ of $z_0$ 
whose closure is disjoint from $\supp(\nu)$, such that
$\nu(U) = 0$.  For $z \in U$, we have
\[
u_\nu(z,\mu) = \int_{X(\CC)} g_\mu(z,w) d\nu(w)
= \int_{X(\CC) \backslash U} g_\mu(z,w) d\nu(w).
\]
Since $g_{\mu}(z,w)$ is uniformly continuous on 
$\overline{U} \times (X(\CC) \backslash U)$, 
it follows that $u_\nu(z,\mu)$ is continuous on $U$, 
and in particular at $z_0$.


The fact that $u_\nu(z,\mu)$ is subharmonic outside $\supp(\nu)$ follows
from the fact that $\Delta u_\nu(z,\mu) = \nu - \mu$, and thus
$- \Delta u_\nu(z,\mu) = \mu$ is a positive distribution 
on the complement of $\supp(\nu)$.  (See \cite{Kl}, Theorem 2.9.11, p.67).
\end{proof}

\medskip

The potential function $u_\nu(z,\mu)$ has additional continuity properties
besides those given by Lemma~\ref{PotentialFunctionLemma}.  For
example, we have the following (compare with \cite[Theorem~3.1.3]{Ransford}).

%

\begin{prop} 
\label{ContinuityPrincipleProp}
Let $\nu$ be a probability measure on $X(\CC)$.
Then for every $z_0 \in K := \supp(\nu)$, we have
\[
\limsup_{z \to z_0} u_\nu(z,\mu) 
\ = \ \limsup \begin{Sb} z \to z_0 \\ z \in K \end{Sb} u_\nu(z,\mu).
\]
\end{prop}


\begin{proof} 
Let $U$ be the complement in $X(\CC)$ of $K$.  
Without loss of generality, we may assume that $U \neq \emptyset$ and
that $z_0 \in \partial U$.
It is easy to see that the desired result is then equivalent to
\[
\limsup \begin{Sb} z \to z_0 \\ z\not\in K \end{Sb} u_\nu(z,\mu)
\ \leq \ 
\limsup \begin{Sb} z^{\prime} \to z_0 \\ z \in K \end{Sb} 
      u_\nu(z^{\prime},\mu)\ .
\]
If $u_\nu(z_0,\mu) = \infty$ then by lower semicontinuity we have $\lim_{z \to
  z_0} u_\nu(z,\mu) = \infty$ and the result is trivial.
Therefore we may assume that $u_\nu(z_0,\mu) < \infty$, in which case $\nu( \{
z_0 \}) = 0$.  It follows that given $\epsilon > 0$, there exists a
closed disk $D$ centered at $z_0$ such that $\nu(D) < \epsilon$.  

Let $F := D \cap K$ (so that $z_0 \in F$), and fix $z \in D \backslash F$.  

\medskip

{\bf Claim:} There exists a constant $C>0$ (independent of $z$) 
and a point $z^{\prime} \in F$ (depending on $z$) such that
\begin{equation}
\label{GreenBoundClaim}
g_\mu(z,w) \ \leq \ g_\mu(z^{\prime},w) + C
\end{equation}
for all $w \in F$.

\medskip

Before proving the claim, let's see how it gives what we want.
Integrating both sides over $F$ against $\nu$, we get
\begin{equation*}
\begin{aligned}
\label{GreenBoundEquation1}
\int_F g_\mu(z,w) d\nu(w) 
&\leq \
\int_F g_\mu(z^{\prime},w) d\nu(w) + C\cdot \nu(F) \\
&\leq \ \int_K g_\mu(z^{\prime},w) d\nu(w) -
\int_{K \backslash F} g_\mu(z^{\prime},w) d\nu(w) + C\cdot \epsilon \\
&= \ u_\nu(z^{\prime},\mu) - \int_{K\backslash F} g_\mu(z^{\prime},w) d\nu(w)
                + C \cdot \epsilon. \\
\end{aligned}
\end{equation*}
Therefore
\begin{equation}
\label{GreenBoundEquation2}
\int_F g_\mu(z,w) d\nu(w) 
+ \int_{K\backslash F} g_\mu(z^{\prime},w) d\nu(w)
\ \leq \ u_\nu(z^{\prime},\mu) + C\cdot \epsilon.
\end{equation}

As $z\to z_0$ in $D \backslash F$, we have $z^{\prime} \to z_0$ in $F$ also
(take $w = z_0$ in (\ref{GreenBoundClaim})).  By the continuity of
$g_\mu(z,w)$ on $D \times (K \backslash F)$, as $z^{\prime} \to z_0$ and $z
\to z_0$ we have
\[
\int_{K \backslash F} \left( g_\mu(z^{\prime},w) - g_\mu(z,w) \right) 
         d\nu(w) \to 0.
\]
Therefore (\ref{GreenBoundEquation2}) gives
\begin{equation*}
\label{GreenBoundEquation3}
\limsup \begin{Sb} z \to z_0 \\ z\not\in K \end{Sb} \int_K g_\mu(z,w)\, d\nu(w) 
\ \leq \ 
\limsup \begin{Sb} z^{\prime} \to z_0 \\ z \in K \end{Sb} 
u_\nu(z^{\prime},\mu) + C\cdot \epsilon.
\end{equation*}
As $\epsilon>0$ was arbitrary, this gives the desired result.

\medskip

It remains to prove the claim.  Choose $\zeta \in X(\CC)$, a small
disk $D_\zeta$ around $\zeta$ with $D_\zeta \cap D = \emptyset$, and a
constant $M_\zeta \in \RR$ so that $|g_\mu(z,\zeta)| \leq M_\zeta$ for
all $z \not\in D_\zeta$.  By \cite[proof of Theorem
3.1.6]{RumelyBook}, there is a constant $C^{\prime}$ (depending only on
$\zeta$ and $D$) such that
\[
-\log [z,w]_\zeta \ \leq \ -\log[z^{\prime},w]_\zeta + C^{\prime}
\]
for all $z,z^{\prime} \in D$ with $z \neq z^{\prime}$.

Using formula (\ref{GreenToCanonicalDistance}), we see that
\[
g_\mu(z,w) - g_\mu(z^{\prime},w) - g_\mu(z,\zeta) + g_\mu(z^{\prime},\zeta) 
\ \leq \ C^{\prime}
\]
and therefore
\[
g_\mu(z,w) - g_\mu(z^{\prime},w) \ \leq \ C^{\prime} + 2M_\zeta,
\]
which proves the claim.
\end{proof}

\medskip

We thus obtain the following analogue of Maria's theorem (see
\cite[Theorem 3.1.6]{RumelyBook}):
 
\begin{cor}
\label{MariaCor}
If $M$ is a real number such that $u_\nu(z,\mu) \leq M$ on $\supp(\nu)$,
then $u_\nu(z,\mu) \leq M$ on all of $X(\CC)$.
\end{cor}

\begin{proof}
We may clearly assume that the complement $U$ of $\supp(\nu)$ is
non-empty.  
By Lemma~\ref{PotentialFunctionLemma},
$u$ is subharmonic on $U$, and by Proposition~\ref{ContinuityPrincipleProp},
for each boundary point $z_0$ of $\supp(\nu)$ we have
\[
\limsup \begin{Sb} z \to z_0 \\ z \in U \end{Sb} u_\nu(z,\mu) = 
\limsup \begin{Sb} z \to z_0 \\ z \in \supp(\nu) \end{Sb} u_\nu(z,\mu) \leq M.
\]
The result now follows from the maximum principle for subharmonic
functions (applied to each connected component of $U$).
\end{proof}

\medskip

Next we have the following result, proved by a standard argument:  

\begin{lemma}
\label{ExistenceOfMinimizerLemma}
There exists an energy-minimizing measure $\nu_0$ for the functional $I_\mu$.
\end{lemma}

\begin{proof}
Let $V_\mu := \inf_{\nu \in \PP} I_\mu(\nu)$, and choose a sequence of
probability measures $\mu_n$ in $\PP$ such that $\lim_{n\to\infty}
I_\mu(\mu_n) = V_\mu$.  Passing to a subsequence if necessary, we may
assume that $\mu_n$ converges weakly to some measure $\nu_0$.  
We claim that
\begin{equation}
\label{ExistenceOfMinimizerEquation}
V_\mu \ = \ \liminf_{n\to\infty} I_\mu(\mu_n) \ \geq \ I_\mu(\nu_0) \ .
\end{equation}
Given this claim, we see that since
$V_\mu = \inf_{\nu \in \PP} I_\mu(\nu)$, we must have $I_\mu(\nu_0) = V_\mu$.

\medskip

To prove (\ref{ExistenceOfMinimizerEquation}), we proceed as in 
the proof of Lemma~\ref{LimInfLemma}:
{\allowdisplaybreaks
\begin{equation*}
\begin{aligned}
& \liminf_{n\to\infty}
\iint_{X(\CC) \times X(\CC)} g_\mu(z,w) \, d\mu_n(z) d\mu_n(w) \notag\\
&\ \ge    \lim_{M\to \infty}
\liminf_{n\to\infty} \iint \min\{M, g_\mu(z,w) \}\, d\mu_n(z) d\mu_n(w) 
               \qquad\text{(since $(*) \ge \min\{M, (*) \}$)} \notag\\
&\ = \lim_{M\to \infty}  \iint \min\{M, g_\mu(z,w) \}\, d\nu_0 (z) d\nu_0(w) 
              \qquad\text{($\mu_n \to \nu_0$ weakly)} \notag\\
&\ = \iint g_\mu(z,w) \, d\nu_0 (z) d\nu_0(w)
             \qquad\text{(monotone convergence theorem).} \notag\\
\end{aligned}
\end{equation*}
}
\end{proof}

\medskip

Recall from \cite[\S 3.1]{RumelyBook} that if $E \subset X(\CC)$ is
compact and $\zeta \not\in E$, then the capacity $\gamma_\zeta(E)$ of
$E$ (with respect to $\zeta$) is defined to be 
$\gamma_\zeta(E) := e^{- V_\zeta(E)}$, where
\[
V_\zeta(E) \ := \ \inf_{\nu \in \PP(E)} \iint -\log [z,w]_\zeta \,
          d\nu(z) d\nu(w) 
\]
and $\PP(E)$ is the set of probability measures supported on $E$.

Similarly, for any compact $E$ we define the $\mu$-capacity of $E$ by
$\gamma_\mu(E) := e^{- V_\mu(E)}$, where
\[
V_\mu(E) \ := \ \inf_{\nu \in \PP(E)} \iint g_\mu(z,w) \, d\nu(z) d\nu(w)\ .
\]

Note that if $E = X(\CC)$ then $I_\mu(\mu) = V_\mu(E)$.

\begin{lemma}
\label{CapacityZeroLemma}
If $E \subset X(\CC)$ is compact and $\zeta \not\in E$, 
then $\gamma_\mu(E) = 0$ if and only if $\gamma_\zeta(E) = 0$.
\end{lemma}

\begin{proof}
This follows immediately from formula (\ref{GreenToCanonicalDistance}), which
implies that for each $\nu \in \PP(E)$,  
\begin{eqnarray*} 
V_{\zeta}(E) & = & \iint -\log [z,w]_{\zeta} \, d\nu(z) d\nu(w) \\ 
& = & \iint \left( g_{\mu}(z,w) - g_{\mu}(z,\zeta) - g_{\mu}(w,\zeta) \right)
         \, d\nu(z) d\nu(w) \\
& = & V_\mu(E) - 2\int_E g_\mu(z,\zeta) \, d\nu(z).
\end{eqnarray*}            
Here $\int_E g_\mu(z,\zeta) \, d\nu(z) < \infty$, 
since $\supp(\nu) \subseteq E$ and $\zeta \not\in E$.
\end{proof}

\begin{remark}
  If $\zeta, \zeta^{\prime} \not\in E$ then it follows from Lemma 
  \ref{CapacityZeroLemma} that $\gamma_\zeta(E) = 0$ if and only if
  $\gamma_{\zeta^{\prime}}(E) = 0$.  In particular, it makes sense to speak of
  a set of capacity zero on $X(\CC)$ without reference to a particular
  base point $\zeta$.  
\end{remark}

\begin{lemma}
\label{CapacityMeasureLemma}
Let $\nu$ be a probability measure on $X(\CC)$ such that $I_\mu(\nu) <
\infty$.  Then $\nu(A) = 0$ for any Borel subset $A \subset X(\CC)$ of
capacity zero.
\end{lemma}

\begin{proof}
Suppose to the contrary that $\nu(A) > 0$.  Then for some
compact subset $A^{\prime} \subseteq A$ we also have $\nu(A^{\prime}) > 0$,
so without loss of generality we may assume that $A$ itself is compact.
Recall that there exists $M \in \RR$ such that $g_\mu(z,w) \geq -M$
for all $z,w \in X(\CC)$.  Then 
\[
\int_{A} \int_{A} (g_\mu(z,w)+M) \, d\nu(z) d\nu(w)
\ \le \ \int_{X(\CC)} \int_{X(\CC)} (g_\mu(z,w)+M) \, d\nu(z) d\nu(w),
\]
so that
\begin{eqnarray*}
\int_A \int_A g_\mu(z,w) \, d\nu(z) d\nu(w) 
& \leq & \int_{X(\CC)} \int_{X(\CC)} g_\mu(z,w) \, d\nu(z) d\nu(w) 
                 + M \cdot (1 - \nu(A)^2) \\
& \le &  I_\mu(\nu) + M \ < \ \infty.
\end{eqnarray*}
Define a probability measure $\nu^{\prime}$ on $A$ by setting $\nu^{\prime} :=
\frac{1}{\nu(A)} \nu |_{A}$.  Then 
\[
I_\mu(\nu^{\prime}) \ \leq \ 
\frac{1}{\nu(A)^2} \left( M + I_\mu(\nu) \right) \ < \ \infty \ , 
\]
so that $\gamma_\mu(A) > 0$, a contradiction.
\end{proof}

\medskip

The following result is an analogue of Frostman's theorem (see \cite[Theorem
3.1.7]{RumelyBook}):

\begin{theorem}
\label{FrostmanTheorem}
Let $\nu_0$ be any probability measure which minimizes the functional
$I_\mu$, and let $V := I_\mu(\nu_0)$.  Then the potential function
$u(z) := u_{\nu_0}(z,\mu)$ on $X(\CC)$ satisfies:

A) $u(z) = V$ for all $z \in X(\CC)$ outside a set of capacity zero.

B) $u(z) \leq V$ for all $z \in X(\CC)$.

\end{theorem}

\begin{proof}
We first show that $u(z) \geq V$ for all $z \in X(\CC)$ outside a set
of capacity zero.  For each $n\geq 1$, put  
$A_n := \{ z \in X(\CC) \, : \, u(z) \leq V - \frac{1}{n} \}$.  
The lower semicontinuity of $u_{\nu_0}(z)$ 
shows that each $A_n$ is closed, and we have
$\cup A_n = A := \{ z \in X(\CC) \, : \, u(z) < V \}$.  
Clearly $A \neq X(\CC)$, since $\int u(z) d\nu_0(z) = V$ and $\nu_0$
is a positive measure.  

If $\zeta$ is any point in the complement of $A$, we claim that
$\gamma_\zeta(A) = 0$.  By Lemma~\ref{CapacityZeroLemma}, it suffices
to show that $\gamma_\mu(A) = 0$.  Suppose on the contrary that
$\gamma_\mu(A) > 0$.  To obtain a contradiction, we first construct
disjoint closed subsets $E_1, E_2$ of $X(\CC)$ as follows.

By \cite[Proposition~3.1.5]{RumelyBook}, we must have
$\gamma_\zeta(A_n) > 0$, and hence $\gamma_\mu(A_n) > 0$, for some
$n$.  Thus for a suitable $\epsilon > 0$ and $n \geq 1$, we have $u(z) < V -
2\epsilon$ on $E_1 := A_n$ and $\gamma_\mu(E_1) > 0$. 

As $\int u(z) d\nu_0(z) = V$, there exists $z_0 \in \supp(\nu_0)$
such that $u_\nu(z_0,\mu) > V - \epsilon$.  Lower semicontinuity implies that
this inequality remains valid in a closed disk $D$ around $z_0$, which
we may assume to be disjoint from $E_1$.  Since $z_0 \in
\supp(\nu_0)$, we have $\nu_0(D) > 0$.
Let $E_2 := D$, and let $W := \nu_0(E_2) > 0$.

Since $\gamma_\mu(E_1) > 0$, there exists a probability measure $\nu^{\prime}$
supported on $E_1$ such that $I_\mu(\nu^{\prime}) < \infty$.  Define a new
measure $\sigma$ on $X(\CC)$ by setting
\[
\sigma := \left \{
\begin{array}{ll}
W \nu^{\prime}   &   {\rm on \;} E_1    \\
- \nu_0 &   {\rm on \;} E_2    \\
0       &   {\rm elsewhere \;} \\
\end{array}
\right.
\]

Then $\sigma(E_1) = W$, $\sigma(E_2) = -W$, and $\sigma(X(\CC)) = 0$.  
Note that for each real number $t \in [0,1]$, $\nu_0 + t\sigma$ is a
probability measure on $X(\CC)$.
As in \cite[proof of Theorem 3.1.6]{RumelyBook}, we calculate that
$I_\mu(\sigma) < \infty$ and 
\begin{equation}
\label{FrostmanEquation}
I_\mu(\nu_0 + t\sigma) - I_\mu(\nu_0) \leq (-2W\epsilon)\cdot t +
I_\mu(\sigma)\cdot t^2.
\end{equation}
For $t$ sufficiently small, the right-hand side of 
(\ref{FrostmanEquation}) is negative, contradicting the fact that
$I_\mu(\nu_0) = V$ is the minimum possible energy of a probability
measure on $X(\CC)$.
This contradiction proves that $A$ has capacity zero, and by
construction we have
$u(z) \geq V$ for all $z \in X(\CC) \backslash A$.

\medskip

Next, we show that $u(z) \leq V$ on $\supp(\nu_0)$.
To see this,
suppose for the sake of
contradiction that $u_\nu(z_0,\mu) > V$ for some $z_0 \in \supp(\nu_0)$.  By
the lower semicontinuity of $u$, there exists $\epsilon > 0$ and a
closed disk $D$ around $z_0$ such that $u(z) > V + \epsilon$ on $D$.

As $z_0 \in \supp(\nu_0)$, the number $T := \nu_0(D)$ is
positive.

We have already seen that $u(z) \geq V$ for all $z \in X(\CC)$,
except on a set $A$ of capacity (and hence, by
Lemma~\ref{CapacityMeasureLemma}, of $\nu_0$-measure) zero.
Therefore
\[
V = \int u(z) d\nu_0(z) \geq V(1-T) + (V + \epsilon)T > V,
\]
a contradiction.  It follows that $u(z) \leq V$ on $\supp(\nu_0)$ as claimed.

\smallskip

Finally, Corollary~\ref{MariaCor} now shows that $u(z) \leq V$ on all of
$X(\CC)$, which proves both A) and B).
\end{proof}

\medskip

We can now prove Theorem~\ref{RiemannSurfaceTheorem}.
\medskip

\begin{proof}
  As in the statement of Theorem~\ref{FrostmanTheorem}, let $\nu_0$ be
  any probability measure which minimizes the functional $I_\mu(\nu)$, and
  let $V = I_\mu(\nu_0)$.  Since $I_\mu(\mu) < \infty$ by (\ref{eq:LF'}),
  we must have $V < \infty$ as well.

\medskip

Integrating over $X(\CC) \times X(\CC)$ and applying Fubini's theorem, we have
\begin{equation} \label{FubiniEquation1}
\int \left( \int g_\mu(z,w) \, d\mu(z) \right) \, d\nu_0(w)
\ = \ \int \left( \int g_\mu(z,w) d\nu_0(w) \right) d\mu(z).
\end{equation}
The interchange of order of integration is justified by the
same reasons as in Lemma \ref{LaplacianCalculationLemma}, 
because $\mu$ is log-continuous.   

Here the left side equals $I_{\mu}(\mu)$, since
$\int g_\mu(z,w) d\mu(z)$ is constant by property $\textrm{(RS3)}^{\prime}$ of
Arakelov Green's functions.  

On the other hand, Theorem~\ref{FrostmanTheorem} shows that 
$u(z) := \int g_\mu(z,w) \, d\nu_0(w) = V$
outside a set $A$ with capacity zero.  
Since $I_\mu(\mu) < \infty$ by assumption,
Lemma~\ref{CapacityMeasureLemma} shows that $\mu(A) = 0$.  Hence the
right side of (\ref{FubiniEquation1}) is $V$.  

Combining these gives  $I_\mu(\mu) = V$.
Therefore $\mu$ is also energy-minimizing.  

To see that $\mu = \nu_0$, note first that sets of capacity zero
have Lebesgue measure zero in any coordinate patch on $X(\CC)$.
(This follows from Lemma \ref{CapacityMeasureLemma}).  
Hence $u(z) = V$ almost everywhere with respect to Lebesgue
measure, and consequently $\Delta u(z) = 0$.
On the other hand, Lemma~\ref{LaplacianCalculationLemma} gives us
the distributional identity $\Delta u(z) = \nu_0 - \mu$.  
Therefore $\mu = \nu_0$ as desired.
\end{proof}

\medskip

\section{Comparison and calculation of various capacities}
\label{HomogeneousSectionalCapacitySection}

As before, we let $\varphi : \PP^1 \rightarrow \PP^1$ be a rational map of
degree $d \ge 2$ defined over a number field $k$, 
and let $F = (F_1,F_2) : \AA^2 \rightarrow \AA^2$  be a lifting of $\varphi$,
where $F_1(x,y)$ and $F_2(x,y)$ are homogeneous polynomials of degree
$d$ with coefficients in $k$ having no common factors over $\kbar$.  

In this section we will prove Theorem \ref{AdelicDeMarcoCapacityTheorem},
the resultant formula for the homogeneous
transfinite diameter of the filled Julia set  $K_{F,v}$ : 
\begin{equation}
d^0_\infty(K_{F,v}) \ = \ |\Res(F)|_v^{-1/d(d-1)} \ .
\end{equation}
We do this by considering various notions of capacity:
the local and global sectional capacities and  
the Chebyshev constant studied in \cite{RL} and \cite{RLV}, 
and the homogenous sectional capacity and 
homogeneous transfinite diameter, 
which are introduced here for the first time.  
The reason for this proliferation of capacities is that
we can compute the sectional capacity, and there
are standard methods for proving inequalities between various other types of 
capacities.  In outline, the plan is to first prove
\begin{equation*}
\text{local sectional capacity} \ = \ |\Res(F)|_v^{-1/d(d-1)} 
\end{equation*} 
by proving an upper bound for the local sectional capacity and 
using the fact that the global sectional capacity is the 
product of the local sectional capacities, and then to show that
for circled sets, 
\begin{eqnarray*}
 \text{local sectional capacity} & = & \text{homogeneous sectional capacity} \\
 & = & \text{homogenous transfinite diameter.} 
\end{eqnarray*}  

We prove more in this section than is strictly needed
for our application to dynamics.  It is our hope that the ideas
developed here will be useful in other contexts as well.  In particular,
it would be interesting to know if the resultant formula for the 
local sectional capacity of the pullback of a ball
(Proposition \ref{Formula1}) generalizes to higher dimensions.   

\subsection{Sectional capacities of polynomial domains.}

\medskip

We can view $F$ as defining a finite map $\tF : \PP^2 \rightarrow \PP^2$
given in homogeneous coordinates by $\tF(X:Y:Z) = (F_1(X,Y):F_2(X,Y):Z^d)$.
Its action on the affine patch $\AA^2$ is given by $F$
and it stabilizes hyperplane $H = \{Z=0\}$, which we identify with $\PP^1$,
where its action is given by $\varphi$.
The map $\tF$ has degree $d^2$, and $\tF^*(H) = d \cdot H$.

\vskip .1 in

The definition of the sectional capacity for sets in $\PP^2$, relative to the
divisor $H$, is as follows.

For each place $v$ of $k$, let $E_v \subset \PP^2(\CC_v)$ be a nonempty
set which is stable under the group of continuous automorphisms
$\Gal^c(\CC_v/k_v) \cong \Gal(\tk_v/k_v)$ and is bounded away from
$H(\CC_v)$ under the $v$-adic metric on $\PP^2(\CC_v)$.  For all but finitely
many $v$ we assume that $E_v = B(0,1) \times B(0,1) \subset \AA^2(\CC_v)$,
the `trivial set' for $v$ with respect to $H$.  We will call these assumptions 
the {\emph{Standard Hypotheses}}.

Put $\EE = \prod_v E_v \subset \AA_k$, where $\AA_k$ is the adele ring
of $k$.

For each $n \ge 0$, identify
the space of sections $\Gamma(n) := H^0(\PP^2,\cO_{\PP^2}(n))$ with the set of
homogeneous polynomials in $k[X,Y,Z]$ of degree $n$.
Consider the basis for $k[X,Y,Z]$ given by the monomials
$\{X^k Y^{\ell} Z^m\}$ ;  equip it with the term order $\prec$ given by the
lexicographic order with $Z \prec X \prec Y$, graded by the degree.
We call this structure the `monic basis';  it is the key ingredient used
in defining local sectional capacities.
(Any other term order graded by the degree would work;  this one is 
most directly compatible with dehomogenization.)

For each place $v$ of $k$, let $\vol_v$ be additive Haar measure 
on $k_v$ (normalized so that $\vol_v(\O_v) = 1$ if $v$ is nonarchimedean, 
and given by Lebesgue measure on $\RR$ or $\CC$ if $v$ is archimedean).  
Let $\vol_{\AA}$ be the additive
Haar measure on the adele ring $\AA_k$ given by the product of the measures $\vol_v$.  
For each $n$, by transport of structure using the monic basis we obtain 
Haar measures $\vol_v$ on the vector spaces
$\Gamma_{\varphi,v}(n) = k_v \otimes_k \Gamma(n)$
and $\vol_{\AA}$ on the $\AA_k$-module
$\Gamma_{\AA}(n) = \AA_k \otimes_k \Gamma(n)$.

To define norms, we dehomogenize at $Z$,
writing $x = X/Z$, $y = Y/Z$,
and identify $\Gamma(n)$ with the space of polynomials in
$k[x,y]$ of total degree $\le n$.  We view these as functions on $\AA^2$.
Put 
\begin{eqnarray*}
\cF_v(n) & = & \{ f \in \Gamma_{\varphi,v}(n) : \|f\|_{E_v} \le 1 \} \ , \\
\cF_{\AA}(n) & = & \left( \prod_v \cF_v(n) \right) \cap \Gamma_{\AA}(n) \ .
\end{eqnarray*}
The local sectional capacity $\Sg(E_v,H)$ is defined by 
\begin{equation*}
-\log(\Sg(E_v,H)) \ = \ 
       \lim_{n \rightarrow \infty} \frac{3!}{n^3} \log(\vol_v(\cF_v(n))) 
\end{equation*}
and the global sectional capacity $\Sg(\EE,H)$ by 
\begin{equation*}
-\log(\Sg(\EE,H)) \ = \ 
  \lim_{n \rightarrow \infty} \frac{3!}{n^3} \log(\vol_{\AA}(\cF_{\AA}(n))) \ .
\end{equation*} 
In \cite{RL} it is shown that under the Standard Hypotheses, 
the limits defining $\Sg(E_v,H)$ and $\Sg(\EE,H)$ exist, and that 
\begin{equation*}
\Sg(\EE,H) \ = \ \prod_v \Sg(E_v,H) \ .
\end{equation*}

\vskip .1 in

We now apply this to polydiscs in $\CC_v^2$ and their pullbacks by $F$.
Given $z = (x,y) \in \CC_v^2$, write $\|z\|_v = \max(|x|_v,|y|_v)$.
(For archimedean $v$, this is a different definition of $\|z\|_v$
than we used in \S\ref{AdelicDynamicsSection}.)

For each $R_v > 0$, put
\begin{equation*}
B_v(R_v) \ := \ B(0,R_v)^2 \ = \
\{ z \in \CC_v^2 : \max(|x|_v,|y|_v) \le R_v \} \ .
\end{equation*}
Thus
\begin{equation*}
F^{-1}(B_v(R_v))
\ = \ \{z \in \CC_v : \max(|F_1(z)|_v,|F_2(z)|_v) \le R_v \} \ .
\end{equation*}
Given a collection of numbers $\vec{R} = \{R_v\}$ with $R_v = 1$ for all but
finitely many $v$, define the adelic sets
\begin{eqnarray*}
\BB(\vec{R}) & = & \prod_v B_v(R_v) \ , \\
F^{-1}(\BB(\vec{R})) & = & \prod_v F^{-1}(B_v(R_v)) \ .
\end{eqnarray*}

\begin{prop} \label{P1} { \ \ \ }

$A)$ For each $v$,  the local sectional capacity $\Sg(B_v(R_v),H)$
equals $R_v^2$.

$B)$ The global sectional capacity $\Sg(\BB(\vec{R}),H)$ equals $\prod_v R_v^2$.

$C)$ The global sectional capacity
$\Sg(F^{-1}(\BB(\vec{R})),H)$ equals $(\prod_v R_v^2)^{1/d}$.
\end{prop}

\begin{proof}  Part A) follows 
from the fact that the  
logarithmic capacity of a ball 
in $\PP^1$ is $\gamma_{\infty}(B(0,R_v)) = R_v$ 
(see, e.g. \cite{RumelyBook}, Example 5.2.15, p.352),
together with the formula for the the sectional capacity of a product set   
$E_v = E_{v,1} \times E_{v,2} \subset \CC_v^2$:
\begin{equation} \label{PFor}
\Sg(E_v,H) = \gamma_{\infty}(E_{v,1}) \cdot \gamma_{\infty}(E_{v,2})
\end{equation}
(see \cite{RL}, Example 4.3, p.558).

Part B) follows from part A) and (\cite{RL}, Theorem 3.1, p.552).

Part C) follows from part B) and functorial properties of the 
global sectional capacity.  By the the pullback formula for finite maps 
(\cite{RLV}, Theorem 10.1, p.54), we have  
\begin{equation*}
\Sg(F^{-1}(\BB(\vec{R})),dH) \ = \ \Sg(\BB(\vec{R}),H)^{d^2} \ ,
\end{equation*}
since $\tF^{-1}(\BB(\vec{R})) = F^{-1}(\BB(\vec{R}))$, 
$\tF^*(H) = dH$, and $\deg(\tF) = d^2$.  
Also, by the homogeneity of the sectional capacity in its
second variable (\cite{RLV}, Theorem C (5), p.9),   
$\Sg(\EE,dH) = \Sg(\EE,H)^{d^3}$ for any $\EE$ in $\PP^2$.
Combining these gives C).
\end{proof}

\vskip .1 in
Determining the local sectional capacity $\Sg(F^{-1}(B_v(R_v)),H)$
is more difficult.
As before, let $\Res(F)$ denote the resultant of $F_1$ and $F_2$.

\begin{prop} \label{Formula1}  For each $v$,
\begin{equation*}
\Sg(F^{-1}(B_v(R_v)),H) \ = \ (R_v^2)^{1/d} \cdot |\Res(F)|_v^{-1/d^2} \ .
\end{equation*}
\end{prop}

Before giving the proof, we will need a lemma.
For each $m$, write  $\Gamma_v^0(m)$ for the
space of homogeneous polynomials in $k_v[x,y]$ of degree $m$.

Take $m = td + d-1$ and
consider the collection of $m+1 = (t+1)d$ polynomials
\begin{equation*}
\{x^i y^j F_1(x,y)^k F_2(x,y)^{\ell} : i + j = d-1, \ k + \ell = t \}
\ \subset \  \Gamma_v^0(m) \ .
\end{equation*}
Let $\Det(m)$ denote the determinant of the matrix expressing these
polynomials in terms of the standard monomials 
$\{ x^m, x^{m-1}y, \ldots, y^m \}$.

\begin{lemma} \label{DetLemma}  For $m = td + d-1$,
$\Det(m) = \pm \Res(F)^{t(t+1)/2}$ \ .
\end{lemma} 
  
\begin{proof}  We will first show that $\Det(m)$ vanishes if and only if
$\Res(F)$ vanishes.  Indeed, $\Det(m) = 0$ if and only if there is
a nontrivial relation of the form  
\begin{equation} 
     \sum_{i=0}^t h_i(x,y) F_1(x,y)^{t-i} F_2(x,y)^{i} \ = \ 0  \label{F1}
\end{equation}
where each $h_i(x,y)$ is homogeneous of degree $d-1$.  

If (\ref{F1}) holds, let $I$ be the least
index for which $h_I(x,y) \ne 0$;  necessarily $I < t$.      
Then $F_2(x,y)$ divides $h_I(x,y) F_1(x,y)^{t-I}$.
Since $F_2(x,y)$ has degree $d$,
it must have an irreducible factor in common with $F_1(x,y)$, 
so $\Res(F) = 0$.  Conversely, if $\Res(F) = 0$
then there is a nontrivial relation 
\begin{equation*}
      h_0(x,y) F_1(x,y) + h_1(x,y) F_2(x,y) \ = 0               \label{F2}
\end{equation*}
where $h_0$ and $h_1$ are homogeneous of degree $d-1$.  Multiplying
through by $F_1(x,y)^{t-1}$ gives a relation of the form (\ref{F2}).

Expand $\Det(m)$ and $\Res(F)$ as polynomials in the coefficients of
$F_1$ and $F_2$.   Comparing degrees and using the fact that $\Res(F)$ is 
irreducible, we see that 
\begin{equation*}
\Det(m) \ = \ C \cdot \Res(F)^{t(t+1)/2}
\end{equation*}
for some constant $C$.  Taking $F_1 = x^d$, $F_2 = y^d$ and evaluating
both sides, we find that $C = 1$ for an appropriate ordering of the terms. 
\end{proof}                 

\vskip .1 in

We now turn to the proof of Proposition \ref{Formula1}.

\begin{proof} 

We give the proof only when $R_v = 1$.  The general case reduces to this
by a scaling argument.

It suffices to prove the upper bound
$\Sg(F^{-1}(B_v(1)),H) \le |\Res(F)|_v^{-1/d^2}$ for each $v$.
If this is known, then by the global equalities
\begin{eqnarray*}
\prod_v \Sg(F^{-1}(B_v(1)),H) \ = \ 1 \ , \\
\prod_v |\Res(F)|_v \ =  1 \ , \phantom{xxx}
\end{eqnarray*}
the local inequality must actually be an equality, for each $v$.

\vskip .1 in
Write $E_v = F^{-1}(B_v(1))$.  Since $E_v$ is bounded, there is a constant
$c_v \in \CC_v$ such that $\|c_v x^i y^j\|_{E_v} \le 1$
for all $i$, $j$ with $i + j \le 2d-1$.

We will study $\vol_v(\cF_v(n))$ by making use of the decomposition
$\Gamma_v(n) = \oplus_{m=0}^n \Gamma_v^0(m)$, which is compatible
with the monomial bases.
For each $m \ge 0$, put
\begin{equation*}
\cF_v^0(m) \ = \ \{f \in \Gamma_v^0(m) : \|f\|_{E_v} \le 1 \} \ .
\end{equation*}
If $m \ge 2d-1$, we can uniquely write
$m = t d + (d-1) + r$ with integers $t \ge 1$, $0 \le r \le d-1$,
and then
\begin{equation*}
\Gamma_v^0(m) \ = \ y^r \cdot \Gamma_v^0(m-r)
      \oplus (\bigoplus_{i=0}^{r-1} k_v \cdot x^{m-i} y^i)
\end{equation*}      

By Lemma \ref{DetLemma},
the polynomials $x^i y^j F_1^k F_2^{\ell}$ with $i+j = d-1$, $k + \ell = t$
form a basis for $\Gamma_v^0(m-r)$, so the corresponding polynomials
$x^i y^{j+r} F_1^k F_2^{\ell}$, together with the monomials
$x^{m-i} y^i$ for $0 \le i < r$, form a basis for $\Gamma_v^0(m)$.
Again by Lemma \ref{DetLemma},
the transition matrix from the monomial basis for $\Gamma_v^0(m)$
to this new basis has determinant $\pm \Res(F)^{t(t+1)/2}$.

For each basis element of the first type, we have
$i + j +r \le 2d-1$, so
\begin{equation*}
c_v \cdot x^i y^{j+r} F_1^k F_2^{\ell} \ \in \ \cF_v^0(m) \ .
\end{equation*}
For each basis element of the second type, 
$x^{m-i} y^i = x^{td} \cdot x^{d-1+r-i} y^i$ with
$(d-1+r-i)+i \le 2d-1$, so
\begin{equation*}
c_v^{t+1} \cdot x^{m-i} y^i \ \in \ \cF_v^0(m) \ .
\end{equation*}

\vskip .1 in
Now suppose $v$ is nonarchimedean. By the ultrametric inequality, we have
\begin{equation} \label{FCC0}
\big( \bigoplus \begin{Sb} i + j = d-1 \\ k + \ell = t \end{Sb}
       \cO_v \cdot c_v x^i y^{j+r} F_1^k F_2^{\ell} \big)
       \oplus
    \big( \bigoplus_{0 \le i < r} \cO_v \cdot c_v^{t+1} x^{m-i} y^i \big)
       \ \subset \ \cF_v^0(m)
\end{equation}
and it follows that
\begin{equation*}
\vol_v(\cF_v^0(m)) \ \ge \
(|c_v|_v)^{(m-r+1)+r(t+1)} \cdot |\Res(F)|_v^{t(t+1)/2} \ .
\end{equation*}
Noting that $m/d > t > (m/d)-1$, we see that
\begin{equation}  \label{FCC1}
\log(\vol_v(\cF_v^0(m)) \ \ge \
     \frac{m^2}{2d^2} \log(|\Res(F)|_v)
                  - O(m).
\end{equation}
By increasing the implied constant, we can assume this
holds for $m \le 2d-1$ as well.

For each $n \ge 0$, the ultrametric inequality shows that
\begin{equation} \label{FCC2}
\bigoplus_{m=0}^n \cF_v^0(m) \ \subset \ \cF_v(n) \ .
\end{equation} 
Using (\ref{FCC1}) and (\ref{FCC2}), it follows that
\begin{eqnarray}
-\log(\Sg(E_v,H)) & = & \lim_{n \rightarrow \infty}
                    \frac{3!}{n^3} \log(\vol_v(\cF_v(n)))  \notag \\
  & \ge & \lim_{n \rightarrow \infty}
     \frac{3!}{n^3} \sum_{m=0}^n \big( \frac{m^2}{2} \cdot \frac{1}{d^2}
                \log(|\Res(F)|_v) - O(m) \big) \notag \\
  & = & \frac{1}{d^2} \log(|\Res(F)|_v) \ .  \label{FCD}
\end{eqnarray}
Thus $\Sg(E_v,H) \le |\Res(F)|_v^{-1/d^2}$, as desired.

\vskip .1 in
If $v$ is archimedean and $k_v \cong \RR$, 
the triangle inequality gives a weaker containment;  it is 
better to formulate the result directly for $\cF_v(n)$.    
Noting that $\dim_k(\Gamma(n)) = (n+1)(n+2)/2$, we obtain  
\begin{eqnarray*}
\lefteqn{ \big( \bigoplus_{i+j < 2d-1} {\scriptscriptstyle \frac{2}{(n+1)(n+2)} } 
        \cdot [-1,1] \cdot c_v x^i y^j \big) \oplus} \\
&  & \phantom{xxx}     
\bigoplus_{m= 2d-1}^{n} \Big[ \big( \bigoplus 
    \begin{Sb} m = td + d-1 + r \\ k + \ell = t \\ i + j = d-1  \end{Sb}
 {\scriptscriptstyle \frac{2}{(n+1)(n+2)} } 
       \cdot [-1,1] \cdot c_v x^i y^{j+r} F_1^k F_2^{\ell}  \big) 
     \\ & &  \phantom{xxxxxxx} \oplus 
       \big( \bigoplus_{0 \le i < r}  {\scriptscriptstyle \frac{2}{(n+1)(n+2)} }
            \cdot [-1,1] \cdot c_v^{t+1} x^{m-i} y^i \big) \Big]
       \ \subset \ \cF_v(n) \ .
\end{eqnarray*}

From this inclusion, we deduce that 
\begin{equation} \label{FCC5}
\log(\vol_v(\cF_v(n)) \ \ge \
     \frac{n^3}{6d^2} \log(|\Res(F)|_v) - O(n^2 \log(n)) \ .
\end{equation}
If $v$ is archimedean and $k_v \cong \CC$ there is a similar containment,  
with $B(0,1)$ replacing $[-1,1]$.  In either case,
a computation like the one in (\ref{FCD}) now shows that
$\Sg(E_v,H) \le |\Res(F)|_v^{-1/d^2}$.
\end{proof}

\medskip

\subsection{The sectional capacity of the filled Julia set.}

\medskip


Let the numbers $r_v$, $R_v$ be as in 
Corollary \ref{Cor1}.  Put
$K_F = \prod_v K_{F,v}$,
$\BB(\vec{r}) =  \prod_v B_v(r_v)$, $\BB(\vec{R})  =  \prod_v B_v(R_v)$.
Then for each $n$, 
\begin{eqnarray} 
(F^{(n)})^{-1}(B_v(r_v)) \ \subseteq \ K_{F,v}
 \ \subseteq \ (F^{(n)})^{-1}(B_v(R_v)) \ ,  \label{FA3A} \\
(F^{(n)})^{-1}(\BB(\vec{r})) \ \subseteq \ K_F
 \ \subseteq \ (F^{(n)})^{-1}(\BB(\vec{R})) \quad \ . \label{FA3}
\end{eqnarray}

\begin{theorem} \label{Thm1} { \ \ \ \ }

$A)$ The global sectional capacity $\Sg(K_F,H) = 1$.  Equivalently,
\begin{equation*}
\prod_v \Sg(K_{F,v},H) \ = \ 1 \ .
\end{equation*}

$B)$ For each $v$, \ $\Sg(K_{F,v},H) = |\Res(F)|_v^{-1/(d(d-1))}$.
\end{theorem}

\begin{proof}
Part (A) follows from Proposition \ref{P1} (C), using (\ref{FA3})
and the fact that $F^{(n)}$ is homogeneous of degree $d^n$.

Part (B) comes out as follows.    
By \cite[Corollary~6.4]{DeMarco}, 
\begin{equation*}
\Res(F^{(n)}) \ = \ \Res(F)^{(d^{2n-1}-d^{n-1})/(d-1)} \ .
\end{equation*}
Hence Proposition \ref{Formula1}, applied to both halves of (\ref{FA3A}), gives 
\begin{eqnarray*}
\Sg(K_{F,v},H) 
       & = & \lim_{n \rightarrow \infty} |\Res(F^{(n)})|_v^{-1/(d^n)^2} \\
       & = & \lim_{n \rightarrow \infty} 
               |\Res(F)|_v^{ -(d^{2n-1}-d^{n-1})/d^{2n}(d-1) } \\
       & = & |\Res(F)|_v^{-1/d(d-1)} \ . 
\end{eqnarray*}      
\end{proof}


\subsection{\bf The homogeneous sectional capacity.}

Identify $\Gamma(n)$ with the space of polynomials 
in $k[x,y]$ of total degree $\le n$.  Instead of considering the asymptotics 
for volumes related to $\Gamma(n)$ in the definition of the sectional capacity, 
one can can consider the corresponding asymptotics for volumes concerning 
homogeneous polynomials alone.  This gives rise to the  
homogeneous sectional capacity.  

Put
\begin{equation*}   
\Gamma^0(n) \ = \ 
\{ f(x,y) \in k[x,y] : \text{$f$ is homogeneous of degree $n$} \}
\end{equation*}
and write $\Gamma_v^0(n) = k_v \otimes_k \Gamma^0(n)$, 
$\Gamma_{\AA}^0(n) = \AA_k \otimes_k \Gamma^0(n)$.  Put
\begin{eqnarray*}
\cF_v^0(n) & = & \{f \in \Gamma_v^0(n) : \|f\|_{E_v} \le 1\} \ , \\
\cF_{\AA}^0(n) & = & \prod_v \cF_v^0(n) \cap \Gamma_{\AA}^0(n) \ .
\end{eqnarray*}
Equipping $\Gamma_v^0(n)$ and $\Gamma_{\AA}^0(n)$ with the 
bases given by the monomials $x^k y^{\ell}$,  and giving those bases
the term order $\prec$ given by the lexicographic order with $x \prec y$,
graded by the degree, we have a situation analogous to that in 
the definition of the sectional capacity.  By transport of structure, the
Haar measure $\vol_v$ on $k_v$ induces a  measure $\vol_v$ on
each $\Gamma_v^0(n)$, and the Haar measure $\vol_{\AA} = \prod_v \vol_v$ on
$\AA_k$ induces a measure on each $\Gamma_{\AA}^0(n)$.  Define the local
homogeneous sectional capacity $\Sg^0(E_v,H)$ by
\begin{equation}  \label{LocalHSC}
-\log(\Sg^0(E_v,H)) \ = \ \lim_{n \rightarrow \infty}
       \frac{2!}{n^2} \log(\vol_v(\cF_v^0(n))) 
\end{equation}
and the global homogeneous sectional capacity $\Sg^0(\EE,H)$ by
\begin{equation} \label{GlobalHSC}
-\log(\Sg^0(E_v,H)) \ = \ \lim_{n \rightarrow \infty}
       \frac{2!}{n^2} \log(\vol_{\AA}(\cF_{\AA}^0(n))) \ . 
\end{equation} 
The existence of these limits follows from
the general existence theorem for sectional capacities of line bundles 
with ``adelically normed sections'' (see \cite{RLV}, Theorem A, p.4), 
applied to the variety $\PP^1$ rather than $\PP^2$.  
The details are as follows.  

One can interpret the set of homogeneous polynomials $\Gamma^0(n)$ 
as the space of global sections $H^0(\PP^1,\cO_{\PP^1}(n))$.  The sup norms 
$\|f\|_{E_v}$ on the spaces $\Gamma_v^0(n)$ satisfy axioms (A1) and (A2) 
of (\cite{RLV}, p.3), because they are implied by Standard Hypotheses for sets 
(\cite{RLV}, Example 1.1, p.13).  Then (\cite{RLV}, Theorem 6.2, p.66) asserts
that the limit (\ref{LocalHSC}) exists, and (\cite{RLV}, Theorem 7.1, p.73) 
tells us that the limit (\ref{GlobalHSC}) exists, and also that
\begin{equation*}
\Sg^0(\EE,H) \ = \ \prod_v \Sg^0(E_v,H)  \ .
\end{equation*}

For each place $v$ (archimedean or nonarchimedean),
given a set $E_v \subset \CC_v^2$, the {\emph{circled}}
set obtained from $E_v$ is
\begin{equation*}
E_v^0 \ = \ \{ w z : z \in E_v, \ w \in \CC_v, |w|_v = 1\} \ .
\end{equation*}
We will call $E_v$ {\emph{circled}} if $E_v = E_v^0$.  Trivially
$\Sg^0(E_v^0,H) = \Sg^0(E_v,H)$.   

\begin{prop} \label{PCircled}
For each bounded, Galois-stable set $E_v \subset \CC_v^2$,
we have $\Sg(E_v,H) \le \Sg^0(E_v,H)$.
If $E_v = E_v^0$, then $\Sg(E_v,H) = \Sg^0(E_v,H)$.
\end{prop}

\begin{proof}

We first show that $\Sg(E_v,H) \le \Sg^0(E_v,H)$.
The proof breaks into two cases, according as 
$\Sg^0(E_v,H) > 0$ or $\Sg^0(E_v,H) = 0$.  

First suppose $\Sg^0(E_v,H) > 0$, and fix $\varepsilon > 0$.  
By the definition of $\Sg^0(E_v,H)$, for all sufficiently large $k$
\begin{equation} \label{FS1} 
-\log(\Sg^0(E_v,H))-\varepsilon \ \le \ \frac{2!}{k^2} \log(\vol_v(F_v^0(n)))
    \ \le \ -\log(\Sg^0(E_v,H))+\varepsilon \ .
\end{equation}

Suppose $v$ is nonarchimedean.  By the ultrametric inequality, for
each $n$ we have
\begin{equation*}
\bigoplus_{k=0}^n \cF_v^0(k) \ \subset \ \cF_v(n) \ .
\end{equation*}
Since the union of the monomial bases for the $\Gamma_v^0(k)$ 
is the monomial basis for $\Gamma_v(n)$, we have 
$\vol_v(\cF_v(n)) \ge \prod_{k=0}^n \vol_v(\cF_v(n))$.
Hence
\begin{eqnarray}
-\log(\Sg(E_v,H))
   & = & \lim_{n \rightarrow \infty} \frac{3!}{n^3} \log(\vol_v(\cF_v(n))
                           \notag  \\
   & \ge & \lim_{n \rightarrow \infty}
            \frac{3!}{n^3} \sum_{k=0}^n \log(\vol_v(\cF_v(k))  \notag \\
   & \ge & \lim_{n \rightarrow \infty}
            \frac{3!}{n^3} \sum_{k=0}^n \frac{k^2}{2!}
                       (-\log(\Sg^0(E_v,H))-\varepsilon) \label{FC1} \\
   & = & -\log(\Sg^0(E_v,H))-\varepsilon \ . \notag
\end{eqnarray}
Since $\varepsilon > 0$ is arbitrary, 
$\Sg(E_v,H) \le \Sg^0(E_v,H)$.

If $v$ is archimedean, then 
$1/(n+1) \oplus_{k=0}^n \cF_v^0(k) \ \subset \ \cF_v(n)$   
by the triangle inequality, so
\begin{equation*}
\frac{3!}{n^3} \log \vol_v(\cF_v(n)
\ \ge \ 
\frac{3!}{n^3} \sum_{k=0}^n \log \vol_v(\cF_v(k) 
       -  \frac{3!}{n^3} \cdot \frac{n(n+1)}{2} \log |n+1|_v \ . 
\end{equation*}
Thus the factor of $1/(n+1)$ washes out in 
the asymptotics as $n \rightarrow \infty$,  
and the same computation as in (\ref{FC1}) carries through.

If $\Sg^0(E_v,H) = 0$, take $M > 0$;  then  
for all sufficiently large $k$, 
\begin{equation} \label{FS2} 
\frac{2!}{k^2} \log(\vol_v(F_v^0(n))) \ \ge \ M\ .
\end{equation}
If $v$ is nonarchimedean, using 
(\ref{FS2}) in place of (\ref{FS1}) in (\ref{FC1}) shows 
$-\log(\Sg(E_v,H)) \ge M$.  Since $M$ is arbitrary, $\Sg(E_v,H) = 0$.
A similar argument, with minor modifications to deal with the 
triangle inequality, applies in the archimedean case.

\vskip .1 in
Now suppose $E_v = E_v^0$.
We claim that $\Sg^0(E_v,H) \le \Sg(E_v,H)$.  

We will give the proof under the assumption that $\Sg^0(E_v,H) > 0$, leaving 
the other case to the reader.  

Fix $n$, and suppose $f \in \cF_v(n)$.  Thus $f$ is a polynomial in $k_v[x,y]$
of total degree at most $n$, with $\|f\|_{E_v} \le 1$.
Decompose $f = \sum_k f_k$ as
the sum of its homogeneous parts of degree $k$.  The $f_k$ can be recovered
from $f$ by finite Fourier analysis:  if $\zeta$ is a primitive $(n+1)$-st
root of unity, then
\begin{equation*}
f_k(z) \ = \ \frac{1}{n+1} \sum_{\ell=0}^n \zeta^{-\ell k} f(\zeta^{\ell} z) \ .
\end{equation*}
Since $E_v$ is circled and $\|f\|_{E_v} \le 1$, for each $\ell$ we have
$\|f(\zeta^{\ell} z)\|_v \le 1$.

If $v$ is archimedean, then $\|f_k\|_{E_v} \le 1$, and
$\cF_v(n) \subset \oplus_{k=0}^n \cF_v^0(k)$.
If $v$ is nonarchimedean, then 
$(n+1) \cF_v(n) \subset \oplus_{k=0}^n \cF_v^0(k)$.
The factor $n+1$ has no effect on the asymptotics of the volumes;
in both the archimedean and nonarchimedean case, for each $\varepsilon
> 0$ we have
\begin{eqnarray*}
-\log(\Sg(E_v,H))
   & = & \lim_{n \rightarrow \infty} \frac{3!}{n^3} \log(\vol_v(\cF_v(n)) \\
   & \le & \lim_{n \rightarrow \infty}
            \frac{3!}{n^3} \sum_{k=0}^n \frac{k^2}{2!}
                       (-\log(\Sg^0(E_v,H))+\varepsilon) \\
   & = & -\log(\Sg^0(E_v,H))+\varepsilon \ ,
\end{eqnarray*}
which yields $\Sg^0(E_v,H) \le \Sg(E_v,H)$.
\end{proof}

\begin{remark}
It can happen that $\Sg(E_v,H) < \Sg^0(E_v,H)$.  
For example, let $k = \QQ$ and let $v$ be the archimedean place.  
Take $E_v = [0,1] \times B(0,1) \subset \CC^2$;  
then $E_v^0 = B(0,1) \times B(0,1)$.  It is well known from the
classical theory of logarithmic capacities that 
$\gamma_{\infty}([0,1]) = 1/4$ and $\gamma_{\infty}(B(0,1)) = 1$.
By formula (\ref{PFor}),
\begin{equation*}
\Sg(E_v,H) = 1/4 \ , \ \qquad \Sg(E_v^0,H) = 1 \ .
\end{equation*}
But then Proposition \ref{PCircled} and the remarks preceding it imply that
$\Sg^0(E_v,H) = \Sg^0(E_v^0,H) = \Sg(E_v^0,H) = 1$.  
Similar examples can be given for nonarchimedean $v$.
\end{remark}

\subsection{The homogeneous Chebyshev constant.}

As before, let $E_v \subset \CC_v^2$ be bounded and
stable under $\Gal^c(\CC_v/k_v)$.  
It will be useful to introduce another quantity equal to the
homogeneous sectional capacity:  the homogeneous Chebyshev constant.  
Its chief virtue is that it is independent of the ground field
used to compute it.

For each field $L$ with $k_v \subset L \subset \CC_v$, put  
$\Gamma_L^0(n) = L \otimes_{k_v} \Gamma_v^0(n)$.  
(The cases of greatest interest are $L = k_v$ and $L = \CC_v$.) 
The basis elements for $\Gamma^0(n)$ are 
$x^n \prec x^{n-1}y \prec \cdots \prec y^n$. 
For each $k = 0, \ldots, n$, define the set of 
``monic $L$-rational homogeneous polynomials of degree $n$ 
with leading term $x^{n-k} y^k$" to be 
\begin{equation*}
\Gamma_L^0(n,k) \ := \ \{ f \in L \otimes_k \Gamma^0(n) : 
f(x,y) = x^{n-k} y^k + \sum_{i < k} a_i x^{n-i} y^i \} 
\end{equation*}
and put 
\begin{equation*}
M_L(n,k) \ = \ \inf_{f \in \Gamma_L^0(n,k)} (\|f\|_{E_v})^{1/n} \ .
\end{equation*}
By (\cite{RLV}, Theorem 6.1, p.64), the number 
$\CH_L^0(E_v,H)$ defined by 
\begin{equation*} 
\log(CH_L^0(E_v,H)) \ = \ \lim_{n \rightarrow \infty} \frac{2!}{n}
\left( \sum_{k=0}^n \log(M_L(n,k)) \right)
\end{equation*}
exists, and is independent of $L$.  We will call it the homogeneous
Chebyshev constant.  By (\cite{RLV}, Theorem 6.2, p.66),
\begin{equation} \label{For3}
\Sg^0(E_v,H) \ = \ \CH_{k_v}^0(E_v,H) \ = \ \CH_{\CC_v} ^0(E_v,H) \ .
\end{equation} 

\subsection{The homogeneous transfinite diameter.}

\medskip
  
We have already introduced the homogeneous transfinite diameter $d_{\infty}^0(E_v)$.
The definition makes sense for an arbitrary set $E_v \subset \CC_v^2$.
Under the assumptions that $E_v$ is bounded and stable
under $\Gal^c(\CC_v/k_v)$, we will show that it coincides with the
homogeneous sectional capacity.

Note that since the kernel $|z_i \wedge z_j|_v$ is invariant when $z_i$ or $z_j$ is
multiplied by $w \in \CC_v$ with $|w|_v = 1$, we clearly have
$d_{\infty}^0(E_v) = d_{\infty}^0(E_v^0)$.

\begin{prop} \label{PropH2}
Let $E_v$ be bounded and stable under $\Gal^c(\CC_v/k_v)$.  Then
\begin{equation*}
\Sg^0(E_v,H) \ = \ d_{\infty}^0(E_v) \ .
\end{equation*}
\end{prop}

\begin{proof}
We claim that we can assume without loss that each $d_n(E_v) > 0$ 
and that each $\cF_v(n) \subset \Gamma_v^0(n)$ is bounded.  

If $d_n(E_v) = 0$ for some $n$, then there would be a finite set of points 
$\xi_1, \ldots, \xi_n$ with 
\begin{equation*}
E_v \ \subset \ \{ w \xi_i : 1 \le i \le n, \ w \in \CC_v\} \ .
\end{equation*}   
This means $d_m^0(E_v) = 0$ for each $m > n$, so $d_{\infty}^0(E_v) = 0$.
Also, for each $m > n$ the polynomial 
$f_m(z) = (\prod_{i=1}^n (z \wedge \xi_i)) \wedge (z \wedge \xi_1)^{m-n} 
\in \Gamma_m^0(n)$ vanishes on $E_v$, which means that 
$\vol_v(\cF_v(m)) = \infty$ so $\Sg^0(E_v,H) = 0$.  

On the other hand, if $\cF_v(n)$ is unbounded for some $n$, then 
(by the local compactness of $k_v$) there would 
be a nonzero polynomial $f(z) \in \Gamma_v^0(n)$ with $\|f\|_{E_v} = 0$.  
Factoring $f(z) = \prod_{i=1}^n (z \wedge \xi_i)$, we find that again 
\begin{equation*}
E_v \ \subset \ \{ w \xi_i : 1 \le i \le n, \ w \in \CC_v\} \ ,
\end{equation*}  
so $d_{\infty}^0(E_v) = \Sg^0(E_v,H) = 0$.  

Henceforth we will assume that each $d_n^0(E_v) > 0$ and each $\cF_v(n)$ is
bounded.

\vskip .1 in 
We will first show that $d_{\infty}^0(E_v) \le \Sg^0(E_v)$.

Fix $n$.  Given $\varepsilon > 0$, choose
$\xi_1, \ldots, \xi_{n+1} \in E_v$ so that
\begin{equation*}
|P_{n+1}(\xi_1, \ldots, \xi_{n+1})|_v
\ \ge \ (d_{n+1}^0(E_v)-\varepsilon)^{(n+1)n/2} \ .
\end{equation*}
Write $\xi_{\ell} = (x_{\ell},y_{\ell})$
and let $D$ be the $(n+1) \times (n+1)$ matrix
whose $\ell$-th column is given by the $x_{\ell}^{n-i} y_{\ell}^i$, 
$i = 0, \ldots, n$.
Then $\det(D) = \pm \prod_{i < j} (\xi_i \wedge \xi_j)$.  (To see this, first
suppose each $x_{\ell}$ is nonzero.  If $x_{\ell}^{n}$ is factored out 
from each column, we obtain a Vandermonde determinant in the 
$y_{\ell}/x_{\ell}$.  When the $x_{\ell}^{n}$ are multiplied through 
in the formula for the Vandermonde,
the formula for $\det(D)$ results.  
The general case follows by continuity.)
Thus $|\det(D)|_v \ge (d_{n+1}^0(E_v)-\varepsilon)^{(n+1)n/2}$.

If $v$ is nonarchimedean, the bounded $\cO_v$-module
$\cF_v(n)$ has an $\cO_v$-basis $\{g_0, \ldots, g_n\}$.
Write $g_{\ell}(z) = \sum_{i=0}^n c_{\ell,i} x^i y^{n-i}$ 
for $i = 0, \ldots, n$, 
and let $C$ be the $(n+1) \times (n+1)$ matrix with rows $c_{\ell,i}$,
$0\le i \le n$.  Then $\vol_v(\cF_v(n)) = |\det(C)|_v$.

For each $g \in \cF_v(n)$ we have $|g(\xi_i)|_v \le 1$.  Identifying 
a polynomial with its vector of coefficients, this says that 
$\cF_v(n) \subset D^{-1} \hcO_v^{n+1}$.  Passing to the
the $\hcO_v$-module $\hcO_v \otimes_{\cO_v} \cF_v(n) = C \hcO_v^{n+1}$
generated by $\cF_v(n)$, we conclude
that $C \hcO_v^{n+1} \subset D^{-1} \hcO_v^{n+1}$, so
$|\det(C)|_v \le |\det(D^{-1})|_v$.  Thus 
\begin{equation*}
\vol_v(\cF_v(n)) \ \le \ (d_{n+1}^0(E_v)-\varepsilon)^{-(n+1)n/2} \ ,
\end{equation*}
and hence  
\begin{eqnarray}
-\log(\Sg^0(E_v,H))
& = & \lim_{n \rightarrow \infty} \frac{2!}{n^2} \log(\vol_v(\cF_v(n))) 
                       \label{FM1} \\
& \le & \lim_{n \rightarrow \infty} -\frac{n+1}{n} 
                  \log(d_{n+1}^0(E_v)- \varepsilon) 
\ = \ -\log(d_{\infty}^0(E_v)-\varepsilon) \notag \ .
\end{eqnarray}
Letting $\varepsilon \rightarrow 0$,
we obtain $d_{\infty}^0(E_v) \le \Sg^0(E_v,H)$.

If $v$ is archimedean, then by the same argument as above we find
\begin{equation*}
\cF_v(n) \ \subset \ D^{-1} B(0,1)^{n+1} \ .
\end{equation*}
If $k_v \cong \CC$, this gives
$\vol_v(\cF_v(n)) \le |\det(D)|_v^{-1} \cdot \pi^{n+1}$ 
(normalizing $| \ |_v$ 
so that $|w|_v = |w|^2$ for $w \in \CC$).  
If $k_v \cong \RR$, then the triangle inequality shows that
\begin{equation*}
\frac{1}{2}(\cF_v(n) \oplus i \cdot \cF_v(n))
        \ \subset \ D^{-1} B(0,1)^{n+1} \ ,
\end{equation*}
which gives 
$\vol_v(\cF_v(n)) \le |\det(D)|_v^{-1} \cdot (2 \sqrt{\pi})^{n+1}$ \ . 
In either case a computation similar to that in (\ref{FM1}) shows that 
$d_{\infty}^0(E_v) \le \Sg^0(E_v,H)$.  

\vskip .1 in
Next we will show that $\Sg^0(E_v,H) \le d_{\infty}^0(E_v)$.

Fixing $n$, we can choose $\xi_0, \ldots, \xi_n \in E_v$ such that
\begin{equation} \label{FE4}
\frac{1}{2} d_{n+1}^0(E_v)^{(n+1)n/2} \ \le \ |P_{n+1}(\xi_0, \ldots, \xi_n)|_v
\ \le \ d_{n+1}^0(E_v)^{(n+1)n/2} \ .
\end{equation}
For each $\ell = 0, \ldots, n$,  put 
\begin{equation} \label{FE5}
f_{\ell}(z) \ = \ \frac{\prod_{i \ne \ell} (z \wedge \xi_i)}
                       {\prod_{i \ne \ell} (\xi_{\ell} \wedge \xi_i)}
    \ = \ \pm \frac{P_{n+1}(\xi_0, \ldots, \xi_{\ell-1},z,
                           \xi_{\ell+1}, \ldots, \xi_n)}
                   {P_{n+1}(\xi_1, \ldots, \xi_{n+1})} \ .
\end{equation}
Writing $z = (x,y)$, we can also expand $f_{\ell}(z)$ as a polynomial
\begin{equation*}
f_{\ell}(z) \ = \ \sum_{i=0}^n a_{\ell,i} x^{n-i} y^i \ \in \CC_v[x,y] \ .
\end{equation*}
By the definition of $d_{n+1}^0(E_v)$, 
\begin{equation}
\sup_{z \in E_v} \label{FE5A}
|P_{n+1}(\xi_0, \ldots, \xi_{\ell-1},z,\xi_{\ell+1}, \ldots, \xi_n)|_v
\ \le \ d_{n+1}^0(E_v)^{(n+1)n/2} \ .
\end{equation}
Combining (\ref{FE4}), (\ref{FE5}), 
and (\ref{FE5A}) shows that $\|f_{\ell}\|_{E_v} \ \le \ 2$.

If $v$ is nonarchimedean, let 
$\cG_v(n) \subset \CC_v \otimes_{k_v} \Gamma_v^0(n)$ 
be the $\hcO_v$-module generated by
$f_0(z), \ldots, f_n(z)$.  
Let $A$ be the $(n+1) \times (n+1)$ matrix with entries $a_{\ell,i}$, and 
let $D$ be the $(n+1) \times (n+1)$ matrix defined above, 
with columns obtained from $\xi_1, \ldots, \xi_{n+1}$.
Since $f_{\ell}(\xi_i) = \delta_{\ell,i}$ it follows that  $A \cdot D = I$. 
Hence 
\begin{equation} \label{FE6}
|\det(A)|_v \ = \ |\det(D)|_v^{-1} \ \ge \ d_{n+1}(E_v)^{-(n+1)n/2} \ .
\end{equation}
The nonsingularity of $A$ means that
the $f_{\ell}(z)$ span $\CC_v \otimes \Gamma_v^0(n)$.
Using elementary row operations, we can transform $A$ to a lower triangular
matrix $B$ with $\det(B) = \det(A)$.  Equivalently, we can find a new
$\hcO_v$-basis $\{g_0(z), \ldots, g_n(z)\}$ for $\cG_v(n)$ having the form
\begin{equation*}
g_{\ell}(z) \ = \ \sum_{i=0}^{\ell} b_{\ell,i} x^{n-i} y^i \ .
\end{equation*}
By the ultrametric inequality, $\|g_{\ell}\|_{E_v} \le 2$ for each $\ell$.
Dividing through by $b_{\ell,\ell}$ we obtain a monic homogeneous
polynomial
\begin{equation*}
h_{\ell}(z) \ = \ x^{n-\ell} y^{\ell} +
    \sum_{i=0}^{\ell-1} \frac{b_{\ell,i}}{b_{\ell,\ell}} x^{n-i} y^i
\end{equation*}
with sup norm $\|h_{\ell}\|_{E_v} \le 2 |b_{\ell,\ell}|_v^{-1}$.
From the definition of $M_{\CC_v}(n,\ell)$, we obtain
\begin{equation} \label{FE7}
M_{\CC_v}(n,\ell)^n \ \le \ 2 |b_{\ell,\ell}|_v^{-1} \ .
\end{equation}
Using (\ref{FE7}) and (\ref{FE6}), together with $|\det(B)|_v = |\det(A)|_v$,
\begin{eqnarray*}
\frac{2!}{n} \sum_{\ell=0}^n \log_v(M_{\CC_v}(n,\ell))
   & \le & \frac{2}{n^2} \sum_{\ell = 0}^n \log(2/|b_{\ell,\ell}|_v) \\
   & = & \frac{2}{n^2} \log(2) - \frac{2}{n^2} \log(|\det(B)|_v) \\
   & \le & \frac{2}{n^2} \log(2) + \frac{n+1}{n} \log(d_{n+1}^0(E_v)) \ .
\end{eqnarray*}
Passing to the limit as $n \rightarrow \infty$, we get
\begin{equation*}
\CH_{\CC_v}^0(E_v) \ \le \ d_{\infty}^0(E_v) \ .
\end{equation*}
Since $\Sg^0(E_v,H) = \CH_{\CC_v}^0(E_v)$, we are done.

If $v$ is archimedean, fix $n$ and let the $f_{\ell}(z) \in \CC[x,y]$ 
be as above;  each $f_{\ell}(z)$ satisfies $\|f_{\ell}\|_{E_v} \le 2$.
If $k_v \cong \CC$, then by the triangle inequality
\begin{equation*}
\frac{1}{n+1} \cdot\frac{1}{2} 
\cdot \bigoplus_{\ell=1}^n B(0,1) f_{\ell} \ \subset \ \cF_v(n) \ .
\end{equation*}
Introducing the matrices $A$ and $D$ as before, 
and identifying a polynomial with its vector of coefficients, we have
$1/(2(n+1)) \cdot A \cdot B(0,1)^{n+1} \subset \cF_v(n)$. 
If $k_v \cong \RR$, then 
\begin{equation*}
\frac{1}{n+1} \cdot \frac{1}{2} 
\cdot \bigoplus_{\ell=1}^n B(0,1) f_{\ell} \ \subset \ 
     \cF_v(n) \oplus i \cF_v(n) \ \subset \ \CC \otimes_{k_v} \Gamma_v^0(n)\ 
\end{equation*}
and $1/(2(n+1)) \cdot A \cdot B(0,1)^{n+1} 
\subset \cF_v(n) \oplus i \cF_v(n)$. 
 
In either case,
\begin{eqnarray*}
\vol_v(\cF_v(n)) 
  & \ge & |2(n+1)|_v^{-(n+1)} \cdot |\det(A)|_v \cdot \pi^{n+1} \\
  & \ge & |2(n+1)|_v^{-(n+1)} \cdot \pi^{n+1} 
               \cdot d_{n+1}^0(E_v)^{-(n+1)n/2} \ .
\end{eqnarray*}  
Thus 
\begin{eqnarray*}
-\log(\Sg^0(E_v,H)) & = & 
        \lim_{n \rightarrow \infty} \frac{2!}{n^2} \log(\vol_v(\cF_v(n))) \\
        & \ge & \lim_{n \rightarrow \infty} \frac{n+1}{n} 
                   (-\log(d_{n+1}^0(E_v))) \ = \ -\log(d_{\infty}^0(E_v)),
\end{eqnarray*}
or equivalently, $\Sg^0(E_v,H) \le d_{\infty}^0(E_v)$.
\end{proof}

\vskip .1 in
As we have noted before, the homogeneous transfinite diameter is defined for 
arbitrary sets $E_v \subset \CC_v^2$.  When $F = (F_1,F_2)$ is defined over 
a number field, then Propositions \ref{PCircled}, \ref{PropH2},
and Theorem \ref{Thm1} tell us that the filled Julia set $K_{F,v}$ satisfies
\begin{equation} \label{FCum} 
\Sg(K_{F,v},H) \ = \ \Sg^0(K_{F,v},H) \ = \ d_{\infty}^0(K_{F,v}) 
               \ = \ |\Res(F)|_v^{-1/d(d-1)} \ .
\end{equation}               

In particular, we have finally proved
Theorem~\ref{AdelicDeMarcoCapacityTheorem}.

\begin{remark}
Note that Corollary~\ref{CapacityProductFormula} follows
directly from Propositions \ref{PCircled}, \ref{PropH2}
and Theorem \ref{Thm1}(A), and in particular
Proposition~\ref{Formula1} is not needed for the proof.
\end{remark}

Finally, we show that the last equality in (\ref{FCum}) holds 
for arbitrary polynomial maps defined over $\CC_v$ (and not just for maps
defined over $k$).

\begin{cor} \label{Formula2}
Let $F_1, F_2 \in \CC_v[x,y]$ be homogeneous polynomials of degree
$d \ge 2$, having no common factor.  Put $F = (F_1,F_2)$ and let $R_v > 0$.  
Then

$A)$   $d_{\infty}^0(F^{-1}(B_v(R_v))) 
             = (R_v^2)^{1/d} \cdot |\Res(F)|_v^{-1/d^2}$ \ , 

$B)$  $d_{\infty}^0(K_{F,v}) \ = \ |\Res(F)|_v^{-1/(d(d-1))}$ \ .
\end{cor}               

\begin{proof}  Part (B) follows from (A) by 
 the proof of Theorem \ref{Thm1}, so it suffices to prove (A).  

First suppose $v$ is nonarchimedean.  Since $\Qbar$ is dense in $\CC_v$,
there are polynomials $\tF_1, \tF_2$ defined over a number field $k$ 
such that the map $\tF = (\tF_1,\tF_2)$ satisfies 
$\tF^{-1}(B_v(R_v)) = F^{-1}(B_v(R_v))$, 
and such that $|\Res(\tF)|_v = |\Res(F)|_v$.  
Let $w$ be the place of $k$
induced by the given embedding $k \hookrightarrow \CC_v$.  The normalized
absolute value $| \ |_w$ on $\CC_w \cong \CC_v$ is a power of $| \ |_v$,
say $| \ |_w = | \ |_v^D$.  For each set $E \subset \CC_w \cong \CC_v$, 
the homogeneous transfinite diameters $d_{\infty}^0(E)_w$ and 
$d_{\infty}^0(E)_v$ computed relative to $| \ |_w$ and $| \ |_v$, 
are related by the same power $D$.  Hence our assertion follows from 
(\ref{FCum}) for $\tF$.  

If $v$ is archimedean, choose a sequence 
$\{s_1, s_2, \ldots \}$ decreasing monotonically to $R_v$ from above,
and another sequence $\{r_1, r_2, \ldots \}$ increasing monotonically to 
$R_v$ from below.  By continuity, for each $n$ we can
choose $\Qbar$-rational polynomials 
$\tF_{n,1}, \tF_{n,2}$ close enough to $F_1$ and $F_2$ that the maps
$\tF_n = (\tF_{n,1},\tF_{n_2})$ satisfy
$\tF_n^{-1}(B_v(r_n)) \subset F^{-1}(B_v(R_v)) \subset \tF_n^{-1}(B_v(s_n))$
and  
\begin{equation*}
\lim_{n \rightarrow \infty} \Res(\tF_n) \ = \ \Res(F) \ .  
\end{equation*}
If $K_v \cong \RR$, 
we also require that the coefficients of the $\tF_{n,i}$ belong to $\RR$.  
Applying (\ref{FCum}) to 
the sets $\tF_n^{-1}(B_v(r_n))$, $\tF_n^{-1}(B_v(s_n))$ and taking a limit,
we obtain (A).  
\end{proof}


\end{document}